\documentclass[10pt]{amsart}
\usepackage[dvipdfmx]{graphicx}
\usepackage{amsmath,amscd,amsthm,amsxtra,yhmath,stackrel}
\usepackage{epsfig,graphics,color,colortbl}
\usepackage{amssymb,latexsym}
\usepackage{mathrsfs}
\usepackage{bm}
\usepackage[poly,all]{xy}
\usepackage{hyperref}
\usepackage{tikz-cd}
\usepackage[draft]{todonotes}
\usepackage{upgreek}
\usepackage{tensor}
\usepackage{ytableau}

\newtheorem{thm}{\bf Theorem}[section]
\newtheorem{df}[thm]{\bf Definition}
\newtheorem{prop}[thm]{\bf Proposition}
\newtheorem{cor}[thm]{\bf Corollary}
\newtheorem{lem}[thm]{\bf Lemma}
\newtheorem{rem}[thm]{\bf Remark}
\newtheorem{ex}[thm]{\bf Example}

\numberwithin{equation}{section}

\newcommand{\mc}{\mathcal}
\newcommand{\mf}{\mathfrak}
\newcommand{\ms}{\mathscr}
\newcommand{\mb}{\bm}

\newcommand{\pf}{\noindent{\bfseries Proof. }}
\newcommand{\ov}{\overline}
\newcommand{\un}{\underline}

\newcommand{\U}{{\mc U}}
\newcommand{\V}{\mc{V}}
\newcommand{\W}{\mc{W}}
\newcommand{\WW}{\ms{W}}
\newcommand{\cP}{\mathscr{P}}

\newcommand{\I}{\mathbb{I}}
\newcommand{\be}{{\bf e}}

\newcommand{\Z}{\mathbb{Z}}
\newcommand{\Q}{\mathbb{Q}}
\newcommand{\C}{\mathbb{C}}

\newcommand{\e}{\epsilon}
\newcommand{\de}{\delta}
\newcommand{\ude}{\updelta}

\newcommand{\td}{\widetilde}

\newcommand{\gl}{\mf{gl}}

\newcommand{\La}{\Lambda}

\newcommand{\la}{\lambda}

\newcommand{\hf}{\frac{1}{2}}

\newcommand{\ot}{\otimes}

\newcommand{\si}{(-1)^{\e_i}}

\newcommand{\bq}{{\bf q}}

\usepackage{stmaryrd}
\usepackage{verbatim}
\usepackage{mathtools}
\usepackage[title]{appendix}
\usepackage{stackengine}
\def\dring#1{\stackinset{c}{}{t}{-2.5pt}%
  {$
  \mkern2mu\scriptscriptstyle{
  {\scalebox{.4}{\begin{picture}(1,1)
	\put(0,0){\makebox(0,0)[c]{$\circ$}}
	\end{picture}}}
  }
  \mkern-6mu{
  \scalebox{.4}{\begin{picture}(1,1)
	\put(0,0){\makebox(0,0)[c]{$\circ$}}
	\end{picture}}
  }
  $}{$#1$}}

\theoremstyle{definition}

\theoremstyle{remark}

\newcommand{\ket}[1]{\left|#1\right\rangle}

\global\long\def\ket#1{\left|#1\right\rangle }%
\global\long\def\boldm{\mathbf{m}}%
\global\long\def\boldmp{\mathbf{m}^{\prime}}%
\global\long\def\bolde{\mathbf{e}}%

\newcommand{\ep}{\epsilon}

\begin{document}
\title[Representations of quantum affine orthosymplectic superalgebras]
{Oscillator representations of quantum affine orthosymplectic superalgebras}

\author{JAE-HOON KWON}

\address{Department of Mathematical Sciences and RIM, Seoul National University, Seoul 08826, Korea}
\email{jaehoonkw@snu.ac.kr}

\author{SIN-MYUNG LEE}

\address{Department of Mathematical Sciences, Seoul National University, Seoul 08826, Korea}
\email{luckydark@snu.ac.kr}

\author{MASATO OKADO}

\address{Osaka Central Advanced Mathematical Institute \&
Department of Mathematics, Osaka Metropolitan University,
3-3-138 Sugimoto, Sumiyoshi-ku, Osaka, 558-8585, Japan}
\email{okado@omu.ac.jp}

\thanks{J.-H.Kwon and S.-M.Lee are supported by the National Research Foundation of Korea(NRF) grant funded by the Korea government(MSIT) (No.2020R1A5A1016126). S.-M.Lee is supported in part by a KIAS Individual Grant (MG095001) at Korea Institute for Advanced Study. M.Okado is supported by JSPS KAKENHI Grant Number~JP19K03426. This work was partly supported by Osaka Central Advanced Mathematical Institute (MEXT Joint Usage/Research Center on Mathematics and Theoretical Physics JPMXP0619217849).}

\begin{abstract}
We introduce a category of $q$-oscillator representations over the quantum affine superalgebras of type $D$ and construct a new family of its irreducible representations. Motivated by the theory of super duality, we show that these irreducible representations naturally interpolate the irreducible $q$-oscillator representations of type $X_n^{(1)}$ and the finite-dimensional irreducible representations of type $Y_n^{(1)}$ for $(X,Y)=(C,D),(D,C)$ under exact monoidal functors. This can be viewed as a quantum (untwisted) affine analogue of the correspondence between irreducible oscillator and irreducible finite-dimensional representations of classical Lie algebras arising from Howe's reductive dual pairs $(\mf{g},G)$, where $\mf{g}=\mf{sp}_{2n}, \mf{so}_{2n}$ and $G=O_\ell, Sp_{2\ell}$.
\end{abstract}

\maketitle
\setcounter{tocdepth}{1}
\section{Introduction}

\subsection{}
Howe duality is a principle which gives a representation theoretic approach to the first fundamental theorem of the classical invariant theory \cite{H,H2}, and has played important roles in representation theory with various applications and generalizations.
It is in general given by a Fock space $\mc{F}$ with an action of a reductive dual pair $(\mf{g},G)$ of a Lie (super)algebra $\mf{g}$ and a classical group $G$, which forms the centralizers of each other and hence provides a one-to-one correspondence between the irreducible representations of $\mf{g}$ and $G$ in $\mc{F}$. 

If there exist two dualities $(\mf{g},G)$ on $\mc{F}$ and $(\mf{g}',G)$ on $\mc{F}'$ for the same $G$, then this yields a correspondence between irreducible representations of $\mf{g}$ and $\mf{g}'$ corresponding to the same irreducible representation of $G$, and explains a similarity between two semisimple tensor categories generated by them (see \cite{CW} and references therein for examples of various dual pairs).

In particular, this enables us to explain a strong analogy between representations of classical Lie algebras constructed by using the spin representation of a Clifford algebra and the oscillator representation of a Weyl algebra. 
For example, the dual pairs $(\mf{sp}_{2n},O_\ell)$ on $S(\C^n)^{\ot \ell}$ and $(\mf{so}_{2n},O_\ell)$ on $\Lambda(\C^n)^{\ot \ell}$ yield a correspondence between infinite-dimensional irreducible representations of $\mf{sp}_{2n}$ called oscillator representations and finite-dimensional irreducible representations of $\mf{so}_{2n}$ called oscillator representations. The dual pairs $(\mf{so}_{2n},{Sp}_{2\ell})$ on $S(\C^n)^{\ot 2\ell}$ and $(\mf{sp}_{2n},{Sp}_{2\ell})$ on $\Lambda(\C^n)^{\ot 2\ell}$ yield a similar correspondence where the types of Lie algebras are exchanged. Here $S(\C^n)$ and $\Lambda(\C^n)$ denote the symmetric and exterior algebras generated by $\C^n$, respectively. 

Moreover, this kind of connection including the above examples can be deduced as a special case of a more general equivalence called super duality  \cite{CL,CLW} between two parabolic BGG categories of classical Lie algebras of infinite rank. In this equivalence, another parabolic BGG category of a Lie superalgebra of infinite rank interpolates these two categories more directly, where a duality with the finite-dimensional representations of $G$ is no longer available.

\subsection{} 
The purpose of this paper is to generalize the correspondence between irreducible oscillator and finite-dimensional representations of classical Lie algebras in a different direction to nonsemisimple categories over the quantum affine algebras. We consider the dual pairs $(\mf{g},G)$, where $\mf{g}=\mf{sp}_{2n}, \mf{so}_{2n}$ and $G=O_\ell, Sp_{2\ell}$, and generalize these cases to the untwisted quantum affine algebra $U_q(\hat{\mf g})$ of type $X_n^{(1)}$ $(X=C,D)$, where $U_q(\hat{\mf g})$ be the quantum affine algebra of Drinfeld-Jimbo (without derivation) associated to the affine Lie algebra $\hat{\mf g}$.

The theory of finite-dimensional $U_q(\hat{\mf g})$-modules has a rich structure and has been developed very much for the last several decades. One of the most remarkable results obtained recently in this area is that the information on the singularity of (normalized) $R$ matrix provides a deeper understanding of the category of finite-dimensional $U_q(\hat{\mf g})$-modules via monoidal categorification of quantum cluster algebras \cite{HL10,KKOP1,KKOP2}.

This work was motivated by \cite{KwO21}, where an analogy between $q$-oscillator and finite-dimensional representations over quantum affine algebras was first observed: the spectral decomposition of the $R$ matrix on the tensor product of two $q$-oscillator representations of type $C_n^{(1)}$ \cite{KO} coincides with that of the tensor product of two spin representation of type $D_n^{(1)}$ \cite{O90} (up to classically irreducible components). Based on this observation, it was conjectured that there exists a more general class of $q$-oscillator representations of type $C_n^{(1)}$ corresponding to finite-dimensional representations of type $D_n^{(1)}$ together with other pairs of quantum affine (super)algebras of types other than $A$. The case of type $A$ was studied in \cite{KL22}.

There is a $q$-analogue of the dual pairs $(\mf{g},G)$ recently obtained in \cite{ST}, where the action of $G$ is replaced by that of the $\iota$quantum group associated to the Lie algebra of $G$, introduced in the study of quantum symmetric pairs \cite{Le}. 
However, instead of considering an affine analogue of \cite{ST} since an $\iota$quantum group is not a Hopf algebra, we study representations of quantum affine superalgebras, which interpolate more directly the representations of the quantum affine algebras of type $C_n^{(1)}$ and $D_n^{(1)}$ as in the theory of super duality for finite types.

We consider the quantum affine superalgebra $\U(\e)$ of type $D$ (more precisely type $D\!D$ in \cite{Ya99}), where $\e=(\e_i)$ is a finite sequence of $0$ and $1$ representing a fundamental system of positive roots for orthosymplectic Lie superalgebra. Here we use another version of the quantum affine superalgebra due to \cite{KOS}, which we call a generalized quantum group of affine type $D$.

The main result in this paper is to construct a category of $\U(\e)$-modules and a family of irreducible representations, which naturally correspond to finite-dimensional representations of type $X^{(1)}_n$ ($X=C,D$). By using the homomorphic image of the quantum affine algebra of type $X_n^{(1)}$ in $\U(\e)$, we also define a category of infinite-dimensional representations of type $X_n^{(1)}$ and construct a family of irreducible representations, which can be viewed as an affine analogue of the semisimple categories generated by the $q$-deformed oscillator representations of $\mf{sp}_{2n}$ and $\mf{so}_{2n}$. We remark that all of these categories and their irreducible objects are constructed in a uniform way. As a byproduct, we show that the irreducible $q$-oscillator representations of type $C_n^{(1)}$ (resp. $D_n^{(1)}$) and the finite-dimensional irreducible representations of type $D_n^{(1)}$ (resp. $C_n^{(1)}$) are interpolated by irreducible $\U(\e)$-modules in this category under certain exact monoidal functors which are called truncation in super duality. 

All the results in this paper can be done using the generalized quantum group of affine type $C$ (or type $CC$ in \cite{Ya99}) (see Section \ref{subsec:conjectures}).

\subsection{}
Let us explain our results in more details. First consider the case for the $q$-oscillator representations of type $C$ and the finite-dimensional representations of type $D$.

Let ${\bm\e}=(1,0,\dots,0,1)$ be an alternating sequence of length $n=2m+1$, and let $\un{\bm\e}$ (resp. $\ov{\bm\e}$) be the subsequence of $\bm\e$ obtained by removing all $1$ (resp. $0$). We observe that there exist natural homomorphisms of algebras $\phi_\e : \U_{{\e}}\longrightarrow\U_{{\bm\e}}$ ($\e=\un{\bm\e},\ov{\bm\e}$), where we put $\U_{{\bm\e}}=\U(\bm\e)$, $\U_{\un{\bm\e}}=U_q(C_{m}^{(1)})$, and $\U_{\ov{\bm\e}}=U_{\td{q}}(D_{m+1}^{(1)})$ where $\td{q}=-q^{-1}$. 

Let $\W_{\e}$ be a $q$-analogue of the symmetric algebra generated by the superspace with a $\Z_2$-graded basis parametrized by $\e$ ($\e={\bm\e},\un{\bm\e},\ov{\bm\e}$), where we regard $\W_\e\subset \W_{\bm\e}$ for $\e=\un{\bm\e},\ov{\bm\e}$. 
We define a $\U_\e$-module structure on $\W_\e$ ($\e={\bm\e},\un{\bm\e},\ov{\bm\e}$) with spectral parameter $x$, say $\mathcal{W}_{\epsilon}(x)$, where the actions of $\U_\e$ ($\e=\un{\bm\e},\ov{\bm\e}$) are given through $\phi_\e$. Note that $\W_{\un{\boldsymbol{\e}}}(x)$ is a $q$-analogue of the oscillator representation of $\mf{sp}_{2m}$ and  $\W_{\ov{\boldsymbol{\e}}}(x)$ is the sum of two spin representations $\W^\pm_{\ov{\boldsymbol{\e}}}(x)$ of type $\mf{so}_{2m+2}$.
We show that $\W_\e(x)^{\ot \ell}$ is semisimple as a module over the subalgebra $\mathring{\U}_\e$ of finite type and classify the irreducible $\mathring{\U}_\e$-modules in $\W_\e(x)^{\ot \ell}$ for $\ell\ge 1$, which generate a semisimple tensor category ${\mc O}^{\,\mf c}_{\rm osc,{\e}}$. An irreducible representation in ${\mc O}^{\,\mf c}_{\rm osc,{\un{\bm\e}}}$ is a $q$-analogue of the oscillator representation of type $C_m$ with the same highest weight, while an irreducible representation in ${\mc O}^{\,\mf c}_{\rm osc,{\ov{\bm\e}}}$ is a finite-dimensional representation of type $D_{m+1}$.

Now we define $\widehat{\mc O}^{\,\mf c}_{\rm osc,{\e}}$ to be the category of $\U_\e$-modules which belong to ${\mc O}^{\,\mf c}_{\rm osc,{\e}}$ as $\mathring{\mathcal{U}}_\epsilon$-modules and satisfy some mild conditions on their weights. 
Note that $\widehat{\mc O}^{\,\mf c}_{\rm osc,{\ov{\bm\e}}}$ is the category of finite-dimensional representations of type $D_{m+1}^{(1)}$.
By using $\phi_\e$ ($\e=\un{\bm\e},\ov{\bm\e}$), we can define natural monoidal functors $\mf{tr}_{\e}$ ($\e=\un{\bm\e},\ov{\bm\e}$) which send an irreducible representation in ${\mc O}^{\,\mf c}_{\rm osc,{\bm\e}}$ to the one in ${\mc O}^{\,\mf c}_{\rm osc,{\e}}$ or zero.
The functors $\mf{tr}_{\e}$ preserves $\W^\pm_\e(x)$ and the normalized $R$ matrix on the tensor product $\W^\pm_\e(x)\ot \W^\pm_\e(y)$, which implies that the spectral decomposition is essentially the same for all $\e$. 
Following \cite{KKKO15,Kas02}, we construct a family of representations $\W_\e({\bm\la})$ in $\widehat{\mc O}^{\,\mf c}_{\rm osc,{\e}}$ with parameters ${\bm\la}\in \mc{P}^+$ independent of $\e$ such that $\mf{tr}_{\e}(\W_{\bm\e}({\bm\la}))=\W_\e({\bm\la})$ ($\e=\un{\bm\e},\ov{\bm\e}$) and any nonzero $\W_\e({\bm\la})$ is irreducible.

We remark that the spectral decompositions in $\widehat{\mc O}^{\,\mf c}_{\rm osc,{\ov{\bm\e}}}$ and $\widehat{\mc O}^{\,\mf c}_{\rm osc,{\un{\bm\e}}}$ were already obtained in \cite{O90} and \cite{KO}, respectively. Our approach explains why these formulas coincide thanks to their connection with $\widehat{\mc O}^{\,\mf c}_{\rm osc,{{\bm\e}}}$.

\subsection{}
Next consider the case for the $q$-oscillator representations of type $D$ and the finite-dimensional representations of type $C$.

Let ${\bm\e}'=(0,1,\dots,1,0)$ be an alternating sequence of length $n=2m+1$, and let $\un{\bm\e}'$ (resp. $\ov{\bm\e}'$) be the subsequence of $\bm\e'$ obtained by removing all $1$ (resp. $0$). In this case, we have homomorphism of algebras $\phi_\e : \U_{{\e}}\longrightarrow\U_{{\bm\e}'}$ ($\e=\un{\bm\e}',\ov{\bm\e}'$), where we put $\U_{{\bm\e}'}=\U(\bm\e')$, $\U_{\un{\bm\e}'}=U_q(D_{m+1}^{(1)})$, and $\U_{\ov{\bm\e}'}=U_{\td{q}}(C_{m}^{(1)})$.
We define a $\U_\e$-module $\WW_\e(x)$ with spectral parameter $x$, where $\WW_{\e}=\W_{\e}^{\ot 2}$ ($\e={\bm\e}',\un{\bm\e}',\ov{\bm\e}'$). 
We show that $\WW_\e(x)^{\ot \ell}$ is semisimple as a module over the subalgebra $\mathring{\U}_\e$ of finite type and classify the irreducible $\mathring{\U}_\e$-modules in $\WW_\e(x)^{\ot \ell}$ for $\ell\ge 1$. 
The categories ${\mc O}^{\,\mf d}_{\rm osc,{\e}}$ and $\widehat{\mc O}^{\,\mf d}_{\rm osc,{\e}}$ are defined in the same way, where $\widehat{\mc O}^{\,\mf d}_{\rm osc,{\ov{\bm\e}}}$ is the category of finite-dimensional representations of type $C_{m}^{(1)}$. 

The proofs of the main results in this case are technically more difficult. We first construct a fundamental type representation $\W_{l,{\e}}(x)$ for $l\ge 0$ ($\e={\bm\e}',\un{\bm\e}',\ov{\bm\e}'$) such that $\mf{tr}_{\e}(\W_{l,{\bm\e}}(x))=\W_{l,{\e}}(x)$ ($\e=\un{\bm\e}',\ov{\bm\e}'$), where $\W_{l,{\ov{\bm\e}'}}(x)$ is isomorphic to a fundamental representation of type $C_m^{(1)}$. Then we prove the irreducibility of $\W_{l,{\e}}(x)\ot \W_{m,{\e}}(y)$ for generic $x,y$, which is known in the case of $\ov{\bm\e}'$ by the theory of crystal base \cite{AK,Kas02}, and compute explicitly the spectral decomposition of the normalized $R$ matrix on $\W_{l,{\e}}(x)\ot \W_{m,{\e}}(y)$ ($\e={\bm\e}',\un{\bm\e}',\ov{\bm\e}'$). It is rather surprising that this spectral decomposition has not been known so far even for $\ov{\bm\e}'$, while only its denominator is known \cite{AK}. Then we obtain the results for $\widehat{\mc O}^{\,\mf d}_{\rm osc,{\e}}$ parallel to the case of $\widehat{\mc O}^{\,\mf c}_{\rm osc,{\e}}$.

\subsection{}
Summarizing, we have constructed two types of categories $\widehat{\mc O}^{\,\mf c}_{\rm osc,{\e}}$ and $\widehat{\mc O}^{\,\mf d}_{\rm osc,{\e}}$ where $(\e,\un{\e},\ov{\e})=({\bm\e},\un{\bm\e},\ov{\bm\e})$ and $({\bm\e}',\un{\bm\e}',\ov{\bm\e}')$ for $\mf{x}=\mf{c}$ and $\mf{d}$, respectively, together with the functors $\mf{tr}_\e$ between them:
\begin{equation}\label{eq:triangle of truncation osc x-0}
\xymatrixcolsep{3pc}\xymatrixrowsep{0.3pc}\xymatrix{
 & \widehat{\mc O}^{\,\mf x}_{\rm osc,{\e}}  \ar@{->}_{\mf{tr}_{\un{\e}}}[dl]\ar@{->}^{\mf{tr}_{\ov{\e}}}[dr] &  \\
 \widehat{\mc O}^{\,\mf x}_{\rm osc,\un{\e}} & &  \widehat{\mc O}^{\,\mf x}_{\rm osc,\ov{\e}}
}\quad 
\xymatrixcolsep{3pc}\xymatrixrowsep{0.8pc}\xymatrix{
 & \W_{\e}({\mb \la})  \ar@{|->}_{\mf{tr}^\e_{\un{\e}}}[dl]\ar@{|->}^{\mf{tr}^\e_{\ov{\e}}}[dr] &  \\
 \W_{\un{\e}}({\mb\la}) & &  \W_{\ov{\e}}({\mb \la})
}
\end{equation}
where $\W_\e({\bm\la})$ is a representation in $\widehat{\mc O}^{\,\mf x}_{\rm osc,{\e}}$, which is irreducible if it is not zero. 
Here the leftmost category in the diagram consists of $q$-oscillator representations of type $C_{m}^{(1)}$ (resp. $D_{m+1}^{(1)}$) and the rightmost one consists of finite-dimensional representations of type $D_{m+1}^{(1)}$ (resp. $C_{m}^{(1)}$) for $\mf{x}=\mf{c}$ (resp. $\mf{x}=\mf{d}$). 

Finally we remark that the irreducible representations $\W_\e({\bm\la})$ of the quantum affine superalgebra in $\widehat{\mc O}^{\,\mf x}_{\rm osc,{\e}}$ for $\e={\bm\e}, {\bm\e}'$ together with those in $\widehat{\mc O}^{\,\mf x}_{\rm osc,{\e}}$ for $\e=\un{\bm\e}, \un{\bm\e}'$ of type $X_n^{(1)}$ are new. The cases of the quantum affine superalgebra with $\e\neq {\bm\e}, {\bm\e}'$ are also discussed. As far as we know, very little seems to be known so far about the representations of quantum affine superalgebras of types other than $A$. 
We expect that there are more connections between the categories in \eqref{eq:triangle of truncation osc x-0} (cf.~\cite{KL20}), and hence our approach would be of help to have a better understanding of representations of quantum affine superalgebras including the $q$-characters of irreducibles, the structure of Grothendieck rings, $T$-systems and so on. 
We also expect that the results in this paper can be extended to other types of quantum affine (super)algebras, which is discussed with more details at the end of the paper.

The paper is organized as follows. In Section \ref{sec:definitions}, we review necessary background for the quantum affine superalgebras. In Section \ref{sec:truncation}, we introduce truncation functors of type $\mf{c}$ and $\mf{d}$. In Section \ref{sec:osc rep}, we review the classical results on the dual pair $(\mf{g},G)$, where $\mf{g}=\mf{sp}_{2n}, \mf{so}_{2n}$ and $G=O_\ell, Sp_{2\ell}$. In Section \ref{sec:category Osc}, we define the categories $\widehat{\mc O}^{\,\mf x}_{\rm osc,{\e}}$. In Sections \ref{sec:irr of type c} and \ref{sec:irr of type d}, we construct irreducible representations in $\widehat{\mc O}^{\,\mf c}_{\rm osc,{\e}}$ and $\widehat{\mc O}^{\,\mf d}_{\rm osc,{\e}}$, respectively. In Section \ref{sec:remarks}, we give some remarks and questions related to other types of quantum affine superalgebras. \smallskip

\noindent {\bf Acknowledgement}
Part of this work has been done while the first two authors were visiting OCAMI in Osaka Metropolitan University during in Japan, January 2023. They would like to thank OCAMI for the invitation and warm hospitality during their visit. We also would like to thank H. Yamane for his kind explanation on his works.

\section{Quantum affine superalgebras of type $D$}\label{sec:definitions}

\subsection{Generalized quantum groups of affine type $D$}\label{subsec:def}
We assume that $n$ is a positive integer such that $n\geq 4$.
We denote by $\mathbb{Z}_+$ the set of nonnegative integers. 
Let $\Bbbk=\mathbb{C}(q^{\hf})$ and $\Bbbk^{\times}=\Bbbk \setminus \left\{0\right\}$, where $q^{\hf}$ is an indeterminate.
We put 
\begin{equation*}
\begin{split}
[m]&=\frac{q^m-q^{-m}}{q-q^{-1}}\quad (m\in \Z_{+}),\\
[m]!&=[m][m-1]\cdots [1]\quad (m\geq 1),\quad [0]!=1,\\
\begin{bmatrix} m \\ k \end{bmatrix}&= \frac{[m][m-1]\cdots [m-k+1]}{[k]!}\quad (0\leq k\leq m).
\end{split}
\end{equation*}

We assume the following notations: 

\begin{itemize}

\item[$\bullet$] $\e=(\e_1,\dots,\e_n)$ : a sequence with $\e_i\in \{0,1\}$ ($i=1,\dots, n$),
 
\item[$\bullet$] $\mathbb{I}= \{1<2<\cdots <n\}$, 

\item[$\bullet$] $P = \Z\La \oplus \Z\de_1\oplus\cdots\oplus \Z\de_n$: a $\Z$-lattice with a symmetric bilinear form $(\,\cdot\,|\,\cdot\,)$ satisfying 
\begin{equation*}\label{eq:extended bilinear form}
(\de_i|\de_j)=(-1)^{\e_i}\delta_{ij}, \quad 
(\de_i|\La)=-\hf,
\end{equation*}

\item[$\bullet$] $I=\{0,1,\ldots,n\}$,

\item[$\bullet$] $\alpha_i=\de_{i+1}-\de_i \in P$ $(i\in I\setminus\{0,n\})$, 
$\alpha_0=\de_2+\de_{1}$, $\alpha_n=-\de_{n}-\de_{n-1}$,

\item[$\bullet$] $I_{\rm even}=\{\,i\in I\,|\,(\alpha_i|\alpha_i)=\pm 2\,\}$, 
$I_{\rm odd}=\{\,i\in I\,|\,(\alpha_i|\alpha_i)=0\,\}$,

\item[$\bullet$] $q_i=\si q^{\si}$ $(i\in \mathbb{I})$, that is,
\begin{equation*}
q_i=
\begin{cases}
q & \text{if $\e_i=0$},\\
\td{q}:=-q^{-1} & \text{if $\e_i=1$},\\
\end{cases} \quad (i\in \mathbb{I}),
\end{equation*}

\item[$\bullet$] ${\bf v}={\bf i}q^{\hf}$ (${\bf i}=\sqrt{-1}$),

\item[$\bullet$] ${\bq}(\,\cdot\,,\,\cdot\,)$: a symmetric biadditive function from $P\times P$ to $\Bbbk^{\times}$ given by
\begin{equation*}
\begin{split}
&\bq(\mu,\nu) = {{\bf v}}^{\sum_{j\in \I}(\ell'\mu_j+\ell\nu_j)}\prod_{i\in \mathbb{I}}q_i^{\mu_i \nu_i} \quad (\mu,\nu\in P),
\end{split}
\end{equation*}
for $\mu=\ell\La+\sum_{i\in \I}\mu_i\de_i$ and $\nu=\ell'\La+\sum_{i\in \I}\nu_i\de_i$.

\end{itemize}

Let $\mathcal{U}_D(\ep)$ be the associative $\Bbbk$-algebra generated by 
$e_i,f_i,k_\mu$  $(i\in I, \mu\in P)$ with the relations below: 
{\allowdisplaybreaks
\begin{gather*}
k_\mu=1 \quad(\mu=0), \quad k_{\mu +\mu'}=k_{\mu}k_{\mu'} \quad (\mu, \mu' \in P),\label{eq:Weyl-rel-1-e} \\ 
k_\mu e_i k_{-\mu}=\bq(\mu,\alpha_i)e_i,\quad 
k_\mu f_i k_{-\mu}=\bq(\mu,\alpha_i)^{-1}f_i\quad (i\in I, \mu\in P), \label{eq:Weyl-rel-2-e} \\ 
e_if_j - f_je_i =\delta_{ij}\frac{k_{i} - k_{i}^{-1}}{q-q^{-1}}\quad (i,j\in I),\label{eq:Weyl-rel-3-e}\\
e_i^2= f_i^2 =0 \quad \text{($i\in I_{\rm odd}$)},\label{eq:Weyl-rel-4}\\
[e_i ,e_j]=[f_i ,f_j]=0 \quad\text{if $(\alpha_i,\alpha_j)=0$},\\
e^2 _0 e_2 - (-1)^{\ep_1}[2]e_0 e_2 e_0 +e_2 e^2 _0 =(e \rightarrow f)=0 
\quad \text{if $\e_1=\e_2$},\\
e^2 _2 e_0 - (-1)^{\ep_2}[2]e_2 e_0 e_2 +e_0 e^2 _2 =(e \rightarrow f)=0 
\quad \text{if $\e_2=\e_3$}, \\
e^2 _i e_j - (-1)^{\ep_i}[2]e_i e_j e_i +e_j e^2 _i =(e \rightarrow f)=0
\quad \text{if $i,j\in I\setminus\{0,n\}$, $|i-j|=1$, $\e_i=\e_{i+1}$}, \\
e^2 _{n-2} e_n - (-1)^{\ep_{n-2}}[2]e_{n-2} e_n e_{n-2} +e_n e^2 _{n-2} =(e \rightarrow f)=0 
\quad \text{if $\e_{n-2}=\e_{n-1}$}, \\
e^2 _{n} e_{n-2} - (-1)^{\ep_{n-1}}[2]e_n e_{n-2} e_n +e_{n-2} e^2 _n =(e \rightarrow f)=0 
\quad \text{if $\e_{n-1}=\e_n$}, \\
e_0e_1e_2-e_1e_0e_2+(-1)^{\ep_2}[2](e_1e_2e_0-e_0e_2e_1)
+e_2e_0e_1-e_2e_1e_0=(e\rightarrow f)=0 \quad \text{if $\e_1\neq \e_2$}, \\
e_0e_2e_3e_2-e_3e_2e_0e_2+(-1)^{\ep_3}[2]e_2e_3e_0e_2
-e_2e_0e_2e_3+e_2e_3e_2e_0=(e \rightarrow f)=0
\quad \text{if $\e_2\neq \e_3$}, \\
e_{i}e_{i-1}e_{i}e_{i+1}-e_{i}e_{i+1}e_{i}e_{i-1}+(-1)^{\ep_i}[2]e_{i}e_{i-1}e_{i+1}e_{i} \\
\hspace{0.1cm}-e_{i-1}e_{i}e_{i+1}e_{i}+e_{i+1}e_{i}e_{i-1}e_{i}=(e \rightarrow f)=0
\quad\text{if $i\in \{2,\dots,n-2\}$, $\e_i\neq \e_{i+1}$}, \nonumber \\
e_{n-2}e_{n-3}e_{n-2}e_{n}-e_{n-2}e_{n}e_{n-2}e_{n-3}+(-1)^{\ep_{n-2}}[2]e_{n-2}e_{n-3}e_{n}e_{n-2} \\
\hspace{0.1cm}-e_{n-3}e_{n-2}e_{n}e_{n-2}+e_{n}e_{n-2}e_{n-3}e_{n-2}=(e \rightarrow f)=0
\quad \text{if $\e_{n-2}\neq \e_{n-1}$},\nonumber \\
e_{n-2}e_{n}e_{n-1}-e_{n-2}e_{n-1}e_{n}+(-1)^{\ep_{n-1}}[2](e_{n-1}e_{n-2}e_{n}-e_{n}e_{n-2}e_{n-1})\\
\hspace{0.1cm}+e_{n}e_{n-1}e_{n-2}-e_{n-1}e_{n}e_{n-2}=(e\rightarrow f)=0 \quad \text{if $\e_{n-1}\neq \e_n$}\nonumber.
\end{gather*}}
Here $k_i=k_{\alpha_i}$ for $i\in I$, and $(e\rightarrow f)$ denotes the relation which is obtained by replacing $e_i$ with $f_i$.

We call $\U_D(\e)$ the {\em generalized quantum group of affine type $D$ associated to $\e$} \cite{Ma} (cf. \cite{KOS}). It has a Hopf algebra structure with the comultiplication $\Delta$ and the counit $S$
\begin{equation}\label{eq:comult-e}
\begin{split}
\Delta(k_\mu)&=k_\mu\otimes k_\mu, \\ 
\Delta(e_i)&= 1\ot e_i + e_i\ot k_i^{-1}, \\
\Delta(f_i)&= f_i\ot 1 + k_i\ot f_i , \\  
S(k_\mu)=k_{-\mu},& \ \ S(e_i)=-e_ik_i, \ \  S(f_i)=-k_i^{-1}f_i,\\
\end{split}
\end{equation}
for $\mu\in P$ and  $i\in I$.

Note that $\U_D(0^n)\cong U_q(D_n^{(1)})$ and $\U_D(1^n)\cong U_{\td{q}}(D_n^{(1)})$ (up to Cartan part), where $U_q(D_n^{(1)})$ is the quantum affine algebra of type $D_n^{(1)}$. In general, one may associate a Dynkin diagram with each $\e$, where
\begin{itemize}
 \item[(1)] the vertices are indexed by $I$, and denoted by $\bigcirc$ if $(\alpha_i|\alpha_i)\neq 0$ and $\bigotimes$ if $(\alpha_i|\alpha_i)= 0$,
 
 \item[(2)] two vertices $i$ and $j$ are connected by an edge whenever $(\alpha_i|\alpha_j)\neq 0$.
\end{itemize}
For example, we have 
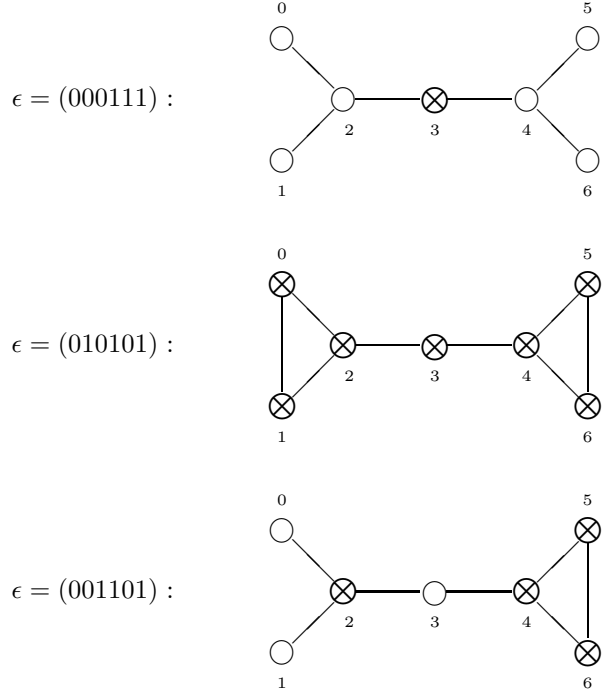
\begin{center}\setlength{\unitlength}{0.16in} \smallskip
\begin{picture}(22,5)
\put(0,2){\makebox(0,0)[c]{$\e=(000111)$ : }}
\put(6,0){\makebox(0,0)[c]{$\bigcirc$}}
\put(6,4){\makebox(0,0)[c]{$\bigcirc$}}
\put(8,2){\makebox(0,0)[c]{$\bigcirc$}}
\put(11,1.95){\makebox(0,0)[c]{$\bigotimes$}}
\put(14,2){\makebox(0,0)[c]{$\bigcirc$}}
\put(16,0){\makebox(0,0)[c]{$\bigcirc$}}
\put(16,4){\makebox(0,0)[c]{$\bigcirc$}}

\put(6.35,0.3){\line(1,1){1.35}} \put(6.35,3.7){\line(1,-1){1.35}}
\put(8.4,2){\line(1,0){2.1}} \put(11.4,2){\line(1,0){2.1}}

\put(14.35,2.3){\line(1,1){1.35}} \put(14.35,1.65){\line(1,-1){1.35}}

\put(6,5){\makebox(0,0)[c]{\tiny $0$}}
\put(6,-1){\makebox(0,0)[c]{\tiny $1$}}
\put(8.2,1){\makebox(0,0)[c]{\tiny $2$}}
\put(11,1){\makebox(0,0)[c]{\tiny $3$}}
\put(14,1){\makebox(0,0)[c]{\tiny $4$}}
\put(16,5){\makebox(0,0)[c]{\tiny $5$}}
\put(16,-1){\makebox(0,0)[c]{\tiny $6$}}
\end{picture}
\end{center}\vskip 10mm

\begin{center}\setlength{\unitlength}{0.16in} \smallskip
\begin{picture}(22,5)
\put(0,2){\makebox(0,0)[c]{$\e=(010101)$ : }}
\put(6,0){\makebox(0,0)[c]{$\bigotimes$}}
\put(6,4){\makebox(0,0)[c]{$\bigotimes$}}
\put(8,2){\makebox(0,0)[c]{$\bigotimes$}}
\put(11,1.95){\makebox(0,0)[c]{$\bigotimes$}}
\put(14,2){\makebox(0,0)[c]{$\bigotimes$}}
\put(16,0){\makebox(0,0)[c]{$\bigotimes$}}
\put(16,4){\makebox(0,0)[c]{$\bigotimes$}}

\put(6.35,0.3){\line(1,1){1.35}} \put(6.35,3.7){\line(1,-1){1.35}}
\put(8.4,2){\line(1,0){2.1}} \put(11.4,2){\line(1,0){2.1}}

\put(6,0.5){\line(0,1){3.1}}\put(16,0.5){\line(0,1){3.1}} 

\put(14.35,2.3){\line(1,1){1.35}} \put(14.35,1.65){\line(1,-1){1.35}}

\put(6,5){\makebox(0,0)[c]{\tiny $0$}}
\put(6,-1){\makebox(0,0)[c]{\tiny $1$}}
\put(8.2,1){\makebox(0,0)[c]{\tiny $2$}}
\put(11,1){\makebox(0,0)[c]{\tiny $3$}}
\put(14,1){\makebox(0,0)[c]{\tiny $4$}}
\put(16,5){\makebox(0,0)[c]{\tiny $5$}}
\put(16,-1){\makebox(0,0)[c]{\tiny $6$}}
\end{picture}
\end{center}\vskip 10mm

\begin{center}\setlength{\unitlength}{0.16in} \smallskip
\begin{picture}(22,5)
\put(0,2){\makebox(0,0)[c]{$\e=(001101)$ : }}
\put(6,0){\makebox(0,0)[c]{$\bigcirc$}}
\put(6,4){\makebox(0,0)[c]{$\bigcirc$}}
\put(8,2){\makebox(0,0)[c]{$\bigotimes$}}
\put(11,1.95){\makebox(0,0)[c]{$\bigcirc$}}
\put(14,2){\makebox(0,0)[c]{$\bigotimes$}}
\put(16,0){\makebox(0,0)[c]{$\bigotimes$}}
\put(16,4){\makebox(0,0)[c]{$\bigotimes$}}

\put(6.35,0.3){\line(1,1){1.35}} \put(6.35,3.7){\line(1,-1){1.35}}
\put(8.4,2){\line(1,0){2.1}} \put(11.4,2){\line(1,0){2.1}}
\put(14.35,2.3){\line(1,1){1.35}} \put(14.35,1.65){\line(1,-1){1.35}}
\put(16,0.5){\line(0,1){3.1}}

\put(6,5){\makebox(0,0)[c]{\tiny $0$}}
\put(6,-1){\makebox(0,0)[c]{\tiny $1$}}
\put(8.2,1){\makebox(0,0)[c]{\tiny $2$}}
\put(11,1){\makebox(0,0)[c]{\tiny $3$}}
\put(14,1){\makebox(0,0)[c]{\tiny $4$}}
\put(16,5){\makebox(0,0)[c]{\tiny $5$}}
\put(16,-1){\makebox(0,0)[c]{\tiny $6$}}
\end{picture}
\end{center}\vskip 5mm

Let $\ring{\mathcal{U}}_D(\e)$ be the subalgebra generated by $e_i$, $f_i$, and $k_\mu$ for $i\in I\setminus \{0\}$ and $\mu\in P$.
If we put 
\begin{equation}\label{eq:fundamental weights}
\varpi_n=\La,\quad \varpi_{n-1}=\La+\td{\de}_n,\quad \varpi_i=2\La+\td{\de}_{n}+\dots+\td{\de}_{i+1}\quad (1\le i\le n-2), 
\end{equation}
where $\td{\de}_a=(-1)^{\e_a}\de_a$ for $a\in \I$, then $(\varpi_i|\alpha_j)=\de_{ij}$ for $i,j\in I\setminus\{0\}$.

For a $\U_D(\e)$-module $V$ and $\la\in P$, we define the $\la$-weight space to be 
\begin{equation}\label{eq:wt space}
 V_\la = \{\,v\in V\,|\,k_{\mu} v= {\bq}(\la,\mu) v \ \ (\mu\in P) \,\},
\end{equation}
and write ${\rm wt}(v)=\la$ for $v\in V_\la$.
Let ${\rm wt}(V)=\{\,\mu\in P\,|\,V_\mu\neq 0\,\}$. 

\subsection{Quantum affine superalgebras of type $D$}

Let us recall the notion of the quantum affine superalgebra $U_D(\e)$ of type $D$ associated with $\e$ \cite{Ya99}. 
The algebra $U_D(\epsilon)$ is the associative $\Bbbk$-algebra with 1 generated by $E_{i}$, $F_{i}$, $K_{\mu}$ $(\mu\in P, i\in I)$ satisfying the following relations:
\allowdisplaybreaks{\small{}
\begin{gather*}
  K_{\mu}=1 \quad(\mu=0),\quad K_{\mu+\mu'}=K_{\mu}K_{\mu'} \quad (\mu, \mu' \in P),\\
  K_{\mu}E_{i}K_{\mu}^{-1}=q^{(\mu|\alpha_{i})}E_{i},\quad K_{\mu}F_{i}K_{\mu}^{-1}=q^{-(\mu|\alpha_{i})}F_{i} \quad (i\in I, \mu\in P),\\
  E_{i}F_{j}-(-1)^{p(i)p(j)}F_{j}E_{i}=(-1)^{\epsilon_{i}}\delta_{ij}\frac{K_{\alpha_{i}}-K_{-\alpha_{i}}}{q-q^{-1}} \quad (i,j\in I),\\
  E_{i}^{2}=F_{i}^{2}=0\quad\text{($i\in I_{\rm odd}$)},\\
  E_{i}E_{j}-(-1)^{p(i)p(j)}E_{j}E_{i}=(E\rightarrow F)=0\quad\text{if $(\alpha_i,\alpha_j)=0$},\\
  E_{0}^{2}E_{2}-[2]E_{0}E_{2}E_{0}+E_{2}E_{0}^{2}=(E\rightarrow F)=0\quad\text{ if $\e_1=\e_2$},\\
  E_{2}^{2}E_{0}-[2]E_{2}E_{0}E_{2}+E_{0}E_{2}^{2}=(E\rightarrow F)=0\quad\text{ if $\e_2=\e_3$},\\
  E_{i}^{2}E_{j}-[2]E_{i}E_{j}E_{i}+E_{j}E_{i}^{2}=(E\rightarrow F)=0\quad \text{if $i,j\in I\setminus\{0,n\}$, $|i-j|=1$, $\e_i=\e_{i+1}$},\\
  E_{n-2}^{2}E_{n}-[2]E_{n-2}E_{n}E_{n-2}+E_{n}E_{n-2}^{2}=(E\rightarrow F)=0\quad\text{ if }\epsilon_{n-2}=\epsilon_{n-1},\\
  E_{n}^{2}E_{n-2}-[2]E_{n}E_{n-2}E_{n}+E_{n-2}E_{n}^{2}=(E\rightarrow F)=0\quad\text{ if }\epsilon_{n-1}=\epsilon_{n},\\
  E_{0}E_{1}E_{2}-E_{1}E_{0}E_{2}+(-1)^{\epsilon_{2}+\epsilon_{3}}(E_{1}E_{2}E_{0}-E_{0}E_{2}E_{1})+E_{2}E_{0}E_{1}-E_{2}E_{1}E_{0}\\ =(E\rightarrow F)=0\quad\text{ if \ensuremath{\epsilon_{1}\neq\epsilon_{2}}},\\
  \left[E_{2},\left[\left[E_{0},E_{2}\right]_{(-1)^{p(0)}q},E_{3}\right]_{(-1)^{(p(0)+p(2))p(3)}q^{-1}}\right]_{(-1)^{p(0)+p(2)+p(3)}}=(E\rightarrow F)=0\quad\text{ if }\epsilon_{2}\neq\epsilon_{3},\\
  \left[E_{i},\left[\left[E_{i-1},E_{i}\right]_{(-1)^{p(i-1)}q},E_{i+1}\right]_{(-1)^{(p(i-1)+p(i))p(i+1)}q^{-1}}\right]_{(-1)^{p(i-1)+p(i)+p(i+1)}}\\
  \hfill=(E\rightarrow F)=0\quad\text{ if }\epsilon_{i}\neq\epsilon_{i+1},i\in[2,n-2],\\
  \left[E_{n-2},\left[\left[E_{n-3},E_{n-2}\right]_{(-1)^{p(n-3)}q},E_{n}\right]_{(-1)^{(p(n-3)+p(n-2))p(n)}q^{-1}}\right]_{(-1)^{p(n-3)+p(n-2)+p(n)}}\\
  \hfill=(E\rightarrow F)=0\quad\text{ if }\epsilon_{n-2}\neq\epsilon_{n-1},\\
  E_{n}E_{n-1}E_{n-2}-E_{n-1}E_{n-0}E_{n-2}+(-1)^{\epsilon_{n-1}+\epsilon_{n-2}}(E_{n-1}E_{n-2}E_{n}-E_{n}E_{n-2}E_{n-1})\\
  \quad+E_{n-2}E_{n}E_{n-1}-E_{n-2}E_{n-1}E_{n}=(E\rightarrow F)=0\quad\text{ if \ensuremath{\epsilon_{n-1}\neq\epsilon_{n}}},
\end{gather*}
}where $[X,Y]_{t}=XY-tYX$ and $p(i)=0$ (resp.$1$) for $i\in I_{\rm even}$ (resp. $i\in I_{\rm odd}$). More precisely, $U_D(\e)$ corresponds to type (DD) in \cite[Section 1.4 and Proposition 6.7.1]{Ya99}.

The two algebras $\U_D(\e)$ and $U_D(\e)$ are closely related as in the case of type $A$ \cite{KL20}. Let us explain it as follows.
Let
\begin{itemize}
\item[$\bullet$] $P^e_{\rm af} = \bigoplus_{j\in \I}\Z\hf\mb{\delta}_j \oplus \Z\mb{\delta} \oplus \Z \Lambda_0$: a $\Z$-lattice
with a symmetric bilinear form 
\begin{equation*}\label{eq:extended bilinear form}
(\mb{\delta}_i|\mb{\delta}_j)=(-1)^{\e_i}\delta_{ij},\quad (\mb{\delta}_i|\mb{\delta})=(\mb{\delta}_i|\Lambda_0)=(\Lambda_0|\Lambda_0)=(\mb{\delta}|\mb{\delta})=0,\quad (\Lambda_0|\mb{\delta})=1,
\end{equation*}

\item[$\bullet$] $\mb{\alpha}_i=\mb{\de}_{i+1}-\mb{\de}_i$ $(i\in I\setminus\{0,n\})$, $\mb{\alpha}_0=\mb{\de}_2+\mb{\de}_{1}+\mb{\de}$, $\mb{\alpha}_n=-\mb{\de}_{n}-\mb{\de}_{n-1}$,

\item[$\bullet$] $Q=\bigoplus_{i\in I}\Z\mb{\alpha}_i$, $Q_{\pm}=\pm\bigoplus_{i\in I}\Z_+\mb{\alpha}_i$.

\end{itemize}
Note that if we put $$\La_n=\La_0-\hf\sum_{i\in \I}(-1)^{\e_i}{\mb\de}_i,$$ 
then we have
\begin{equation*}
(\mb{\delta}_i|\Lambda_n)=-\hf,\quad (\Lambda_n |\Lambda_n)=\frac{1}{4}\sum_{i\in\mathbb{I}}(-1)^{\epsilon_i}.
\end{equation*}  
Hence we may regard $P\subset P^e_{\rm af}$ by identifying $\de_i$ and $\mb{\de}_i$ ($i\in \I$), and $\La$ with $\La_n$, respectively. 
Let ${\rm cl} : P^e_{\rm af} \longrightarrow P^{e}_{\mathrm{af}}/\mathbb{Z}\boldsymbol{\delta}$ be the projection under this identification.

We define $\U_D(\e)^e_{\rm af}$ to be the $\Bbbk$-algebra as in $\U_D(\e)$, which is defined in the same way as $\mathcal{U}_D (\epsilon)$ where $P$ and $\bq(\,\cdot\,,\,\cdot\,)$ are now replaced by $P^e_{\rm af}$ and $\bq(\mu,\nu)=\bq({\rm cl}(\mu),{\rm cl}(\nu))$ 
for $\mu,\nu\in P^e_{\rm af}$.

Let $\Sigma$ be the $\Bbbk$-bialgebra generated by $\sigma_i$ $(i\in \I)$ such that $\sigma_i\sigma_j=\sigma_j\sigma_i$ and $\sigma_j^2=1$ for $i,j\in \I$, with the comultiplication $\Delta(\sigma_i)=\sigma_i\otimes \sigma_i$ for $i\in \I$.
We define ${\U_D(\e)}^e_{\rm af}[\sigma]$ to be the semidirect product of ${\U_D(\e)}^e_{\rm af}$ and $\Sigma$, where $\Sigma$ acts on ${\U_D(\e)}^e_{\rm af}$ by
\begin{equation*}\label{eq:sigma-rel-2}
\begin{split}
&\sigma_j k_\mu =k_\mu,\quad
\sigma_je_i=(-1)^{\e_j(\mb{\de}_j|\mb{\alpha}_i)}e_i,\quad 
\sigma_jf_i=(-1)^{\e_j(\mb{\de}_j|\mb{\alpha}_i)}f_i,
\end{split}
\end{equation*}
for $j\in \I$, $\mu\in P^e_{\rm af}$ and $i\in I$.  
We define ${U_D(\e)}^e_{\rm af}[\sigma]$ in the same manner.

Suppose that $\e\neq (0^n), (1^n)$.
Let $\mathbb{I}=\mathbb{I}^{(1)}\sqcup\cdots\sqcup\mathbb{I}^{(l)}$ be a partition of $\I$
such that
\begin{enumerate}
\item $\epsilon_{i}$ is constant on each $\mathbb{I}^{(k)}$
for $1\leq k\leq l$,

\item $\e_i\neq \e_j$ if $i\in \I^{(k)}$ and $j\in \I^{(k+1)}$ for $1\le k\le l-1$,

\item if we set $i_k =\max \mathbb{I}^{(k)}$ for $1\leq k\leq l$ and $i_0=0$, we have
\begin{equation*}\mathbb{I}^{(k)}=\mathbb{I}\cap\{i_{k-1}+1,i_{k-1}+2,\dots, i_k\}.
\end{equation*}
\end{enumerate}
Put $\sigma_{\leq j}=\sigma_{1}\sigma_{2}\cdots\sigma_{j}$ for $j\in\mathbb{I}$
and $\varsigma_{i}=\sigma_{i}\sigma_{i+1}$ for $1\leq i\leq n-1$. 

For $X=E, F, K$ and $i\in I$, we define $\tau(X_{i})\in {\U_D(\e)}^e_{\rm af}$  as follows:
\begin{enumerate}
\item If $i\in I_{\mathrm{even}}\setminus\{0,n\}$ with $(\epsilon_{i},\epsilon_{i+1})=(0,0)$,
then we put
\[
\tau(E_{i})=e_{i},\quad\tau(F_{i})=f_{i},\quad\tau(K_{i})=k_{i}.
\]

\item If $i\in I_{\mathrm{odd}}\setminus\{0,n\}$ with $(\epsilon_{i},\epsilon_{i+1})=(0,1)$,
then we put
\[
\tau(E_{i})=e_{i}\sigma_{\leq i},\quad\tau(F_{i})=f_{i}\sigma_{\leq i}\varsigma_{i},\quad\tau(K_{i})=k_{i}\varsigma_{i}.
\]

\item If $i\in I_{\mathrm{even}}\setminus\{0,n\}$ with $(\epsilon_{i},\epsilon_{i+1})=(1,1)$,
then we have $\{i,i+1\}\subset\mathbb{I}^{(k)}$ for a unique $k$ and put
\[
\tau(E_{i})=e_{i}\varsigma_{i}^{i-i_{k-1}},\quad\tau(F_{i})=-f_{i}\varsigma_{i}^{i-i_{k-1}-1},\quad\tau(K_{i})=k_{i}\varsigma_{i}.
\]

\item If $i\in I_{\mathrm{odd}}\setminus\{0,n\}$ with $(\epsilon_{i},\epsilon_{i+1})=(1,0)$,
then we have $i=i_{k}$ for unique $k$ and put
\[
\tau(E_{i})=e_{i}\sigma_{\leq i}\varsigma_{i}^{i_{k}-i_{k-1}},\quad\tau(F_{i})=(-1)^{i_{k}-i_{k-1}}f_{i}\sigma_{\leq i}\varsigma_{i}^{i_{k}-i_{k-1}-1},\quad\tau(K_{i})=k_{i}\varsigma_{i}.
\]
\item If $i=0$, then we put $\tau(X_{0})$ as in the case of $\tau(X_1)$, that is,  if $\tau(E_{1})=e_{1}\sigma$ for a monomial $\sigma\in\Sigma$,
then we put $\tau(E_{0})=e_{0}\sigma$, and define $\tau(F_0)$ and $\tau(K_0)$ in the same way.

\item If $i=n$, then we put $\tau(X_{n})$ as in the case of $\tau(X_{n-1})$.
\end{enumerate}\smallskip



Let $\ms{X}$ (resp. $\ms{Y}$) be the subalgebra of $U_D(\e)^e_{\rm af}[\sigma]$
(resp. $\U_D(\e)^e_{\rm af}[\sigma]$) generated
by $E_{i}$, $F_{i}$, $K_{i}^{\pm 1}$ ($i\in I$), $K_{\pm\Lambda_{0}}$
and $\sigma_{j}$ ($j\in\mathbb{I}$) (resp. $e_{i}$, $f_{i}$, $k_{i}^{\pm 1}$,
$k_{\pm\Lambda_{0}}$ and $\sigma_{j}$).
\begin{prop}\label{prop:tau iso}
The map $\tau$ extends to an isomorphism of $\Bbbk$-algebras $\tau:\mathscr{X}\longrightarrow\mathscr{Y}$ such that 
$\tau(K_{\Lambda_{0}})=k_{\Lambda_{0}}$ and $\tau(\sigma_{j})=\sigma_{j}$ for $j\in \I$.
\end{prop}
\pf Since the relations for the type $A$ are proved in \cite{KL20},
it remains to verify that
\begin{gather*}
  \tau\left(K_{i}E_{j}K_{i}^{-1}-q^{(\bm{\alpha}_{i},\bm{\alpha}_{j})}E_{j}\right)=\tau\left(K_{i}F_{j}K_{i}^{-1}-q^{-(\bm{\alpha}_{i},\bm{\alpha}_{j})}F_{j}\right)=0,\\
  \tau\left(E_{i}F_{j}-(-1)^{p(i)p(j)}F_{j}E_{i}\right)=0,
\end{gather*}
for $\{i,j\}=\{0,1\}$, and
\begin{align*}
& \tau\left(E_{0}E_{1}E_{2}-E_{1}E_{0}E_{2}+(-1)^{\epsilon_{2}+\epsilon_{3}}(E_{1}E_{2}E_{0}-E_{0}E_{2}E_{1})+E_{2}E_{0}E_{1}-E_{2}E_{1}E_{0}\right)\\
& =\tau\left(F_{0}F_{1}F_{2}-F_{1}F_{0}F_{2}+(-1)^{\epsilon_{2}+\epsilon_{3}}(F_{1}F_{2}F_{0}-F_{0}F_{2}F_{1})+F_{2}F_{0}F_{1}-F_{2}F_{1}F_{0}\right)=0,
\end{align*}
which are straightforward to check.
\qed

\subsection{Classical limits of $\ring{\U}_D(\e)$-modules}\label{subsec:cl-lim-U_D(e)-modules}

Let $V=V_{\ov 0}\oplus V_{\ov 1}$ be a complex superspace and let $\gl(V)$ be the Lie superalgebra of linear endomorphisms of $V$.
Let $B$ be a nondegenerate supersymmetric bilinear form on $V$. 
Let $\mf{osp}(V)$ be the subalgebra of $\gl(V)$, which is defined to be the $\Z_2$-graded subspace whose homogeneous components consist of $T$ satisfying $B(T(v),w)+(-1)^{|T||v|}B(v,T(w))=0$ for $v,w\in V$, where $|T|$ and $|v|$ denote the degrees of $T$ and $v$. Similarly, $\mf{spo}(V)$ is defined with respect to a nondegenerate skew supersymmetric bilinear form on $V$.

Let $M$ and $N$ be the number of $0$'s and $1$'s in $\e$ and let $V$ be the superspace such that $\dim V_{\ov 0}=2M$ and $\dim V_{\ov 1}=2N$.
We put
\begin{equation}\label{eq:osp(e)}
 \mf{osp}(\e)=
\begin{cases}
 \mf{spo}(V) & \text{if $\e_n=1$},\\
 \mf{osp}(\Pi\, V) & \text{if $\e_n=0$},\\
\end{cases}
\end{equation}
where $\Pi$ denotes the parity change functor. 
We take the fundamental system for $\mf{spo}(V)$ (the set of simple roots for $\mf{spo}(V)$) to be the one corresponding to $\e$. For example, when $\e_n=1$, we take the fundamental system whose $\upvarepsilon \updelta$ sequence (cf.~\cite[Section 1.3]{CW}) is obtained by replacing $0$ and $1$ in $\e$ with $\updelta$ and $\upvarepsilon$, respectively.

Now we explain how to take the classical limit of a $\ring{\U}_D(\e)$-module as an $\mf{osp}(\e)$-module. 
Let $\ring{Q}=\bigoplus_{i=1}^n\Z\alpha_i=\bigoplus_{i=1}^n\Z\mb{\alpha}_i\subset P$.
Note that any $\la\in {P}$ can be written uniquely as
\begin{equation}\label{eq:lambda in P of finite type}
 \la=\ell\varpi_n + \sum_{i\in \I}l_i\de_i,
\end{equation}
for some $\ell, l_1,\dots,l_n\in\Z$. 

Let $V$ be a $\ring{\U}_D(\e)$-module with weight space decomposition $V=\bigoplus_{\la\in {P}}V_\la$ in \eqref{eq:wt space} with respect to $\langle\,k_\mu\,|\,\mu\in P\,\rangle$. 
Let $\ring{\ms{X}}$ (resp. $\ring{\ms{Y}}$) be the subalgebra of ${\ms X}$ (resp. $\ms{Y}$) generated by $E_i, F_i, K_i^{\pm 1}$ (resp. $e_i, f_i, k_i^{\pm 1}$) for $i\in I\setminus\{0\}$, and $\sigma_j$ for $j\in \I$.
We may extend $V$ to a $\ring{\ms Y}$-module by defining $\sigma_jv=(-1)^{l_j}v$ for $v\in V_\la$ with $\la\in {\rm wt}(V)$ given as in \eqref{eq:lambda in P of finite type}. 
Let $V^\tau=\{\,v^\tau\,|\,v\in V\,\}$ denote the $\mathring{\ms X}$-module obtained from $V$ by applying $\tau$ in Proposition \ref{prop:tau iso}. 
\begin{lem}\label{lem:tau pullback of wt space}
For $\la \in {P}$ and $v\in V_\la$, we have
\[
K_\beta v^\tau= (-1)^{(\ell\varpi_n|\beta)} q^{(\la|\beta)}v^\tau \quad (\beta\in \ring{Q}).
\]
\end{lem}
\pf Let $\la$ be as in \eqref{eq:lambda in P of finite type} and $\beta=\sum_{i=1}^n a_i\alpha_i=\sum_{j\in \I}b_j\de_j$. 
Let $\la^\circ=\la-\ell\varpi_n$.
Note that 
\begin{equation*}
 {\bq}(\la,\alpha_i)=
\begin{cases}
 \bq(\la^\circ,\alpha_i) & \text{if $i\neq n$},\\
 {\bq}(\ell\varpi_n,\alpha_n)\bq(\la^\circ,\alpha_n) & \text{if $i = n$}.
\end{cases}
\end{equation*}
Here ${\bq}(\ell\varpi_n,\alpha_n)={\bf v}^{-2\ell}=(-q)^{-\ell}$ and hence ${\bq}(\ell\varpi_n,\beta)=(-q)^{-a_n\ell}=(-q)^{(\ell\varpi_n|\beta)}$.
For $v\in V_\la$, we have
\begin{equation*}
\begin{split}
 K_\beta v^\tau= (\tau(K_\beta) v)^\tau &
 =\left(k_\beta\prod_{j\in \I}\sigma_j^{b_j}v\right)^\tau
 =\left({\bq}(\la,\beta)\prod_{j\in \I}\sigma_j^{b_j}v\right)^\tau\\
 &=(-q)^{(\ell\varpi_n|\beta)}\bq(\la^\circ,\beta)\prod_{j\in \I}(-1)^{\e_j l_jb_j}v^\tau\\
 &=(-1)^{(\ell\varpi_n|\beta)} q^{(\ell\varpi_n|\beta)+(\la^\circ|\beta)} v^\tau=(-1)^{(\ell\varpi_n|\beta)} q^{(\la|\beta)}v^\tau.
\end{split}
\end{equation*}
\qed\smallskip

Suppose further that $V$ is a highest weight $\ring{\U}_D(\e)$-module with highest weight $\mu\in {P}$ such that $\mu=\ell\varpi_n + \sum_{i\in \I}m_i\de_i$ for some $\ell, m_1,\dots,m_n\in\Z$.
Then we have
\begin{equation}\label{eq:dominated by lambda}
 {\rm wt}(V)\subset \ell\varpi_n + \sum_{i\in \I}\Z\de_i.
\end{equation}
By replacing $K_{n}$ (or $K_{\bm{\alpha}_n}$) and $F_n$ with $(-1)^\ell K_n$ and $(-1)^\ell F_n$, respectively, we have for $v\in V_\la$ with $\la\in {\rm wt}(V)$
\begin{equation*} 
 K_\beta v^\tau= q^{(\la|\beta)}v^\tau \quad (\beta\in \ring{Q}),
\end{equation*}
which yields a usual weight space decomposition of $V^\tau$ as an $\ring{\ms X}$-module.

Now we consider the classical limit of $V^\tau$ as a module over the subalgebra $\langle\, E_i, F_i, K_i^{\pm 1}\,|\, i\in I\setminus\{0\}\,\rangle\subset \ring{U}_D(\e)$ (cf.~\cite{Ja}).
Let ${\bf A}=\Z[q,q^{-1}]$. Let $V^\tau_{\bf A}$ be the ${\bf A}$-span of $f_{i_1}\ldots f_{i_s}v$ for $s\geq 0$ and $i_1,\ldots,i_s\in I\setminus\{0\}$, where $v$ is a highest weight vector.  
Then $V^\tau_{\bf A}$ is invariant under $E_i$, $F_i$, $K_i^{\pm 1}$ and $\{K_i\}:=\frac{K_i-K_i^{-1}}{q-q^{-1}}$ for $i\in I\setminus\{0\}$.
Let $\ov{V^\tau}=V^\tau_{{\bf A}}\otimes_{{\bf A}}\mathbb{\C}$, where $\C$ is an ${\bf A}$-module such that $f(q)\cdot z = f(1)z$ for $f(q)\in {\bf A}$ and $z\in \C$.
Then $\ov{V^\tau}$ is invariant under the $\mathbb{C}$-linear endomorphisms ${\rm E}_i$, ${\rm F}_i$ and ${\rm H}_i$ induced from $E_i$, $F_i$ and $\{K_i\}$, respectively for $i\in I\setminus\{0\}$. 

The endomorphisms ${\rm E}_i$, ${\rm F}_i$ and ${\rm H}_i$ satisfy the defining relations for the enveloping algebra $U(\mf{osp}(\e))$ in \eqref{eq:osp(e)} (cf.~\cite[Section 10.5]{Ya94}). 
Hence we have the following.

\begin{lem}\label{lem:classical limit of hw module}
Under the above hypothesis, $\ov{V^\tau}$ becomes an $U(\mf{osp}(\e))$-module with highest weight $\mu$.
\end{lem}

\subsection{Quasi $R$ matrix for $\U_D(\e)$}

Let $\U_D(\e)^+$ (resp. $\U_D(\e)^-$) be the subalgebra of $\U_D(\e)$ generated by $e_i$ (resp. $f_i$) for $i\in I$. Note that $\U_D(\e)^\pm$ is naturally graded by $Q_\pm$.
The following can be proved by the same arguments as in \cite[Sections 2.2 and 2.4]{KL20} and Proposition \ref{prop:tau iso}.

\begin{prop}\label{prop:bilinear form}
\mbox{}
\begin{itemize}
 \item[(1)] For $i\in I$, there exists a unique $\Bbbk$-linear map ${}_{i}r$ on $\U_D(\e)^-$ such that 
\begin{itemize}
\item[(i)] ${}_{i}r(1)=0$ and ${}_{i}r(f_j)=\delta_{ij}$ for $j\in I$,

\item[(ii)] ${}_{i}r(xy) = {}_{i}r(x)y + \bq(\beta,\alpha_i)x{}_{i}r(y)$ for homogeneous $x\in \U_D(\e)^-_\beta$ and $y\in \U_D(\e)^-$.
\end{itemize}
\item[(2)] There exists a unique nondegenerate symmetric bilinear form $(\, \cdot \,,\, \cdot\, )$ on $\U_D(\e)^-$ such that $(1,1)=1$ and $(f_ix,y) = (x, {}_ir(y))$ for $i\in I$ and $x,y\in \U_D(\e)^-$.
\end{itemize}
\end{prop}

For $\beta\in Q_+$, let ${\bf B}_\beta$ be a basis of $\U_D(\e)^-_{-\beta}$ and let ${\bf B}^*_{\beta}=\{\,b^*\,|\,b\in {\bf B}_\beta\,\}$ be the dual basis of ${\bf B}_\beta$ with respect to $(\,\cdot\,,\,\cdot\,)$ in Proposition \ref{prop:bilinear form}(2). Let $ - $ denote the involution on $\U_D(\e)$ as a $\C$-algebra given by $\ov{e_i}=e_i$, $\ov{f_i}=f_i$, $\ov{k_\mu}=k_{-\mu}$ ($i\in I$, $\mu\in P$) and $\ov{q}=q^{-1}$.
Then the existence of the quasi $R$ matrix $\Theta$ in $\U(\e)^+\widehat{\ot}\ \U(\e)^-:= \bigoplus_{\xi\in Q}\prod_{\xi=\mu+\nu}
\U_D(\e)^+_\mu \ot \U_D(\e)^-_\nu$ can be proved by almost the same arguments as in \cite[Theorem 4.1.2]{Lu93} (cf.~\cite{KL20}).
\begin{thm}\label{thm:quasi R matrix}\mbox{}
\begin{itemize}
\item[(1)]
There is a unique $\Theta_\beta\in \U_D(\e)^+_\beta\ot\U_D(\e)^-_{-\beta}$ for each $\beta\in Q_+$ such that $\Theta_0=1\ot 1$ and 
$\Theta\, \Delta(u) = \ov{\Delta}(u)\, \Theta$,
for $u\in \U_D(\e)$, where $\ov{\Delta}(u)=\ov{\Delta(\ov{u})}$ and $\Theta = \sum_{\beta\in Q_+}\Theta_\beta\in \U(\e)^+\widehat{\ot}\ \U(\e)^-$.  
\item[(2)] For $\beta=\sum_{i\in I}a_i \mb{\alpha}_i\in Q_+$, we have
\begin{equation*}
\Theta_\beta = (q-q^{-1})^{{\rm ht}(\beta)}\sum_{b\in {\bf B}_\beta}b^+\ot b^*,
\end{equation*}
where ${\rm ht}(\beta)=\sum_{i\in I}a_i$ and $b^+$ is obtained from $b$ by replacing $f_i$ with $e_i$ for $i\in I$.

\item[(3)] Let $\ov{\Theta}=\sum_{\beta\in Q_+}\ov{\Theta}_\beta$, where $\ov{\Theta}_\beta:= (- \ot -)(\Theta_\beta)$. Then $\ov{\Theta}{\Theta} ={\Theta}\ov{\Theta}=1$.
\end{itemize}
\end{thm}

\subsection{$\U_D(\e)$-modules $\W(x)$ and $\W^{\ot 2}(x)$}\label{subsec:W and W2}

Let
\begin{equation*}
\Z^n_+(\e)=\{\,{\bf m}=(m_1,\ldots,m_n)\,|\,\e_i=0 \Rightarrow m_i\in\Z_+,\  \e_i=1 \Rightarrow m_i\in \{0,1\}\,\}.
\end{equation*}
Let
\begin{equation*}
\begin{split}
\W &= \bigoplus_{{\bf m}\in \Z^n_+(\e)}\Bbbk|{\bf m}\rangle 
\end{split}
\end{equation*}
be the $\Bbbk$-space spanned by $\ket{\bf m}$ for ${\bf m}\in \Z^n_+(\e)$, and let 
\begin{equation*}
\begin{split}
\W^+ &= \bigoplus_{\text{$|{\bf m}|$:even}}\Bbbk|{\bf m}\rangle, \quad 
\W^- = \bigoplus_{\text{$|{\bf m}|$:odd}}\Bbbk|{\bf m}\rangle,
\end{split}
\end{equation*}
where $|{\bf m}|=m_1+\dots + m_n$. Let ${\bf 0}=(0,\dots,0)$ and $\{\,\be_i\,|\,1\le i\le n\,\}$ be the standard basis of $\Z^n$.

\begin{prop}\label{prop:osc W for e}
Suppose that $\epsilon_1=\epsilon_n=1$.
\begin{itemize}
 \item[(1)] For each $x\in\Bbbk^\times$, the following action of the generators defines a $\U_D(\e)$-module structure on $\W$, which we denote by $\W(x)$:
\begin{align*}
&e_{0} |{\bf m}\rangle = x|{\bf m}+{\bf e}_{1}+{\bf e}_{2} \rangle,\\
&f_{0} |{\bf m}\rangle = x^{-1}[m_2]|{\bf m}-{\bf e}_{1}-{\bf e}_{2} \rangle,\\
&k_{0} |{\bf m}\rangle =q_1^{m_1-1}q_2^{m_2}|{\bf m} \rangle, \\
&e_{i} |{\bf m}\rangle = [m_{i}]|{\bf m}-{\bf e}_{i}+{\bf e}_{i+1} \rangle \quad(1\leq i < n), \\
&f_{i} |{\bf m}\rangle = [m_{i+1}]|{\bf m}+{\bf e}_{i}-{\bf e}_{i+1} \rangle \quad(1\leq i < n), \\
&k_{i} |{\bf m}\rangle =q^{-m_i} _i q^{m_{i+1}} _{i+1}  |{\bf m}\rangle \quad(1\leq i < n),\\
&e_{n} |{\bf m}\rangle = [m_{n-1}]|{\bf m}-{\bf e}_{n-1}-{\bf e}_{n} \rangle,\\
&f_{n} |{\bf m}\rangle = |{\bf m}+{\bf e}_{n-1}+{\bf e}_{n} \rangle,\\
&k_{n} |{\bf m}\rangle =q_n^{1-m_n}q_{n-1}^{-m_{n-1}}|{\bf m} \rangle,\\
&k_{\La}\ket{{\bf m}} = {\bf v}^{\sum_{i\in \I}m_i}\ket{{\bf m}},
\end{align*}
for ${\bf m}=(m_1,\dots,m_n)\in \Z^n_+(\e)$, where we understand $\ket{\bf m'}=0$ unless $\ket{\bf m'}$ on the right-hand side belongs to $\Z^n_+(\e)$.

 \item[(2)] The subspace $\W^\pm$ is an irreducible submodule of $\W(x)$, which we denote by $\W^\pm(x)$. 
As a $\U_D(\e)$-module, we have $\W(x)=\W^+(x)\oplus \W^-(x)$.
\end{itemize}
\end{prop}

\begin{rem}\label{rem:wt of m}
{\rm
We have $k_{\de_i}\ket{{\bf m}} = {\bf v} q_i^{m_i}\ket{{\bf m}}$ for $i\in \I$ and $\ket{{\bf m}}\in \W(x)$, and for $\ket{{\bf m}}\in \W(x)$
\begin{equation}\label{eq:wt of m}
 {\rm wt}(\ket{{\bf m}})=\La +\sum_{i\in \I}m_i\de_i.
\end{equation}
In particular, $\W^+(x)$ and $\W^-(x)$ are highest weight $\ring{\U}_D(\e)$-modules with highest weights $\La$ and $\La+\de_n$, respectively.
}
\end{rem}

\begin{prop}\label{prop:osc W2 for e}
Let $x\in\Bbbk^\times$ be given.
If $\ep_1=\ep_n=0$, then the following action of the generators defines a $\U_D(\e)$-module structure on $\W^{\ot 2}$, which we denote by $\W^{\ot 2}(x)$:
{\small
\begin{align*}
&e_{0}\ket{\boldm}\otimes\ket{\boldmp} 
 =x\ket{\boldm+\bolde_{1}}\otimes\ket{\boldmp+\bolde_{2}}-xq_{1}^{-m_{1}}q_{2}^{-m_{2}^{\prime}}q^{-1}\ket{\boldm+\bolde_{2}}\otimes\ket{\boldmp+\bolde_{1}},\\
&f_{0}\ket{\boldm}\otimes\ket{\boldmp} \\
& =-x^{-1}q_{1}^{m_{1}^{\prime}}q_{2}^{m_{2}}q[m_{1}][m_{2}^{\prime}]\ket{\boldm-\bolde_{1}}\otimes\ket{\boldmp-\bolde_{2}}+x^{-1}[m_{1}^{\prime}][m_{2}]\ket{\boldm-\bolde_{2}}\otimes\ket{\boldmp-\bolde_{1}},\\
&k_{0}\ket{\boldm}\otimes\ket{\boldmp} 
 =q_{1}^{m_{1}+m_{1}^{\prime}}q_{2}^{m_{2}+m_{2}^{\prime}}q^{2}\ket{\boldm}\otimes\ket{\boldmp},\\
&e_{i}\ket{\boldm}\otimes\ket{\boldmp} 
 =q_{i}^{m_{i}^{\prime}}q_{i+1}^{-m_{i+1}^{\prime}}[m_{i}]\ket{\boldm-\bolde_{i}+\bolde_{i+1}}\otimes\ket{\boldmp}+[m_{i}^{\prime}]\ket{\boldm}\otimes\ket{\boldmp-\bolde_{i}+\bolde_{i+1}},\\
&f_{i}\ket{\boldm}\otimes\ket{\boldmp} 
 =[m_{i+1}]\ket{\boldm+\bolde_{i}-\bolde_{i+1}}\otimes\ket{\boldmp}+q_{i}^{-m_{i}}q_{i+1}^{m_{i+1}}[m_{i+1}^{\prime}]\ket{\boldm}\otimes\ket{\boldmp+\bolde_{i}-\bolde_{i+1}},\\
&k_{i}\ket{\boldm}\otimes\ket{\boldmp} 
 =q_{i}^{-m_{i}-m_{i}^{\prime}}q_{i+1}^{m_{i+1}+m_{i+1}^{\prime}}\ket{\boldm}\otimes\ket{\boldmp},\\
& e_{n}\ket{\boldm}\otimes\ket{\boldmp} \\
& =-q_{n-1}^{m_{n-1}^{\prime}}q_{n}^{m_{n}}q[m_{n-1}][m_{n}^{\prime}]\ket{\boldm-\bolde_{n-1}}\otimes\ket{\boldmp-\bolde_{n}}+[m_{n-1}^{\prime}][m_{n}]\ket{\boldm-\bolde_{n}}\otimes\ket{\boldmp-\bolde_{n-1}},\\
& f_{n}\ket{\boldm}\otimes\ket{\boldmp} 
 =\ket{\boldm+\bolde_{n-1}}\otimes\ket{\boldmp+\bolde_{n}}-q_{n-1}^{-m_{n-1}}q_{n}^{-m_{n}^{\prime}}q^{-1}\ket{\boldm+\bolde_{n}}\otimes\ket{\boldmp+\bolde_{n-1}},\\
&k_{n}\ket{\boldm}\otimes\ket{\boldmp} 
 =q_{n-1}^{-m_{n-1}-m_{n-1}^{\prime}}q_{n}^{-m_{n}-m_{n}^{\prime}}q^{-2}\ket{\boldm}\otimes\ket{\boldmp},\\
& k_{\La}\ket{\boldm}\otimes\ket{\boldmp} = {\bf v}^{\sum_{i\in \I}(m_i+m'_i)}\ket{\boldm}\otimes\ket{\boldmp}, 
\end{align*}}
for ${\bf m}=(m_1,\dots,m_n), {\bf m}'=(m'_1,\dots,m'_n)\in \Z^n_+(\e)$.
\end{prop}

\begin{rem}{\rm
As in Remark \ref{rem:wt of m}, we have 
$$
k_{\de_i}\ket{\boldm}\otimes\ket{\boldmp} = {\bf v}^2 q_i^{m_i+m'_i}\ket{\boldm}\otimes\ket{\boldmp},
$$ 
for $i\in \I$ and $\ket{\boldm}\otimes\ket{\boldmp}\in \W^{\ot 2}(x)$, and hence
\begin{equation}\label{eq:wt of m tensor m'}
 {\rm wt}(\ket{\boldm}\otimes\ket{\boldmp})=2\La +\sum_{i\in \I}(m_i+m'_i)\de_i.
\end{equation}
}
\end{rem}

\section{Truncation functors}\label{sec:truncation}

We introduce two kinds of monoidal or tensor functors which send representations of generalized quantum group of affine type $D$ or quantum affine superalgebras of type $D$ to those of quantum affine algebras of type $C_m^{(1)}$ and $D_{m+1}^{(1)}$ (cf.~\cite{KY} for type $A$). They play a crucial role in later sections.

\subsection{Truncation of type $\mf{c}$}\label{subsec:truncation c}

We assume the following notations:

\begin{itemize}
\item[$\bullet$] $\bm{\e}=(\e_1,\dots,\e_n)=(1,0,1\dots,0,1)$ with $n=2m+1$ ($m\ge 2$),

\item[$\bullet$] $\U_{\bm{\e}}=\U_D(\bm{\e})$ with $I_{\bm\e}=I$,

\item[$\bullet$] $\un{\bm{\e}}=(0^m)$, $\ov{\bm\e}=(1^{m+1})$,

 \item[$\bullet$] $I_{\un{\bm\e}}=\{0,1,\dots,m\}$, $I_{\ov{\bm\e}}=\{0,1,\dots,m+1\}$,

\item[$\bullet$] $\U_{\un{\bm\e}}=U_q(C_m^{(1)})=\langle\, e_i,f_i, k_i^{\pm 1}\,|\, i\in I_{\un{\bm\e}}\, \rangle$, 

\item[$\bullet$] $\U_{\ov{\bm\e}}=U_{\td{q}}(D_{m+1}^{(1)})=\langle\, e_i,f_i, k_i^{\pm 1}\,|\, i\in I_{\ov{\bm\e}}\, \rangle$,

\end{itemize}
where $U_q(X_N^{(1)})$ denotes the quantum affine algebra of type $X_N^{(1)}$ over $\Bbbk$. We remark that $\U_{\e}$ ($\e=\un{\bm \e}, \ov{\bm\e}$) should not be confused with $\U_D(\e)$.

Let us define 
$\hat{e}_i, \hat{f}_i, \hat{k}_i\in \U_{\bm\e}$ for $i\in I_{\un{\bm\e}}$, and 
$\check{e}_j, \check{f}_j, \check{k}_j \in \U_{\bm\e}$ for $j\in I_{\ov{\bm\e}}$ as follows:
\begin{equation}\label{eq:hat{e} for c}
\begin{split}
 \hat{e}_i&=\frac1{[2]^{\delta_{i0}+\delta_{i{m}}}}[e_{2i},e_{2i+1}]_{c_i},\\
 \hat{f_i}&=\frac1{[2]^{\delta_{i0}+\delta_{i{m}}}}[f_{2i+1},f_{2i}]_{c_i^{-1}},\\
 \hat{k}_i&=k_{2i}k_{2i+1}, 
\end{split}
\end{equation}
where $c_i=-q^{2\eta}$ for $i=0,m$, and $c_i=-q^{\eta}$ otherwise $(\eta=\pm 1)$, and
{\allowdisplaybreaks
\begin{equation}\label{eq:check{e} for c}
\begin{split}
 \check{e}_j&=
\begin{cases}
 [e_0,e_2]_{d} & \text{if $j=0$},\\ 
 [e_{2j-1},e_{2j}]_{d} & \text{if $1\le j\le m$},\\
 [e_{2m-1},e_{2m+1}]_{d} & \text{if $j=m+1$},
\end{cases}\\
 \check{f}_j&=
\begin{cases}
 [f_2,f_0]_{d^{-1}} & \text{if $j=0$},\\ 
 [f_{2j},f_{2j-1}]_{d^{-1}} & \text{if $1\le j\le m$},\\
 [f_{2m+1},f_{2m-1}]_{d^{-1}} & \text{if $j=m+1$},
\end{cases}\\
 \check{k}_j&=
\begin{cases}
 k_0k_2 & \text{if $j=0$},\\ 
 k_{2j-1}k_{2j} & \text{if $1\le j\le m$},\\
 k_{2m-1}k_{2m+1} & \text{if $j=m+1$},
\end{cases}
\end{split}
\end{equation}}
where $d=q^{\pm 1}$. Here $[X,Y]_{t}=XY-tYX$ for $t\in \Bbbk^\times$.

\begin{prop}\label{prop:reduction}\mbox{}
\begin{itemize}
\item[(1)] There exists a homomorphism of $\Bbbk$-algebras $\phi_{\un{\bm\e}} : \U_{\un{{\bm \e}}} \longrightarrow \U_{\bm \e}$ such that $\phi_{\un{\bm\e}}(e_i)=\hat{e}_i$, $
\phi_{\un{\bm\e}}(f_i)=\hat{f}_i$, and $\phi_{\un{\bm\e}}(k_i)=\hat{k}_i$  for $i\in I_{\un{\bm\e}}$.

\item[(2)] There exists a homomorphism of $\Bbbk$-algebras $\phi_{\ov{\bm\e}} : \U_{\ov{{\bm \e}}} \longrightarrow \U_{\bm \e}$ such that $\phi_{\ov{\bm\e}}(e_j)=\check{e}_j$, $
\phi_{\ov{\bm\e}}(f_j)=\check{f}_j$, and $
\phi_{\ov{\bm\e}}(k_j)=\check{k}_j$ for $j\in I_{\ov{\bm\e}}$.
\end{itemize}
\end{prop}
\pf It is straightforward to check that the image of the generators satisfy the defining relations for $\U_{\un{\bm \e}}$ and $\U_{\ov{\bm \e}}$, except the Serre relations. Since the type $A$ cases are proved in \cite[Theorem~4.3]{KY},  we only need to verify the following two identities:
\begin{align*}
& \widehat{e}^{2}_0 \widehat{e}_1 - (q^2+q^{-2})\widehat{e}_0 \widehat{e}_1 \widehat{e}_0 + \widehat{e}_1 \widehat{e}^{2}_0 =0, \\
& \widehat{e}^{3}_1 \widehat{e}_0 -(q^2 +1+q^{-2})\widehat{e}^{2}_1 \widehat{e}_0 \widehat{e}_1 +(q^2 +1+q^{-2})\widehat{e}_1 \widehat{e}_0 \widehat{e}^{2}_1 -\widehat{e}_0 \widehat{e}^{3}_1 =0.
\end{align*}
The proof of the first one is given in Appendix~\ref{sec:app-pf-Serre-rel}, and we leave the other one to the reader.
\qed\smallskip

Let
\begin{itemize}
 \item[$\bullet$] $\I_{\bm\e}=\I$, $\I_{\un{\bm\e}}=\{2,4,\dots,2m\}$, $\I_{\ov{\bm\e}}=\{1,3,\dots,2m+1\}$,

 \item[$\bullet$] $P_{\e}= \Z\La \oplus \bigoplus_{a\in \I_\e}\Z\de_a$, 
 $P_{\ge 0,\e}=\Z\La \oplus \bigoplus_{a\in \I_\e}\Z_+\de_a$\ \ for $\e={\bm\e},\un{\bm\e},\ov{\bm\e}$,
 
 \item[$\bullet$] $\{\,\alpha_i^\e\,|\,i\in I_\e\,\}$: the simple roots for $\U_\e$\ \ for $\e={\bm\e},\un{\bm\e},\ov{\bm\e}$, where $\alpha_i^{\bm\e}=\alpha_i$ in Section \ref{subsec:def}.
\end{itemize}
By Proposition \ref{prop:reduction}, we may regard $\{\,\alpha_i^\e\,|\,i\in I_\e\,\}\subset P_{\bm\e}$ ($\e=\un{\bm\e},\ov{\bm\e}$) by identifying
\begin{equation*}
\begin{split}
\alpha^{\un{\bm\e}}_i=
\begin{cases}
 2\de_2 & \text{if $i=0$},\\
 -\de_{2i}+\de_{2i+2} & \text{if $1\le i<m$},\\
 -2\de_{2m} & \text{if $i=m$},\\ 
\end{cases}\quad
\alpha^{\ov{\bm\e}}_i=
\begin{cases}
 \de_1+\de_3 & \text{if $i=0$},\\
 -\de_{2i-1}+\de_{2i+1} & \text{if $1\le i\le m$},\\
 -\de_{2m-1}-\de_{2m+1} & \text{if $i=m+1$}.
\end{cases}
\end{split}
\end{equation*}

For $\e={\bm\e},\un{\bm\e},\ov{\bm\e}$,
let $\mc{C}(\e)$ be the category of $\U_{\e}$-modules $V$ such that 

\begin{itemize}
 \item[(1)] $V$ has a weight space decomposition in the sense that  $V=\bigoplus_{\la\in P_\e}V_\la$ with
\begin{equation}\label{eq:weight space}
 V_\la=\{\,v\,|\,k_\mu v = {\bq}(\la,\mu)v\ (\mu\in P_\e)\,\},
\end{equation}
 
 \item[(2)] ${\rm wt}(V)\subset P_{\ge 0,\e}$.
 
\end{itemize}
For $\ell\in\Z$, we define $\mc{C}^\ell(\e)$ to be the subcategory of $\mc{C}(\e)$ such that 
\begin{equation}
 {\rm wt}(V)\subset \ell\La + \sum_{a\in \I_\e}\Z_+\de_a.
\end{equation}
Then we have $\mc{C}(\e)=\bigoplus_{\ell\in \mathbb{Z}}\mc{C}^\ell(\e)$. 
Note that $\mc{C}(\e)$ is closed under taking submodules, quotients and tensor products, while $\mc{C}^\ell(\e)$ is closed under taking submodules and quotients. 
For example, $\W(x)$ ($x\in \Bbbk^\times$) belongs to $\mc{C}(\e)$ by \eqref{eq:wt of m} and hence so does $\W(x)^{\otimes \ell}$ for $\ell\ge 1$.

Let $\e=\un{\bm\e},\ov{\bm\e}$ be given.
For a $\U_{\bm\e}$-module $V$ in $\mc{C}({\bm\e})$, define $\mf{tr}^{\bm\e}_{\e}(V)$, or simply $\mf{tr}_{\e}(V)$, as follows: 
\begin{equation}\label{eq:truncation-obj}
\mf{tr}_{\e}(V) = 
\bigoplus_{\substack{\mu\in {\rm wt}(V) \\ ({\rm pr}(\mu)|\de_{a})=0\ (a\, \not\in\, \I_\e)}}V_\mu,
\end{equation}
where ${\rm pr}: P_{\ge 0,{\bm\e}} \longrightarrow \bigoplus_{i\in \I_{\bm\e}}\Z\de_i$ denotes the canonical projection.
For $V, W\in \mc{C}({\bm\e})$ and $f \in {\rm Hom}_{\U({\bm\e})}(V,W)$, we have a $\Bbbk$-linear map
\begin{equation}\label{eq:truncation-mor}
\xymatrixcolsep{2pc}\xymatrixrowsep{3pc}
\xymatrix{
\mf{tr}_\e(f) : \mf{tr}_\e(V) \ \ar@{->}[r] &  \mf{tr}_\e(W)
}
\end{equation}
given by $\mf{tr}_\e(f)(v)=f(v)$ for $v\in \mf{tr}_\e(V)$. 

\begin{prop}\label{prop:truncation}
Under the above hypothesis,
\begin{itemize}
\item[(1)] $\mf{tr}_\e(V)$ is invariant under the action of $\phi_\e(\U_\e)$, and hence a $\U_\e$-module in $\mc{C}(\e)$,

\item[(2)] $\mf{tr}_\e(f) : \mf{tr}_\e(V) \longrightarrow \mf{tr}_\e(W)$ is $\mathcal{U}(\e)$-linear,

\item[(3)] $\mf{tr}_{\e}(V\ot W)$ is naturally isomorphic to $\mf{tr}_\e(V)\otimes \mf{tr}_\e(W)$ as a $\mathcal{U}(\e)$-module.

\end{itemize}
Hence \eqref{eq:truncation-obj} and \eqref{eq:truncation-mor} define a well-defined functor 
\begin{equation*}\label{eq:tr}
\xymatrixcolsep{2pc}\xymatrixrowsep{3pc}
\xymatrix{
\mf{tr}_\e : \mc{C}({\bm\e})  \ar@{->}[r] &  \mc{C}(\e)},
\end{equation*}
which is exact and monoidal.
\end{prop}
\pf It can be checked by the same arguments as in \cite[Propositions 4.4]{KY} due to the condition ${\rm wt}(V)\subset P_{\ge 0,{\bm\e}}$.
\qed\smallskip

\begin{rem}{\rm
The functor $\mathfrak{tr}_\epsilon$ does no depend on the choice of $\eta$ and $d$ in \eqref{eq:hat{e} for c} and \eqref{eq:check{e} for c}, and the same holds in the next subsection.
}
\end{rem}

\subsection{Truncation of type $\mf{d}$}\label{subsec:truncation d}

Let us introduce another type of truncation.
Let
\begin{itemize}
\item[$\bullet$] $\bm{\e}'=(\e_1,\dots,\e_n)=(0,1,0\dots,1,0)$ with $n=2m+1$ ($m\ge 2$),

\item[$\bullet$] $\U_{\bm{\e}'}=\U_D(\bm{\e}')$ with $I_{\bm{\e}'}=I$,

\item[$\bullet$] $\un{\bm{\e}}'=(0^{m+1})$, $\ov{\bm\e}'=(1^m)$,

 \item[$\bullet$] $I_{\un{\bm\e}'}=\{0,1,\dots,m+1\}$, $I_{\ov{\bm\e}'}=\{0,1,\dots,m\}$,

\item[$\bullet$] $\U_{\un{\bm\e}'}=U_{q}(D_{m+1}^{(1)})=\langle\, e_i,f_i, k_i^{\pm 1}\,|\, i\in I_{\un{\bm\e}'}\, \rangle$, 

\item[$\bullet$] $\U_{\ov{\bm\e}'}=U_{\td{q}}(C_{m}^{(1)})=\langle\, e_i,f_i, k_i^{\pm 1}\,|\, i\in I_{\ov{\bm\e}'}\, \rangle$.

\end{itemize}
In this case, we have an analogue of Proposition \ref{prop:reduction} as follows. 
First, let 
\begin{equation}\label{eq:hat{e} for d}
\begin{split}
 \hat{e}_i&=\left(-\frac{1}{[2]}\right)^{\delta_{i0}+\delta_{im}}[e_{2i+1},e_{2i}]_{c_{i}},\\
 \hat{f_i}&=\left(-\frac{1}{[2]}\right)^{\delta_{i0}+\delta_{im}}[f_{2i},f_{2i+1}]_{c_{i}^{-1}},\\
 \hat{k}_i&=k_{2i}k_{2i+1}, 
\end{split}
\end{equation}
for $i\in I_{\ov{\bm\e}'}$, where $c_{i}=-q^{2\eta}$ for $i=0,\,m$ and $c_i=q^{\eta}$ otherwise ($\eta=\pm 1$).
Similarly, let 
{\allowdisplaybreaks
\begin{equation}\label{eq:check{e} for d}
\begin{split}
 \check{e}_j&=
\begin{cases}
 [e_2,e_0]_{d} & \text{if $j=0$},\\ 
 [e_{2j},e_{2j-1}]_{d} & \text{if $1\le j\le m$},\\
 [e_{2m+1},e_{2m-1}]_{d} & \text{if $j=m+1$},
\end{cases}\\
 \check{f}_j&=
\begin{cases}
 [f_0,f_2]_{d^{-1}} & \text{if $j=0$},\\ 
 [f_{2j-1},f_{2j}]_{d^{-1}} & \text{if $1\le j\le m$},\\
 [f_{2m-1},f_{2m+1}]_{d^{-1}} & \text{if $j=m+1$},
\end{cases}\\
 \check{k}_j&=
\begin{cases}
 k_0k_2 & \text{if $j=0$},\\ 
 k_{2j-1}k_{2j} & \text{if $1\le j\le m$},\\
 k_{2m-1}k_{2m+1} & \text{if $j=m+1$},
\end{cases}
\end{split}
\end{equation}}
for $j\in I_{\un{\bm\e}'}$, where $d=-q^{\pm1}$. Note that \eqref{eq:hat{e} for d} and \eqref{eq:check{e} for d} can be obtained from \eqref{eq:hat{e} for c} and \eqref{eq:check{e} for c} (up to scalar multiplication) by replacing $q$ with $\td{q}$.

\begin{prop}\label{prop:reduction'}\mbox{}
\begin{itemize}
\item[(1)] There exists a homomorphism of $\Bbbk$-algebras $\phi_{\ov{\bm\e}'} : \U_{\ov{{\bm \e}}'} \longrightarrow \U_{{\bm \e}'}$ such that $\phi_{\ov{\bm\e}'}(e_i)=\hat{e}_i$, $
\phi_{\ov{\bm\e}}(f_i)=\hat{f}_i$, and $\phi_{\ov{\bm\e}'}(k_i)=\hat{k}_i$  for $i\in I_{\ov{\bm\e}'}$.

\item[(2)] There exists a homomorphism of $\Bbbk$-algebras $\phi_{\un{\bm\e}'} : \U_{\un{{\bm \e}}'} \longrightarrow \U_{{\bm \e}'}$ such that $\phi_{\un{\bm\e}'}(e_j)=\check{e}_j$, $
\phi_{\un{\bm\e}'}(f_j)=\check{f}_j$, and $
\phi_{\un{\bm\e}'}(k_j)=\check{k}_j$ for $j\in I_{\un{\bm\e}'}$.
\end{itemize}
\end{prop}

Let \begin{itemize}
\item[$\bullet$] $\I_{\bm\e'}=\I$, $\I_{\un{\bm\e}'}=\{1,3,\dots,2m+1\}$, $\I_{\ov{\bm\e}'}=\{2,4,\dots,2m\}$, 

\item[$\bullet$] $P_{\e}= \Z\La \oplus \bigoplus_{a\in \I_\e}\Z\de_a$, 
 $P_{\ge 0,\e}=\Z\La \oplus \bigoplus_{a\in \I_\e}\Z_+\de_a$\ \ for $\e={\bm\e}',\un{\bm\e}',\ov{\bm\e}'$,
 
\item[$\bullet$] $\{\,\alpha_i^\e\,|\,i\in I_\e\,\}$: the simple roots for $\U_\e$\ \ for $\e={\bm\e}',\un{\bm\e}',\ov{\bm\e}'$, 
where $\alpha_i^{\bm\e'}=\alpha_i$ in Section \ref{subsec:def}, $\alpha_i^{\un{\bm\e}'} = \alpha_i^{\ov{\bm\e}}$ for $i\in I_{\un{\bm\e}'}$ and $\alpha_j^{\ov{\bm\e}'} = \alpha_j^{\un{\bm\e}}$ for $j\in I_{\ov{\bm\e}'}$ in Section \ref{subsec:truncation c}.

\end{itemize}

We define the category $\mc{C}(\e)$ and the functor $\mf{tr}^{\bm\e'}_\e$ or simply $\mf{tr}_\e$
\begin{equation}\label{eq:tr-2}
\xymatrixcolsep{2pc}\xymatrixrowsep{3pc}
\xymatrix{
\mf{tr}_\e : \mc{C}({\bm\e}')  \ar@{->}[r] &  \mc{C}(\e)} \quad (\e=\un{\bm \e}', \ov{\bm \e}')
\end{equation}
in the same way as in $\mf{tr}_\e$ for $\e=\un{\bm \e}, \ov{\bm \e}$.
Then Proposition \ref{prop:truncation} also holds in this case. 

\subsection{Truncations of $\W(x)$ and $\W^{\ot 2}(x)$}
Let $x\in \Bbbk^\times$ be given.

First, let $\W_{\bm\e}(x)$ denote the $\U_{\bm\e}$-module $\W(x)$ in Proposition \ref{prop:osc W for e} with respect to ${\bm\e}$ (cf. Section \ref{subsec:truncation c}), and let
\begin{equation}\label{eq:W for e}
 \W_{\e}(x)=\mf{tr}_{\e}(\W_{\bm\e}(x)),\quad \W_{\e}^\pm(x)=\mf{tr}_{\e}(\W^\pm_{\bm\e}(x))\quad (\e=\un{\bm\e},\ov{\bm\e}). 
\end{equation}
Note that $\W_{\e}(x)$ should not be confused with $\W(x)$ as a $\U_D(\e)$-module. 

We remark that $\W^\pm_{\un{\bm\e}}(x)$ is isomorphic to the $q$-oscillator representation of $U_q(C_m^{(1)})$ introduced in \cite[Proposition 3]{KO}, and $\W^\pm_{\ov{\bm\e}}(x)$ is isomorphic to a spin representation of $U_{\td{q}}(D_{m+1}^{(1)})$.\smallskip

Next, let $\W^{\ot 2}_{\bm\e'}(x)$ denote the $\U_{\bm\e'}$-module $\W^{\ot 2}(x)$ in Proposition \ref{prop:osc W2 for e} with respect to ${\bm\e'}$ (cf. Section \ref{subsec:truncation d}), and let
\begin{equation}\label{eq:W for e'}
 \W^{\ot 2}_{\e}(x)=\mf{tr}_{\e}(\W^{\ot 2}_{\bm\e'}(x))\quad (\e=\un{\bm\e}',\ov{\bm\e}').
\end{equation}
In this case, $\W^{\ot 2}_{\ov{\bm\e}'}(x)$ is a finite-dimensional representation of $U_{\td{q}}(C_m^{(1)})$, while $\W^{\ot 2}_{\un{\bm\e}'}(x)$ is an infinite-dimensional representation of $U_{q}(D_{m+1}^{(1)})$.

\section{Oscillator representations of $\mf{sp}_{2n}$ and $\mf{so}_{2n}$}\label{sec:osc rep}

The spaces $\W(x)$ and $\W^{\ot 2}(x)$ in Section \ref{subsec:W and W2}, when restricted to representations of $U_q(C_m)$ and $U_q(D_{m+1})$, are $q$-analogues of infinite-dimensional representations of $\mf{sp}_{2m}$ and $\mf{so}_{2m+2}$ so-called oscillator representations, respectively. Each of them belongs to a nice semisimple tensor category generated by a family of irreducible representations, which arise naturally from a viewpoint of Howe duality. 
In this section, we give a brief review of these representations. 

Let us assume the following notations:
\begin{itemize}
\item[$\bullet$] $\mf{g}=\mf{sp}_{2n}$, $\mf{so}_{2n}$ ($n\ge 4$).

\item[$\bullet$] $\texttt{P} = \bigoplus_{i=1}^n\Z\upvarpi_i$: the weight lattice of $\mf{g}$, where $\upvarpi_i$ is the $i$-th fundamental weight, 
 
\item[$\bullet$] $\ude_1=-\upvarpi_1$, $\ude_i=\begin{cases}
  \upvarpi_{i-1}-\upvarpi_{i} & \text{ if }\mathfrak{g}=\mathfrak{sp}_{2n}\\
  \upvarpi_{i-1}-\upvarpi_{i}-\delta_{i,n-1}\upvarpi_n & \text{ if }\mathfrak{g}=\mathfrak{so}_{2n}
\end{cases}\quad$ for $1<i\le n$,

\item[$\bullet$] $\upalpha_i$: the simple root for $\mf{g}$, where $\upalpha_i=\ude_{i+1}-\ude_{i}$ for $1\le i< n$, and $\upalpha_n=-2\ude_n$ (resp. $-\ude_{n}-\ude_{n-1}$) if $\mf{g}=\mf{sp}_{2n}$ (resp. $\mf{g}=\mf{so}_{2n}$).
\end{itemize}
Note that
\begin{equation*}
\begin{split}
\upvarpi_i&=
\upvarpi_n + \ude_{n}+\dots+\ude_{i+1}\quad \text{if $1\le i<n$},
\quad (\mf{g}=\mf{sp}_{2n}),\\ 
\upvarpi_i&=
\begin{cases}
\upvarpi_n+\ude_n  & \text{if $i=n-1$},\\
2\upvarpi_n + \ude_{n}+\dots+\ude_{i+1} & \text{if $1\le i<n-1$},\\
\end{cases}\quad (\mf{g}=\mf{so}_{2n}), 
\end{split}
\end{equation*} 
and hence $\ude_1+\dots+\ude_n=-r\upvarpi_n$ where $r=1$ for $\mf{sp}_{2n}$ and $r=2$ for $\mf{so}_{2n}$. For $\upvarpi\in \texttt{P}$, let $V(\upvarpi)$ be the irreducible representation of $\mf{g}$ with highest weight $\upvarpi$.

\subsection{Howe duality for $(\mf{sp}_{2n},O_\ell)$ and $(\mf{so}_{2n},O_\ell)$}\label{subsec:Howe duality}

For a positive integer $\ell$, let $O_\ell$ be the complex orthogonal group and   
\begin{equation*}
\mc{P}(O_\ell) =\{\,\la=(\la_1,\dots,\la_\ell)\in \Z_+^\ell\,|\, \la_1\ge \dots\ge \la_\ell,\ \la'_1+\la'_2\le \ell\,\},
\end{equation*}
where $\la'=(\la'_1,\la'_2,\dots)$ is the conjugate of $\la$ as a partition. 
Then $\mc{P}(O_\ell)$ parametrizes the finite-dimensional irreducible representations of $O_\ell$, say $V_{O_\ell}(\la)$ (cf.~\cite{BT}). 
For $\la\in \mc{P}(O_\ell)$, we define the weight $\upvarpi_{\la,\mf{g}}\in {\texttt{P}}$  by
\begin{equation}\label{eq:Lambda highest weight}
\begin{split}
\upvarpi_{\la,\mf{sp}_{2n}} 
&=
\begin{cases}
-\ell \upvarpi_n +\sum_{i=1}^n\la_i\ude_{n-i+1} & \text{if $\ell(\la)\le n$},\\
0 & \text{otherwise},
\end{cases}\\
\upvarpi_{\la,\mf{so}_{2n}}
&= 
\begin{cases}
\ell \upvarpi_n +\sum_{i=1}^n\la'_i\ude_{n-i+1} & \text{if $\ell(\la')\le n$},\\
0 & \text{otherwise}.
\end{cases}
\end{split}
\end{equation}
In particular, when $\ell=1$, we have $\mc{P}(O_1)=\{\,(0), (1)\,\}$ and 
\begin{equation}\label{eq:hw for spin}
\upvarpi_{(0),\mf{g}} = 
\begin{cases}
 -\upvarpi_n, & \text{if $\mf{g}=\mf{sp}_{2n}$}, \\
 \upvarpi_n, & \text{if $\mf{g}=\mf{so}_{2n}$},
\end{cases}
\quad\quad
\upvarpi_{(1),\mf{g}} = 
\begin{cases}
 -\upvarpi_n+\ude_n, & \text{if $\mf{g}=\mf{sp}_{2n}$}, \\
 \upvarpi_n+\ude_n=\upvarpi_{n-1}, & \text{if $\mf{g}=\mf{so}_{2n}$}.
\end{cases}
\end{equation}
Let
\begin{equation*}
 \mc{P}(O_\ell)_{\mf{g}}=\{\,\la\,|\,\la\in \mc{P}(O_\ell),\ \upvarpi_{\la,\mf{g}} \neq 0 \,\},
\end{equation*}
and  
\begin{equation*}
 V(\la,{\mf{g}}) = V(\upvarpi_{\la,\mf{g}} )
\end{equation*}
denote the irreducible $\mf{g}$-module with highest weight $\upvarpi_{\la,\mf{g}}$ for $\la\in \mc{P}(O_\ell)_{\mf g}$. 
For later use, let us also write $V_{\mf{g}}^{+}=V((0),{\mf{g}})$ and $V_{\mf{g}}^{-}=V((1),{\mf{g}})$.

We remark that $V(\la,\mf{sp}_{2n})$ is an infinite-dimensional representation called an {\em oscillator representation}, while $V(\la,\mf{so}_{2n})$ is finite-dimensional, and $ V(\la,{\mf{g}})$ for $\la\in \mc{P}(O_\ell)_{\mf{g}}$ satisfy the following duality (see for example, ~\cite[Theorem 3.2(a)]{HTW} and references therein, or \cite[Theorems 3.2 and 5.3]{W99} given in terms of the Lie algebras of infinite rank).

\begin{thm}\label{thm:Howe duality}
Let 
\begin{equation*}
W_{\mf g}=
\begin{cases}
S(\C^n) & \text{if ${\mf g}=\mf{sp}_{2n}$},\\
\Lambda(\C^n) & \text{if ${\mf g}=\mf{so}_{2n}$}. 
\end{cases}
\end{equation*}
Then there exists a $(\mf{g},O_\ell)$-action on $W_{\mf{g}}^{\ot \ell}$ $(\ell\ge 1)$ which gives the following multiplicity-free decomposition as a $(\mf{g},O_\ell)$-module:
\begin{equation*}
W_{\mf{g}}^{\ot \ell} = 
\bigoplus_{\la\in \mc{P}(O_\ell)_{\mf{g}}} V(\la,{\mf{g}})\ot V_{O_\ell}(\la).
\end{equation*}
\end{thm}\smallskip

In particular, we have $W_{\mf{g}}=V_{\mf{g}}^+\oplus V_{\mf{g}}^-$.

\subsection{Howe duality for $(\mf{so}_{2n},{Sp}_{2\ell})$ and $(\mf{sp}_{2n},{Sp}_{2\ell})$}\label{subsec:Howe duality-2}

For a positive integer $\ell$, let $Sp_{2\ell}$ be the complex symplectic group and   
\begin{equation*}
\mc{P}(Sp_{2\ell}) =\{\,\la=(\la_1,\dots,\la_\ell)\in \Z_+^\ell\,|\, \la_1\ge \dots\ge \la_\ell \,\},
\end{equation*}
which parametrizes the finite-dimensional irreducible representations of $Sp_{2\ell}$, say $V_{Sp_{2\ell}}(\la)$ (cf.~\cite{BT}). 
For $\la\in \mc{P}(Sp_{2\ell})$, we define the weight 
$\upvarpi'_{\la,\mf{g}} \in \texttt{P}$ by
\begin{equation}\label{eq:Lambda highest weight-2}
\begin{split}
\upvarpi'_{\la,\mf{so}_{2n}}
&=
\begin{cases}
-\ell \upvarpi_n +\sum_{i=1}^n\la_i\ude_{n-i+1} & \text{if $\ell(\la)\le n$},\\
0 & \text{otherwise},
\end{cases}\\
\upvarpi'_{\la,\mf{sp}_{2n}}
&= 
\begin{cases}
\ell \upvarpi_n +\sum_{i=1}^n\la'_i\ude_{n-i+1} & \text{if $\ell(\la')\le n$},\\
0 & \text{otherwise}.
\end{cases}
\end{split}
\end{equation}
Let
\begin{equation*}
 \mc{P}(Sp_{2\ell})_{\mf{g}}=\{\,\la\,|\,\la\in \mc{P}(Sp_{2\ell}),\  \upvarpi'_{\la,\mf{g}} \neq 0 \,\},
\end{equation*}
and    
\begin{equation*}
 V'(\la,{\mf{g}}) = V(\upvarpi'_{\la,\mf{g}})
\end{equation*}
denote the irreducible $\mf{g}$-module with highest weight $\upvarpi'_{\la,\mf{g}}$ for $\la\in \mc{P}(Sp_{2\ell})_{\mf g}$. 

In this case, $V'(\la,\mf{so}_{2n})$ is infinite-dimensional also called an oscillator representation, while $V'(\la,\mf{sp}_{2n})$ is finite-dimensional, and $ V'(\la,{\mf{g}})$ for $\la\in \mc{P}(Sp_{2\ell})_{\mf{g}}$ satisfy the following duality (see for example, ~\cite[Theorem 3.2(b)]{HTW}, \cite[Theorems 3.4 and 5.2]{W99}).

\begin{thm}\label{thm:Howe duality-2}
Let
\begin{equation*}
W'_{\mf g}=
\begin{cases}
S(\C^n)^{\ot 2}  & \text{if ${\mf g}=\mf{so}_{2n}$},\\
\Lambda(\C^n)^{\ot 2} & \text{if ${\mf g}=\mf{sp}_{2n}$}. 
\end{cases}
\end{equation*}
Then there exists a $(\mf{g},Sp_{2\ell})$-action on $(W'_{\mf g})^{\ot \ell}$ $(\ell\ge 1)$ which gives the following multiplicity-free decomposition as a $(\mf{g},Sp_{2\ell})$-module:
{\rm \begin{equation*}
(W'_{\mf g})^{\ot \ell} = 
\bigoplus_{\la\in \mc{P}(Sp_{2\ell})_{\mf{g}}}V'(\la,{\mf{g}}) \ot V_{Sp_{2\ell}}(\la).
\end{equation*}}
\end{thm}\smallskip

In particular, we have 
\begin{equation*}
\begin{split}
W'_{\mf{g}}
 &=\bigoplus_{(l)\in \mc{P}(Sp_{2})_{\mf g}}V'((l),{\mf{g}})^{\oplus (l+1)}\\
\end{split}
\end{equation*}
where $\mc{P}(Sp_{2})_{\mf{so}_{2n}}=\{\,(l)\,|\,l\in\Z_+\,\}$ and $\mc{P}(Sp_{2})_{\mf{sp}_{2n}}=\{\,(0),(1),\dots,(n)\,\}$.

\subsection{Categories ${O}^{\mf c}_{\rm osc}$ and ${O}^{\mf d}_{\rm osc}$}

Let $(\mf{g},G_\ell)$ denote one of the pairs in Sections \ref{subsec:Howe duality} and \ref{subsec:Howe duality-2}, where $\mf{g}=\mf{so}_{2n}, \mf{sp}_{2n}$ and $G=O, Sp$. Let $L^\la$ denote either $V(\la,\mf{g})$ or $V'(\la,\mf{g})$ for $\la\in \mc{P}(G_\ell)_{\mf g}$.

Suppose that $\mu\in \mc{P}(G_{\ell})_{\mf{g}}$ and $\nu\in \mc{P}(G_{\ell'})_{\mf{g}}$ are given. 
By Theorems \ref{thm:Howe duality} and \ref{thm:Howe duality-2}, $L^\mu\ot L^\nu \subset W_{\mf{g}}^{\ot{(\ell+\ell')}}$ or $ {W'_{\mf{g}}}^{\ot{(\ell+\ell')}}$ is semisimple, where 
\begin{equation}\label{eq:LR for osc}
L^\mu\ot L^\nu = \bigoplus_{\la\in \mc{P}(G_{\ell+\ell'})_{\mf{g}}}(L^\la)^{\oplus c_{\mu\nu}^\la(\mf{g},G)},
\end{equation}  
with the multiplicity $c_{\mu\nu}^\la(\mf{g},G)$ ($G=O,Sp$) given by
\begin{equation}\label{eq:LR}
c_{\mu\nu}^\la(\mf{g},G)=\dim{\rm Hom}_{G_\ell\times G_{\ell'}}\left( V_{G_\ell}(\mu)\ot V_{G_{\ell'}}(\nu), V_{G_{\ell+\ell'}}(\la) \right),
\end{equation}
which follows from considering reciprocity pairs in \cite[Table II and (3.16)]{HTW}.
In particular, it follows from \eqref{eq:LR} that 
\begin{equation}\label{eq:LR equal}
c_{\mu\nu}^\la(\mf{sp}_{2n},G)=c_{\mu\nu}^\la(\mf{so}_{2n},G), 
\end{equation}
whenever they are nonzero. Note that $c_{\mu\nu}^\la(\mf{g},G)$ is the multiplicity of a finite-dimensional irreducible representation for $(\mf{g},G)=(\mf{so}_{2n},O)$ and $(\mf{sp}_{2n},Sp)$.

\begin{ex}{\rm 
Let us give some examples of \eqref{eq:LR for osc} for oscillator representations.

(1) Suppose that $\mf{g}=\mf{sp}_{2n}$. Then it follows from \eqref{eq:LR equal} and the Littlewood-Richardson rule for finite-dimensional $\mf{so}_{2n}$-modules (see for example, \cite{Na}) that
\begin{equation}\label{eq:decom of two tensor of spin}
\begin{split}
&  V_{\mf{sp}_{2n}}^+\ot  V_{\mf{sp}_{2n}}^+ = \bigoplus_{k\ge 0}V_{\mf{sp}_{2n}}^{(2k,0)}, \\
&  V_{\mf{sp}_{2n}}^+\ot  V_{\mf{sp}_{2n}}^- = \bigoplus_{k\ge 0}V_{\mf{sp}_{2n}}^{(2k+1,0)}, \\
&  V_{\mf{sp}_{2n}}^-\ot  V_{\mf{sp}_{2n}}^- = V_{\mf{sp}_{2n}}^{(1,1)} \oplus \bigoplus_{k\ge 1}V_{\mf{sp}_{2n}}^{(2k,0)}.
\end{split}
\end{equation}

(2) Suppose that $\mf{g}=\mf{so}_{2n}$. Then we have for $l_1,l_2\in \Z_+$
\begin{equation}\label{eq:decom of two tensor of fundamentals}
\begin{split}
&  V'((l_1),{\mf{so}_{2n}})\ot V'((l_2),{\mf{so}_{2n}}) 
= \bigoplus_{r\ge 0}\bigoplus_{0\le k\le \min(l_1,l_2)}V'((l_1+l_2+r-s,r+s),{\mf{so}_{2n}}),
\end{split}
\end{equation}
by \eqref{eq:LR equal} and the Littlewood-Richardson rule for finite-dimensional $\mf{sp}_{2n}$-modules.}
\end{ex}

Let ${O}^{\mf c}_{\rm osc}$ be the category of $\mf{sp}_{2n}$-modules 
$V$ such that $V=\bigoplus_{\ell\ge 1}V_{\ell}$, where $V_\ell$ is a direct sum of $V(\la,\mf{sp}_{2n})$'s for $\la\in \mc{P}(O_\ell)_{\mf{sp}_{2n}}$ with finite multiplicity for each $\la$, and $V_\ell=0$ for all sufficiently large $\ell$. 
Similarly, we define ${O}^{\mf d}_{\rm osc}$ be the category of $\mf{so}_{2n}$-modules in the same manner, where $V(\la,\mf{sp}_{2n})$ for $\la\in \mc{P}(O_\ell)_{\mf{sp}_{2n}}$ is replaced by $V'(\la,\mf{so}_{2n})$'s for $\la\in \mc{P}(Sp_{2\ell})_{\mf{so}_{2n}}$.
Note that $\dim V_\mu<\infty$ for each weight of $\mu$.

If we further assume that ${\rm wt}(V)$ is finitely dominated for $V\in {O}^{\mf x}_{\rm osc}$ for $\mf{x}=\mf{c}, \mf{d}$, then ${O}^{\mf x}_{\rm osc}$ is a full subcategory of the parabolic BGG category of $\mf{g}$ with respect to the Levi subalgebra associated to $\{\upalpha_1,\dots\upalpha_{n-1}\}$, and it is closed under tensor product by \eqref{eq:LR for osc}. If we replace $\mf{sp}_{2n}$ with $\mf{so}_{2n}$ (resp. $\mf{so}_{2n}$ with $\mf{sp}_{2n}$) in the definition of ${O}^{\mf c}_{\rm osc}$ (resp. ${O}^{\mf d}_{\rm osc}$), then the corresponding category is that of finite-dimensional representations.

\begin{prop}
The categories ${O}^{\mf c}_{\rm osc}$ and ${O}^{\mf d}_{\rm osc}$ are semisimple tensor categories. 
\end{prop} 
\pf It follows immediately from \eqref{eq:LR for osc}.
\qed

\section{Quantum affine analogue of oscillator representations}\label{sec:category Osc}

In this section, we introduce a category of $U_q(X)$-modules ($X=C_{m}^{(1)}, D_{m+1}^{(1)}$), which can be viewed as a quantum affine analogue of ${O}^{\mf x}_{\rm osc}$ (${\mf x}=\mf{c}, \mf{d}$) introduced in Section \ref{sec:osc rep}, respectively. More generally, we introduce a category $\widehat{\mc O}^{\,\mf x}_{\rm osc,\e}$ of $\U_\e$-modules for any $\e$ in Section \ref{sec:truncation}.

\subsection{Subalgebras of finite type}\label{subsec:notation for osc}

We keep the notations in Section \ref{sec:truncation}.
Throughout this section, we assume that $\e={\bm\e},\un{\bm\e},\ov{\bm\e}$, ${\bm\e}',\un{\bm\e}',\ov{\bm\e}'$ unless otherwise specified. Let 
\begin{itemize}
 \item[$\bullet$] $\ring{\U}_\e$ : the subalgebra of $\U_\e$ generated by $e_i$, $f_i$, $k_\mu$ for $i\in I_\e\setminus \{0\}$ and $\mu\in P_\e$,
 
 \item[$\bullet$] $\dring{\U}_\e$ : the subalgebra of $\U_\e$ generated by $e_i$, $f_i$, $k_{\de_a}^{\pm 1}$ for $i\in I_\e\setminus \{0,n_\e\}$ and $a\in \I_\e$,
\end{itemize}
where $n_\e$ is the largest index in $I_\e$.
We regard $P_\e$ and $\bigoplus_{a\in \I_\e}\Z\de_a$ as the weight lattices for $\ring{\U}_\e$ and $\dring{\U}_\e$, respectively. 

The subalgebra $\dring{\U}_\e$ is a generalized quantum group of finite type $A$, which is isomorphic to the quantum superalgebra associated to a general linear Lie superalgebra under $\tau$ (and an extension by $\sigma_i$'s) in Proposition \ref{prop:tau iso}. Let us briefly recall its polynomial representations (cf.~\cite{BKK,KY}).

Let $\I_\e=\{i_1,i_2,\dots\}$ with $i_1>i_2>\dots$. For a partition $\nu=(\nu_i)_{i\ge 1}$, let $H^\e_{\nu}$ denote the tableau of shape $\nu$ with entries in $\I_\e$ defined inductively as follows:
\begin{itemize}
\item[(1)] If $\e_{i_1} =0$ (resp. $\e_{i_1}=1$), then fill the first row (resp. column) of $\nu$ with $i_1$.

\item[(2)] Suppose that a subdiagram $\eta$ of $\nu$ is filled with $i_1,\dots,i_k$ ($k\ge 1$). Then fill the first row (resp. column) of $\nu/\eta$ with $i_{k+1}$ if $\e_{i_{k+1}} =0$ (resp. $\e_{i_{k+1}}=1$).
\end{itemize} 
Let $\cP_\e$ be the set of $\nu$ such that $H^\e_{\nu}$ is well-defined (cf. \cite{CW}).
For $\nu\in\cP_\e$, let $V_\e(\nu)$ be the irreducible highest weight $\dring{\U}(\e)$-module with the highest weight $\sum_{i\in \I_\e} m_i\de_i$, where $m_i$ is the number of $i$'s in the tableau $H^\e_{\nu}$.
Recall that a tensor product $V_\e(\eta)\otimes V_\e(\zeta)$ is a finite direct sum of  $V_\e(\nu)$'s for $\eta,\zeta \in \cP_\e$.

For $m\ge 2$, let 
\begin{equation}\label{eq:e alternating}
{\bm\e}^{(m)}=(\underbrace{1,0,1,\dots,1,0,1}_{2m+1}).
\end{equation}
Choose an injective map $\iota_m: \I_{{\bm\e}^{(m)}}\longrightarrow \I_{{\bm\e}^{(m+1)}}$ such that $\iota_m(k)<\iota_m(l)$ for $k<l$, $\e_i=\e_{\iota_m(i)}$ for all $i$, and $\{1,2,2m+2,2m+3\}\subset \iota_m\left(\I_{{\bm\e}^{(m)}}\right)$.
Identifying ${\bm\e}^{(m)}$ as a subsequence of ${\bm\e}^{(m+1)}$ corresponding to the subsequence $\iota_m(\I_{{\bm\e}^{(m)}})$, we have a homomorphisms of $\phi_m: \dring{\U}_{{\bm\e}^{(m)}}\longrightarrow \dring{\U}_{{\bm\e}^{(m+1)}}$ by \cite[Theorem 4.3]{KY}, which naturally extends to a homomorphism $\phi_m:  {\U}_{{\bm\e}^{(m)}}\longrightarrow{\U}_{{\bm\e}^{(m+1)}}$.
Hence we may define a functor $\mf{tr\,}^{{\bm\e}^{(m+1)}}_{{\bm\e}^{(m)}} : \mc{C}({\bm\e}^{(m+1)})\longrightarrow \mc{C}({\bm\e}^{(m)})$ in the same way as in Proposition \ref{prop:truncation}, which is also defined on the $\ring{\U}_{{\bm\e}^{(m+1)}}$-modules whose weights satisfy the same conditions.

In general, for $\bm{\e}_1=\bm{\e}^{(m_1)}$ and $\bm{\e}_2=\bm{\e}^{(m_2)}$ with $m_2<m_1$, we have a functor $\mf{tr}^{\bm{\e}_1}_{\bm{\e}_2}$ by composing $\mf{tr\,}^{{\bm\e}^{(m+1)}}_{{\bm\e}^{(m)}}$'s. Similarly, we may define functors $\mf{tr\,}^{\un{\bm{\e}_1}}_{\un{\bm{\e}_2}}$ and $\mf{tr\,}^{\ov{\bm{\e}_1}}_{\ov{\bm{\e}_2}}$. 
The same works for the case of 
\begin{equation}\label{eq:e' alternating}
 {\bm\e'}^{(m)}=(\underbrace{0,1,0,\dots,0,1,0}_{2m+1}).
\end{equation}

\subsection{Categories ${\mc O}^{\mf x}_{\rm osc,\e}$ and $\widehat{\mc O}^{\,\mf x}_{\rm osc,\e}$}\label{subsec:fin osc category type x}
Let $\mf{x}=\mf{c}, \mf{d}$, and assume that 
\begin{equation}\label{eq:case of e}
 \e=
\begin{cases}
 {\bm\e},\,\un{\bm\e},\,\ov{\bm\e} & \text{when $\mf{x}=\mf{c}$},\\
 {\bm\e}',\,\un{\bm\e}',\,\ov{\bm\e}' & \text{when $\mf{x}=\mf{d}$}.
\end{cases}
\end{equation}
For $x\in \Bbbk^\times$, let 
\begin{equation}\label{eq:WW(x)}
 \WW_\e(x)=
 \begin{cases}
 \W_\e(x) & \text{when $\mf{x}=\mf{c}$}, \\
 \W^{\ot 2}_\e(x) & \text{when $\mf{x}=\mf{d}$}.
 \end{cases}
\end{equation}
We simply write $\WW_\e=\WW_\e(x)$ and $\W_\e^\pm=\W_\e^\pm(x)$ when we regard them as $\mathring{\U}_\e$-modules.

\begin{lem}\label{lem:semisimplicity of type A}
For $\ell\ge 1$, $\WW_\e^{\otimes \ell}$ is semisimple as a $\dring{\U}_{\e}$-module, and it is a direct sum of $V_\e(\nu)$'s for $\nu\in \cP_\e$.
\end{lem}
\pf For $l\ge 0$, let $(\W_\e)_l$ be the subspace of $\W_\e$ spanned by $\ket{\bf m}$ with $|{\bf m}|=l$. Then we have $(\W_\e)_l\cong V_\e((l))$. This proves the claim.
\qed\smallskip

\begin{lem}\label{lem:semisimplicity of finite type D}
For $\ell\ge 1$, $\WW_{\e}^{\otimes\ell}$ is semisimple as a $\mathring{\U}_{\e}$-module.
\end{lem}
\pf {\em Case 1}. Suppose that $\mf{x}=\mf{c}$.
We define an anti-involution $\eta$ on $\mathring{\U}_{\bm\e}$
by
\begin{align*}
\eta(k_\mu) & =k_\mu,\\
\eta(e_{i}) & =\begin{cases}
(-q^{2})^{\epsilon_{i}-\epsilon_{i+1}}q_{i}f_{i}k_{i}^{-1} & \text{if }1\leq i<n,\\
q_{n}f_{n}k_{n}^{-1} & \text{if }i=n,
\end{cases}\\
\eta(f_{i}) & =\begin{cases}
(-q^{2})^{\epsilon_{i+1}-\epsilon_{i}}q_{i}^{-1}k_{i}e_{i} & \text{if }1\leq i<n,\\
q_{n}^{-1}k_{n}e_{n} & \text{if }i=n,
\end{cases}
\end{align*}
for $\mu\in P_{\bm\e}$ and $i\in I_{\bm\e}\setminus\{0\}$. 
Then we have $(\eta\otimes\eta)\circ\Delta=\Delta\circ\eta$.
We define a nondegenerate symmetric bilinear form on $\WW_{\bm{\e}}=\W_{\bm{\e}}$ by 
\begin{equation}\label{eq:bilinear form on W}
 (\ket{\mathbf{m}},\ket{\mathbf{m}^{\prime}})=\delta_{\mathbf{m},\mathbf{m}^{\prime}}q^{-\sum_{i\in \I_{\bm{\e}}}\frac{m_{i}(m_{i}-1)}{2}}\prod_{i\in \I_{\bm{\e}}}[m_{i}]!,
\end{equation}
for ${\bf m}$, ${\bf m}'\in \Z^n_+(\e)$.
We may assume that $\Bbbk=\Q(q)$ if we ignore the action of $k_\La$. If we set $\mathcal{L}_{\infty}=\sum A_{\infty}\ket{\mathbf{m}}$,
where $A_{\infty}=\{f\in\Bbbk\,|\,f\text{ is regular at }q=\infty\}$,
then $(\mathcal{L}_{\infty},\mathcal{L}_{\infty})\subset A_{\infty}$.
Hence we obtain the $\mathbb{Q}$-valued induced form on $\mathcal{L}_{\infty}/q^{-1}\mathcal{L}_{\infty}$,
which is positive-definite since the image of $\{\,\ket{\mathbf{m}}\,|\,\mathbf{m}\in\mathbb{Z}_{+}^{n}(\bm{\e})\,\}$
in $\mathcal{L}_{\infty}/q^{-1}\mathcal{L}_{\infty}$ is an orthonormal
basis.

We have $(uv,w)=(v,\eta(u)w)$ for $u\in\ring{\U}_{\bm{\e}}$ and $v,\,w\in\WW_{\bm\e}$, which is straightforward to check. Therefore the semisimplicity of $\WW_{\e}^{\ot \ell}$ ($\e=\bm{\e}, \un{\bm{\e}}, \ov{\bm{\e}}$) follows from the arguments in the proof of \cite[Theorem 2.12]{BKK} (see also the proof of \cite[Lemma 6.5]{KL22}).\smallskip

{\em Case 2}. Suppose that $\mf{x}=\mf{d}$. In this case, we define an anti-involution $\eta'$ on $\mathring{\U}_{\bm\e'}$
by
\begin{align*}
\eta^\prime(k_\mu) & =k_\mu,\\
\eta^\prime(e_{i}) & =\begin{cases}
(-q^{2})^{\epsilon_{i}-\epsilon_{i+1}}q_{i}f_{i}k_{i}^{-1} & \text{if }1\leq i<n,\\
(-q^{-2})q_{n}f_{n}k_{n}^{-1} & \text{if }i=n,
\end{cases}\\
\eta^\prime(f_{i}) & =\begin{cases}
(-q^{2})^{\epsilon_{i+1}-\epsilon_{i}}q_{i}^{-1}k_{i}e_{i} & \text{if }1\leq i<n,\\
(-q^2)q_{n}^{-1}k_{n}e_{n} & \text{if }i=n,
\end{cases}
\end{align*}
for $\mu\in P_{\bm\e'}$ $i\in I_{\bm\e'}\setminus\{0\}$, and a nondegenerate symmetric bilinear form on $\WW_{\bm{\e}'}=\W_{\bm{\e}'}^{\ot 2}$ by 
\[
(\ket{\mathbf{m}}\otimes\ket{\mathbf{m}^\prime},\ket{\mathbf{m}^{\prime\prime}}\otimes\ket{\mathbf{m}^{\prime\prime\prime}})=(\ket{\mathbf{m}},\ket{\mathbf{m}^{\prime\prime}})(\ket{\mathbf{m}^\prime},\ket{\mathbf{m}^{\prime\prime\prime}}),
\]
where the bilinear form on the righthand side is given in \eqref{eq:bilinear form on W}. Then we apply the same argument as in {\em Case 1} to conclude that $\WW_\e^{\otimes \ell}$ is semisimple. \qed\smallskip

Let $\ell\ge 1$ be given, and let $G_\ell$ denote $O_\ell$ or $Sp_{2\ell}$.
For $\la\in \mc{P}(G_{\ell})$, let  
\begin{equation*}
\La_{\la,\e}= 
\begin{cases}
 r\ell \La + \sum_{i\in \I_\e} m_i\de_{i} & \text{if $\la\in \cP_\e$},\\
 0 & \text{otherwise},
\end{cases}
\end{equation*}
where $r=1$ for $\mf{x}=\mf{c}$, $r=2$ for $\mf{x}=\mf{d}$, and $m_i$ is the number of occurrences of $i$ in $H^{\epsilon}_\la$.
Put 
$$\mc{P}(G_\ell)_{\e}=\{\,\la\in \mc{P}(G_\ell)\,|\,\La_{\la,\e}\neq 0\,\}=\mc{P}(G_\ell)\cap\cP_\e.$$ 

\begin{df}{\rm 
For $\la\in \mc{P}(G_\ell)_{\e}$, let $\mc{V}_{\e}^\la$ be the irreducible highest weight $\ring{\U}_\e$-module with highest weight $\La_{\la,\e}$. }
\end{df}
For example, when $\mf{x}=\mf{c}$ we have
\begin{equation}\label{eq:W=V for level one}
\V^{(0)}_\e= \W_\e^+,\quad \V^{(1)}_\e=\W_\e^-.
\end{equation}

\begin{rem}\label{rem:q-osc finite type}
{\rm
(1) Suppose that $\mf{x}=\mf{c}$ and let $\e=\un{\bm\e}, \ov{\bm\e}$. 
We assume that $\mf{g}=\mf{sp}_{2m}$ for $\e=\un{\bm\e}$ and $\mf{g}=\mf{so}_{2(m+1)}$ for $\e=\ov{\bm\e}$.
If we identify $\texttt{P}$ in Section~\ref{sec:osc rep} with a subgroup of $P_\e$ by letting $\upvarpi_m=-\La$ (resp. $\upvarpi_{m+1}=\La$) and $\ude_i=\de_{2i}$ for $1\le i\le m$ (resp. $\ude_i=\de_{2i-1}$ for $1\le i\le m+1$) when $\e=\un{\bm\e}$ (resp. $\ov{\bm\e}$), then we have $\mathcal{P}(O_{\ell})_\epsilon=\mathcal{P}(O_{\ell})_{\mathfrak{g}}$ and $\La_{\la,\e}= \upvarpi_{\la,\mf{g}}$ \eqref{eq:Lambda highest weight} for $\la\in \mc{P}(O_\ell)_{\e}$.
Hence $\mc{V}_{\e}^\la$ is a $q$-analogue of $V(\la,{\mf{g}})$.

(2) Suppose that $\mf{x}=\mf{d}$ and let $\e=\un{\bm\e}', \ov{\bm\e}'$. 
We assume that $\mf{g}=\mf{so}_{2(m+1)}$ for $\e=\un{\bm\e}'$ and $\mf{g}=\mf{sp}_{2m}$ for $\e=\ov{\bm\e}'$.
If we identify $\texttt{P}$ with a subgroup of $P_\e$ by letting $\upvarpi_{m+1}=-\La$ (resp. $\upvarpi_m=2\La$) and $\ude_i=\de_{2i}$ for $1\le i\le m$ (resp. $\ude_i=\de_{2i-1}$ for $1\le i\le m+1$) when $\e=\un{\bm\e}'$ (resp. $\ov{\bm\e}'$), then we have $\mathcal{P}(Sp_{2\ell})_\epsilon=\mathcal{P}(Sp_{2\ell})_{\mathfrak{g}}$ and $\La_{\la,\e}= \upvarpi'_{\la,\mf{g}}$ \eqref{eq:Lambda highest weight-2} for $\la\in \mc{P}(Sp_{2\ell})_{\e}$.
Hence $\mc{V}_{\e}^\la$ is a $q$-analogue of $V'(\la,{\mf{g}})$.
}
\end{rem}

\begin{lem}\label{lem:semisimplicity of osc}
A highest weight $\ring{\U}_\e$-submodule in $\WW_\e^{\ot \ell}$ is isomorphic to $\mc{V}^{\la}_\e$ for some $\la\in \mc{P}(G_\ell)_\e$.
\end{lem}
\pf We prove the case only when $\mf{x}=\mf{c}$ since the proof for $\mf{x}=\mf{d}$ is almost the same.

Let $V$ be a highest weight $\ring{\U}_\e$-submodule of $\W_\e^{\ot \ell}$ with highest weight vector $v$, which is irreducible by Lemma \ref{lem:semisimplicity of finite type D}.

{\em Case 1}. Suppose that $\e=\un{\bm\e}, \ov{\bm\e}$. By Remark \ref{rem:wt of m} and \eqref{eq:hw for spin}, the classical limit $\overline{\mathcal{W}^{\tau}_\epsilon}$ of $\W_\e=\W_\e^+\oplus \W_\e^-$ is isomorphic to $W_{\mf g}=V^+_{\mf g}\oplus V^-_{\mf g}$. 
By Lemma~\ref{lem:classical limit of hw module} and Theorem~\ref{thm:Howe duality}, the classical limit $\overline{V^{\tau}}$ of $V$ is equal to $V^\la_{\mf g}$ for some $\lambda \in\mc{P}(O_\ell)_{\e}$ (see Section~\ref{subsec:cl-lim-U_D(e)-modules}), and hence the weight of $v$ is $\La_{\la,\e}$ by Remark~\ref{rem:q-osc finite type}.

{\em Case 2}. Suppose that $\e={\bm\e}$.
By Lemma \ref{lem:semisimplicity of type A}, we have as a $\dring{\U}_{\bm\e}$-module $V=\bigoplus_{\nu\in S}V_{\bm\e}(\nu)^{\oplus m_\nu}$ for some $S\subset \cP_{\bm\e}$ and $m_\nu>0$ for $\nu\in S$. If $v\in V_{\bm\e}(\nu)$ for some $\nu\in S$, then $|\nu|=\sum_{i\ge 1}\nu_i$ is minimal among $\eta\in S$.

Note that $\mf{tr}_{\ov{\bm\e}}(V)=\bigoplus_{\nu\in \ov{S}}V_{\ov{\bm\e}}(\nu)^{\oplus m_\nu}$ for some $\ov{S}\subset S$ since $\mf{tr}_{\ov{\bm\e}}(V_{\bm\e}(\nu))=V_{\ov{\bm\e}}(\nu)$ if $\nu\in \cP_{\ov{\bm\e}}$ and $0$ otherwise (\cite[Proposition 4.5]{KY}). 
Since $\mf{tr}^{{\bm\e}^{(m+1)}}_{{\bm\e}^{(m)}}(\W_{{\bm\e^{(m+1)}}})\cong \W_{{\bm\e^{(m)}}}=\W_{{\bm\e}}$ (cf.~Section \ref{subsec:notation for osc}), we may assume that $m$ is large enough and hence $\nu\in \ov{S}$. Let $w$ be a $\dring{\U}_{\ov{\bm\e}}$-highest weight vector in $V_{\ov{\bm\e}}(\nu)\subset \mf{tr}_{\ov{\bm\e}}(V)$. By minimality of $\nu$, $w$ is also a $\ring{\U}_{\ov{\bm\e}}$-highest weight vector. By {\em Case 1}, the $\ring{\U}_{\ov{\bm\e}}$-weight of $w$ is $\La_{\nu,\ov{\bm \e}}$, which implies that the weight of $v$ is $\La_{\nu,\bm{\e}}$. 
\qed

\begin{thm}\label{thm:trunc sends simple to simple}
Let $(\epsilon,\underline{\epsilon},\overline{\epsilon})=(\boldsymbol{\epsilon},\boldsymbol{\underline{\epsilon}},\boldsymbol{\overline{\epsilon}})\text{ or }(\boldsymbol{\epsilon}',\boldsymbol{\underline{\epsilon}}',\boldsymbol{\overline{\epsilon}}')$
and $\la\in \mc{P}(G_\ell)_{\e}$ be given. Then 
\begin{itemize}
 \item[(1)] the set of weights in $\V^\la_{\e}$ satisfies
 \begin{equation*}
{\rm wt}(\V^\la_{\e})\subset r\ell\La + \sum_{i\in \I_{\e}}\Z_+\de_i, 
\end{equation*}
 \item[(2)] the image of $\mc{V}^{\la}_{\e}$ under $\mf{tr}_{\e'}$ is given by
 \begin{equation*}
\mf{tr}_{\e'}\left(\mc{V}^{\la}_{\e}\right)\cong
\begin{cases}
\mc{V}^{\la}_{\e'}  & \text{if $\la\in \mc{P}(G_\ell)_{\e'}$},\\
0 & \text{otherwise},
\end{cases}\quad (\e'=\un{\e}, \ov{\e}).
\end{equation*}
\end{itemize}
\end{thm}
\pf The proof is similar to that of \cite[Theorem 6.17]{KL22}.
We give a self-contained proof for the reader's convenience. Also, we prove the case only when $\mf{x}=\mf{c}$ since the proof for $\mf{x}=\mf{d}$ is similar.

We have seen in the proof of Lemma \ref{lem:semisimplicity of osc} that
for $\e=\un{\bm\e}, \ov{\bm\e}$
\begin{equation}\label{eq:finite type 1}
\W_{\e}^{\ot \ell} \cong 
\bigoplus_{\la\in \mc{P}(O_\ell)_{\e}}
{\mc{V}_{{\e}}^\la}^{\oplus d^\la},
\end{equation}
where $d^\la$ is the dimension of $V_{O_\ell}(\la)$. 
Again by Lemma \ref{lem:semisimplicity of osc}, we have
\begin{equation*}\label{eq:finite type 2}
\W_{\bm\e}^{\ot \ell} \cong \bigoplus_{\la\in S}\ {\mc{V}_{\bm\e}^\la}^{\oplus d^{\la}_{\bm\e}},
\end{equation*}
for some $S\subset \mc{P}(O_\ell)_{\bm\e}$ and $d_{\bm\e}^\la\in \mathbb{N}$.  

{\em Step 1}.
Consider the case when $\e=\ov{\bm\e}$. Let $\la\in S$ be given
and let $V=\mc{V}_{\bm\e}^\la$ with a highest weight vector $v$.

Suppose that $\la\in \mc{P}(O_\ell)_{\ov{\bm\e}}$.
Then we have $V_{\ov{\bm\e}}(\la)\subset \mf{tr}_{\ov{\bm\e}}(V)$ since $v\in V_{\bm\e}(\la)$ as a $\dring{\U}_{\bm\e}$-module and $\mf{tr}_{\ov{\bm\e}}(V_{\bm\e}(\la))=V_{\ov{\bm\e}}(\la)\neq 0$. 
This implies $\mc{V}_{\ov{\bm\e}}^{\la}\subset \mf{tr}_{\ov{\bm\e}}(V)$ by the minimality of $\la$ in $V$ (cf.~Lemma \ref{lem:semisimplicity of osc}).
To prove the equality, we consider the classical limit  $\ov{V^\tau}$, which is a highest weight $U(\mf{osp}(\bm{\e}))$-module with highest weight $\La_{\la,\bm{\e}}$ by Lemma \ref{lem:classical limit of hw module}. 
We may define $\mf{tr}_{\ov{\bm\e}}\left(\ov{V^\tau}\right)$ in the same way as in \eqref{eq:truncation-obj}, which is a $U(\mf{so}_{2(m+1)})$-module.
Now by the same method in the proof of \cite[Lemma 3.5]{CL} we conclude that $\mf{tr}_{\ov{\bm\e}}\left(\ov{V^\tau}\right)$ is also a highest weight $U(\mf{so}_{2(m+1)})$-module with highest weight $\La_{\la,\ov{\bm\e}}$, and hence irreducible since it is integrable. 
This implies that $\mf{tr}_{\ov{\bm\e}}(V)=\mc{V}_{\ov{\bm\e}}^{\la}$.

Suppose that $\la\not\in \mc{P}(O_\ell)_{\ov{\bm\e}}$.
Choose ${\bm\e}_0={\bm\e}^{(m_0)}$ for some $m_0>m$ such that $\la\in \mc{P}(O_\ell)_{\ov{{\bm\e}_0}}$. 
We have 
$\V^\la_{{\bm \e}_0}\subset \W_{{\bm \e}_0}(x)^{\ot \ell}$ and 
$\mf{tr}^{{\bm\e}_0}_{{\bm \e}}\left(\V^\la_{{\bm \e}_0}\right)=V$ (cf.~\cite[Lemma 6.16]{KL22}).
Moreover, by the previous arguments, we have 
$\mf{tr}^{\bm{\e}_0}_{\ov{\bm{\e}_0}}\left(\V^\la_{\bm{\e}_0}\right)\cong \V^\la_{\ov{\bm{\e}_0}}$.
Since
$\mf{tr}^{\ov{\bm{\e}_0}}_{\ov{\bm{\e}}}(\V^\la_{\ov{\bm{\e}_0}})=0$, we have 
\begin{equation*}
0= \mf{tr}^{\ov{\bm{\e}_0}}_{\ov{\bm\e}}(\V^\la_{\ov{\bm{\e}_0}})=\mf{tr}^{\ov{\bm{\e}_0}}_{\ov{\bm\e}}\circ \mf{tr}^{\bm{\e}_0}_{\ov{\bm{\e}_0}}\left(\V^\la_{\bm{\e}_0}\right)
 =\mf{tr}^{\bm\e}_{\ov{\bm\e}}\circ \mf{tr}^{\bm{\e}_0}_{\bm\e}\left(\V^\la_{\bm{\e}_0}\right)
 =\mf{tr}^{\bm\e}_{\ov{\bm\e}}\left(V\right).
\end{equation*}

{\em Step 2}.
Suppose that $\la\in \mc{P}(O_\ell)_{\bm\e}$ is given. 
Consider again $\bm{\e}_0=\bm{\e}^{(m_0)}$ as chosen above.
We have $\V^\la_{\bm{\e}_0}\subset \W_{\bm{\e}_0}^{\ot \ell}$ by {\em Step 1},
and 
$\V^\la_{\bm\e} = \mf{tr}^{\bm{\e}_0}_{\bm\e}\left(\V^\la_{\bm{\e}_0}\right)\subset \W_{\bm\e}(x)^{\ot \ell}$, that is, $\la\in S$. Therefore we have $S=\mc{P}(O_\ell)_{\bm\e}$ and $d^\la_{\bm\e}=d^\la$ by {\em Step 1}. In particular, we have ${\rm wt}(\V^\la_{\bm\e})\subset \ell\La + \sum_{i\in \I_{\bm\e}}\Z_+\de_i$.
\smallskip 

{\em Step 3}. Consider the case when $\e=\un{\bm\e}$.
Let $\la\in \mc{P}(O_\ell)_{\un{\bm\e}}\subset \mc{P}(O_\ell)_{{\bm\e}}$ be given.
Let $v$ be a highest weight vector in ${\mc{V}_{\bm\e}^\la}$.
We have $\V^\la_{\bm\e}\subset \W_{\bm\e}(x)^{\ot \ell}$ by {\em Step 2}, and may assume that $v\in V_{\bm\e}(\la)$ by Lemma \ref{lem:semisimplicity of type A}.

If $w$ is a $\dring{\U}_{\un{\bm\e}}$-highest weight vector in $V_{\un{\bm\e}}(\la)=\mf{tr}^{\bm\e}_{\un{\bm\e}}\left(V_{{\bm\e}}\right)$, then we have by minimality of $\la$ in ${\mc{V}_{\bm\e}^\la}$ (cf.~Lemma \ref{lem:semisimplicity of osc}) that $w$ is a $\ring{\U}_{\un{\bm\e}}$-highest weight vector.
By Lemma \ref{lem:semisimplicity of osc}, the highest weight $\ring{\U}_{\un{\bm\e}}$-module generated by $w$ is isomorphic to $\mc{V}^\la_{\un{\bm\e}}$ and
\begin{equation}\label{eq:tr to 0}
 \mc{V}^\la_{\un{\bm\e}} \subset 
 \mf{tr}^{\bm\e}_{\un{\bm\e}}\left( \mc{V}^\la_{\bm\e}\right).
\end{equation}
It follows from \eqref{eq:finite type 1} and {\em Step 2} that the equality holds in \eqref{eq:tr to 0} and 
$\mf{tr}^{\bm\e}_{\un{\bm\e}}\left( \mc{V}^\la_{\bm\e}\right)=0$ for $\la\in \mc{P}(O_\ell)_{\bm\e}\setminus \mc{P}(O_\ell)_{\un{\bm\e}}$.
\qed

\begin{cor}\label{cor:decomp of W for e}
For each $\e$,
we have the following decomposition:
\begin{equation*} 
\WW_{\e}^{\ot \ell} \cong \bigoplus_{\la\in \mc{P}(G_\ell)_{\e}}\ {\mc{V}_{\e}^\la}^{\oplus d^{\la}},
\end{equation*}
where $d^\la=\dim V_{G_\ell}(\la)$.
\end{cor}

\begin{df}
{\rm
Let $\mc{O}^{\mf x}_{\text{osc},\e}$ be the category of $\ring{\U}_{\e}$-modules 
$V$ such that $V=\bigoplus_{\ell\ge 1}V_{\ell}$, where $V_\ell$ is a direct sum of $\mc{V}_{\e}^\la$'s for $\la\in \mc{P}(G_\ell)_{\e}$ with finite multiplicity for each $\la$, and $V_\ell=0$ for all sufficiently large $\ell$. 

} 
\end{df}

\begin{prop}\label{prop:semisimplicity of O_osc}
The category $\mc{O}^{\mf x}_{\rm osc,\e}$ is a semisimple tensor category. 
\end{prop} 
\pf Let $\mu\in \mc{P}(G_\ell)_\e$ and $\nu\in \mc{P}(G_{\ell'})_\e$ be given. 
It follows from Corollary \ref{cor:decomp of W for e} that 
$\V_{\e}^\mu\ot \V_{\e}^\nu \subset \WW_{\e}^{\ot{(\ell+\ell')}}$ is semisimple, where 
\begin{equation*}
\mc{V}_{\e}^\mu\ot \V_{\e}^\nu = \bigoplus_{\la\in \mc{P}(O_{\ell+\ell'})_{\e}}
(\V_{\e}^\la)^{\oplus c_{\mu\nu}^\la(\e)},
\end{equation*}  
for some $c_{\mu\nu}^\la(\e)\in\Z_+$. 
%
%
This implies that $V_\ell\ot W_{\ell'} \in \mc{O}^{\mf c}_{\rm osc,\e}$ for $V,W\in \mc{O}^{\mf c}_{\rm osc,\e}$, and hence $V\ot W\in \mc{O}^{\mf c}_{\rm osc,\e}$.
\qed\smallskip

Now we define the following.

\begin{df}\label{def:aff osc c and d}
{\rm
Let $\widehat{\mc O}^{\,\mf x}_{\rm osc,\e}$ be the category of $\U_\e$-modules $V$ such that $V$ belongs to ${\mc O}^{\,\mf x}_{\rm osc,\e}$ when restricted as a $\ring{\U}_\e$-module.
}
\end{df}
The category $\widehat{\mc O}^{\,\mf x}_{\rm osc,\e}$ is a full subcategory of $\mc{C}(\e)$ by Theorem \ref{thm:trunc sends simple to simple}(1), which is closed under taking submodules, and quotients. It is also closed under tensor product by Proposition \ref{prop:semisimplicity of O_osc}. For $x\in \Bbbk^\times$, we have $\WW_\e(x)^{\ot\ell}\in \widehat{\mc O}^{\,\mf x}_{\rm osc,\e}$ ($\ell\ge 1$) by Corollary \ref{cor:decomp of W for e}. 
Note that a $\U_\e$-module in $\widehat{\mc O}^{\,\mf x}_{\rm osc,\e}$ is  finite-dimensional when $\e=\ov{\bm\e}, \ov{\bm\e}'$. If we further assume that ${\rm wt}(V)$ is finitely dominated, then $\widehat{\mc O}^{\,\mf x}_{\rm osc,{\e}}$ ($\e=\un{\bm\e}, \un{\bm\e}'$) is a subcategory of $\widehat{\mc O}$ introduced in \cite{MY}.

\begin{prop}\label{prop:truncation on affine osc}
Suppose that $(\e,\un{\e},\ov{\e})=({\bm\e},\un{\bm\e},\ov{\bm\e})$ and $({\bm\e}',\un{\bm\e}',\ov{\bm\e}')$ for $\mf{x}=\mf{c}, \mf{d}$, respectively. Then the restrictions of $\mf{tr}^\e_{\un\e}$ and $\mf{tr}^\e_{\ov\e}$ on $\widehat{\mc O}^{\,\mf x}_{\rm osc,{\e}}$ give well-defined exact monoidal functors
\begin{equation*}
\xymatrixcolsep{3pc}\xymatrixrowsep{0.3pc}\xymatrix{
 & \widehat{\mc O}^{\,\mf x}_{\rm osc,{\e}}  \ar@{->}_{\mf{tr}^\e_{\un{\e}}}[dl]\ar@{->}^{\mf{tr}^\e_{\ov{\e}}}[dr] &  \\
 \widehat{\mc O}^{\,\mf x}_{\rm osc,\un{\e}} & &  \widehat{\mc O}^{\,\mf x}_{\rm osc,\ov{\e}}
}.
\end{equation*}
\end{prop}
\pf It follows from Theorem \ref{thm:trunc sends simple to simple}(2).
\qed

\subsection{Universal $R$ matrix}
Let $z$ be an indeterminate. For a $\U_\e$-module $V$ in $\mc{C}(\e)$, we let $V_{\rm aff} = \Bbbk[z,z^{-1}]\ot V$ be the affinization of $V$ defined as a $\U_\e$-module in a usual way, that is,
$k_\mu$, $e_i$, and $f_i$ act as $1\ot k_\mu$, $z^{\delta_{i,0}}\ot e_i$, and $z^{-\delta_{i,0}}\ot f_i$, respectively ($\mu\in P_\e$, $i\in I_\e$).
For $x\in \Bbbk^\times$, we let $V_x = V_{\rm aff}/(z-x)V_{\rm aff}$,
where we understand the map sending $g(z)\ot v$ to $zg(z)\ot v$ for $g(z)\ot v\in V_{\rm aff}$ as an isomorphism of $\U_\e$-module. 
For example, if $V=\W_\e(1)$, then $V_x=\W_\e(x)$ for $x\in \Bbbk^\times$.

Let $V, W$ be $\U_\e$-modules in $\mc{C}(\e)$, and let
$V_{\rm aff} = \Bbbk[z_1,z_1^{-1}]\ot V$ and 
$W_{\rm aff} = \Bbbk[z_2,z_2^{-1}]\ot W$, where $z_1,z_2$ are indeterminates. 
Take the completions $V_{\rm aff}\, \widehat{\ot}\, W_{\rm aff}$ and $W_{\rm aff}\, \td{\ot}\, V_{\rm aff}$ defined as in \cite[Section 7]{Kas02} with respect to a suitably extended weight lattice (for example $P_{\rm af}^e$ when $\e=\bm{\e}$). 
Let $\Pi_{\bq} : V_{\rm aff}\, \widehat{\ot}\, W_{\rm aff} \longrightarrow V_{\rm aff}\, \widehat{\ot}\, W_{\rm aff}$ be given by 
\begin{equation*}
\Pi_{\bq}(v\ot w)=\bq(\mu,\nu)\, v\ot w,
\end{equation*}
for $v=f(z_1)\ot v_0\in V_{\rm aff}$ and $w=g(z_2)\ot w_0\in W_{\rm aff}$ with $v_0\in V_\mu$ and $w_0\in W_\nu$. 
Let $s :V_{\rm aff}\, \widehat{\ot}\, W_{\rm aff} \longrightarrow W_{\rm aff}\, \td{\ot}\, V_{\rm aff}$ be the map given by $s(v\ot w)=w\ot v$.
Then as in \cite[Theorem 32.1.5]{Lu93}, we have an isomorphism of $\U_\e$-modules
\begin{equation}\label{eq:univ R matrix}
\xymatrixcolsep{2pc}\xymatrixrowsep{3pc}\xymatrix{
\mc{R}^{\rm univ}:= \Theta \circ \Pi_{\bq} \circ s : 
V_{\rm aff}\, \widehat{\ot}\, W_{\rm aff}\ \ar@{->}[r] & \ W_{\rm aff}\, \td{\ot}\, V_{\rm aff}
}
\end{equation}
(cf.~\cite[Theorem 3.7]{KL20}).

\section{Irreducible representations in $\widehat{\mc O}^{\,\mf c}_{\rm osc,\e}$}\label{sec:irr of type c}

In this section, we construct a family of irreducible representations in $\widehat{\mc O}^{\,\mf c}_{\rm osc,\e}$ and explain a connection between $q$-oscillator representations of type $C_m^{(1)}$ and finite-dimensional irreducible representations of type $D_{m+1}^{(1)}$.
%
We assume the notations in Sections \ref{subsec:truncation c} and \ref{subsec:fin osc category type x}.

\subsection{Normalized $R$ matrix on $\W^{\sigma_1}_{\e}(x)\ot \W^{\sigma_2}_{\e}(y)$}\label{subsec:norm R matrix for osc C} 

For $\sigma\in \{+,-\}$, let $\W_\e^\sigma=\W_\e^\sigma(1)$ and let $v^\sigma_\e$ be a $\ring{\U}_\e$-highest weight vector in $\W_\e^\sigma$ given by
\begin{equation*}
 v^\sigma_\e=
\begin{cases}
 \ket{{\bm 0}} & \text{if $\sigma=+$},\\
 \ket{{\be}_{i_\e}} & \text{if $\sigma=-$},
\end{cases}
\end{equation*}
where $i_\e$ is the largest entry in $\I_\e$.

\begin{prop}\label{prop:irreducibility of W tensor W}
For $\sigma_1,\sigma_2\in \{+,-\}$ and generic $x_1,x_2\in \Bbbk^\times$, $\W_\e^{\sigma_1}(x_1)\ot \W_\e^{\sigma_2}(x_2)$ is an irreducible $\U_\e$-module.
\end{prop}
\pf The cases when $\e=\ov{\bm\e}$ and $\un{\bm\e}$ are known (see \cite{AK} and \cite{KO}, respectively).
It suffices to show the case when $\e=\bm{\e}$. 
Suppose that $\sigma_1=\sigma_2=+$. 
Note that 
$\W^+_{\un{\bm\e}}(x_1)\ot \W^+_{\un{\bm\e}}(x_2) \subset \W^+_{{\bm\e}}(x_1)\ot \W^+_{{\bm\e}}(x_2)$, and
\begin{equation*}
 \W_{{\e}}^+ (x_1)\ot \W_{{\e}}^+ (x_2) \cong \bigoplus_{k\ge 0}\V_{{\e}}^{(2k,0)} \quad (\e=\bm{\e}, \un{\bm \e}),
\end{equation*}
as a $\mathring{\U}_\e$-module by \eqref{eq:decom of two tensor of spin} and Theorem \ref{thm:trunc sends simple to simple} (cf. Remark~\ref{rem:q-osc finite type}) since $\mc{P}(O_2)_{\un{\bm\e}}=\mc{P}(O_2)_{\bm\e}$.  

Let $V$ be a nonzero $\U_{\bm\e}$-submodule of $\W^+_{{\bm\e}}(x_1)\ot \W^+_{{\bm\e}}(x_2)$.
Given any nonzero $v\in V$, we may obtain a nonzero vector in $\W^+_{\un{\bm\e}}(x_1)\ot \W^+_{\un{\bm\e}}(x_2)$ by applying the action of $e_{2i-1},f_{2i-2}\in\mathcal{U}_{\boldsymbol{\epsilon}}$ ($i=1,2,\dots,m+1$), which implies that $V$ includes $\W^+_{\un{\bm\e}}(x_1)\ot \W^+_{\un{\bm\e}}(x_2)$ since it is irreducible. 
For $k\ge 0$, we have $\V_{\un{\bm{\e}}}^{(2k,0)}\subset \V_{\bm{\e}}^{(2k,0)}$ and hence $\V_{\bm{\e}}^{(2k,0)}\subset V$. Hence $V=\W^+_{{\bm\e}}(x_1)\ot \W^+_{{\bm\e}}(x_2)$. The proof for the case when $\sigma_1\neq \sigma_2$ is similar.
\qed\smallskip

For $\bm\sigma=(\sigma_1,\sigma_2)$ with $\sigma_1, \sigma_2\in \{+,-\}$, let 
\begin{equation*}
 S^{\bm\sigma}=
\begin{cases}
 \{\,(2k,0)\,|\,k\in\Z_+\,\} & \text{if $\bm\sigma=(+,+)$},\\
 \{\,(2k+1,0)\,|\,k\in\Z_+\,\} & \text{if $\bm\sigma=(+,-)$ or $(-,+)$},\\
 \{(1,1)\}\cup \{\,(2k,0)\,|\,k\ge 1\,\} & \text{if $\bm\sigma=(-,-)$}.
\end{cases}
\end{equation*}
We have as a $\ring{\U}_\e$-module
\begin{equation}\label{eq:decomp of two tensor of W}
 \W_{{\e}}^{\sigma_1}\ot \W_{{\e}}^{\sigma_2} \cong \bigoplus_{\la\in S^{\bm\sigma}_\e}\V_\e^\la,
\end{equation}
where $S^{\bm\sigma}_\e=S^{\bm\sigma}\cap \mc{P}(O_2)_\e$, by \eqref{eq:decom of two tensor of spin} and Theorem \ref{thm:trunc sends simple to simple} (cf. Remark~\ref{rem:q-osc finite type}).
Put $\la^0=(0,0)$, $(1,0)$, and $(1,1)\in S^{\bm\sigma}_\e$ when $\bm\sigma=(+,+)$, $(\pm,\mp)$, and $(-,-)$, respectively.

Put $z=z_1 /z_2$. Let $\W_{\e}^{\sigma_{i}}(z_{i}):=(\W_{\e}^{\sigma_{i}})_{{\rm aff}}=\Bbbk[z_{i}^{\pm1}]\ot \W_{\e}^{\sigma_{i}}$ for $i=1,2$.
The universal $R$ matrix \eqref{eq:univ R matrix} induces a map
\begin{equation*}
\xymatrixcolsep{4pc}\xymatrixrowsep{3pc}\xymatrix{
\mc{R}^{\rm univ}_{\bm\sigma,\e}(z):
\W_{\e}^{\sigma_1}(z_1)\, \widehat{\ot}\,\W_{\e}^{\sigma_2}(z_2) \ \ar@{->}[r] &\ \W_{\e}^{\sigma_2}(z_2) \, \widetilde{\ot}\, \W_{\e}^{\sigma_1}(z_1) },
\end{equation*}
(see \cite[(7.6)]{Kas02}), where
$\mc{R}^{\rm univ}_{{\bm\sigma},\e}(z)\big\vert_{\V_\e^{\la^0}}= a(z){\rm id}_{\V_\e^{\la^0}}$
for some nonzero $a(z)\in \Bbbk\llbracket z\rrbracket$. 
Define the normalized $R$ matrix by
\begin{equation}\label{eq:norm R}
\mc{R}_{{\bm\sigma},\e}(z)= a(z)^{-1}\mc{R}^{\rm univ}_{{\bm\sigma},\e}(z).
\end{equation}
Then we have a unique $\Bbbk[z_1^{\pm 1},z_2^{\pm 1}]\ot \U_\e$-linear map
\begin{equation}\label{eq:R matrix on two tensor of W}
\xymatrixcolsep{3pc}\xymatrixrowsep{3pc}\xymatrix{
\mc{R}_{{\bm\sigma},\e}(z) : \W_{\e}^{\sigma_1}(z_1) \, {\ot}\,\W_{\e}^{\sigma_2}(z_2) \ \ar@{->}[r] &  \
\Bbbk(z_1,z_2)\ot_{\Bbbk[z_1^{\pm 1},z_2^{\pm 1}]}\left(\W_{\e}^{\sigma_2}(z_2)\, {\ot}\,\W_{\e}^{\sigma_1}(z_1) \right)
}
\end{equation}
such that 
\begin{equation}\label{eq:normalization}
 \mc{R}_{{\bm\sigma},\e}(z)\big\vert_{\V_\e^{\la^0}}= {\rm id}_{\V_\e^{\la^0}},
\end{equation}
since $\Bbbk(z_1,z_2)\ot_{\Bbbk[z_1^{\pm 1},z_2^{\pm 1}]}\left(\W_{\e}^{\sigma_1}(z_1)\, {\ot}\,\W_{\e}^{\sigma_2}(z_2)\right)$ is an irreducible $\Bbbk(z_1,z_2)\ot \U_\e$-module by Proposition \ref{prop:irreducibility of W tensor W}.

\begin{lem}\label{lem:tr of norm R}
Under the above hypothesis, we have
\begin{equation*}
  \mf{tr}^{\bm\e}_{\e}\left(\mc{R}_{{\bm\sigma},\bm{\e}}(z)\right)
 =\mc{R}_{{\bm\sigma},\e}(z)\quad (\e=\un{\bm\e}, \ov{\bm\e}).
\end{equation*} 
\end{lem}
\pf It follows from \eqref{eq:W for e}, \eqref{eq:decomp of two tensor of W}, Theorem \ref{thm:trunc sends simple to simple}, and the uniqueness of \eqref{eq:R matrix on two tensor of W}.
\qed\smallskip

For $\la\in S^{{\bm\sigma}}_\e$, we define a $\mathring{\U}_\e$-linear map
$\mc P^{{\bm\sigma}}_{\la,\e} : 
\W_{\e}^{\sigma_1} \otimes \W_{\e}^{\sigma_2} \longrightarrow \W_{\e}^{\sigma_2} \otimes \W_{\e}^{\sigma_1}$ 
as follows:

\begin{itemize}
 \item[(1)] Let $v_\la$ and $v'_\la$ be the $\mathring{\U}_{\un{\bm\e}}$-highest weight vectors of $\mc{V}^{\la}_{\un{\bm\e}}$ in $\W_{\un{\bm\e}}^{\sigma_1} \otimes \W_{\un{\bm\e}}^{\sigma_2}$ and $\W_{\un{\bm\e}}^{\sigma_2}\otimes \W_{\un{\bm\e}}^{\sigma_1}$ respectively \cite[Proposition 5]{KO}.
 
 \item[(2)] Let 
 $\mc P^{{\bm\sigma}}_{\la,\un{\bm\e}}$ be the map determined by
\begin{equation*}
\begin{split}
&
\mc P^{{\bm\sigma}}_{\la,\un{\bm\e}}(v_\mu)=\de_{\la\mu}v'_\la.
\end{split}
\end{equation*}
 
 \item[(3)] Let $\mc P^{{\bm\sigma}}_{\la,{\bm\e}}$ and $\mc P^{{\bm\sigma}}_{\la,\ov{\bm\e}}$ to be the unique maps satisfying 
\begin{equation*}
 \mf{tr}^{\bm\e}_{\un{\bm\e}}\left(\mc P^{{\bm\sigma}}_{\la,{\bm\e}}\right)
 =\mc P^{{\bm\sigma}}_{\la,\un{\bm\e}},\quad 
 \mf{tr}^{\bm\e}_{\ov{\bm\e}}\left(\mc P^{{\bm\sigma}}_{\la,{\bm\e}}\right)
 =\mc P^{{\bm\sigma}}_{\la,\ov{\bm\e}}.
\end{equation*}
\end{itemize}

\begin{thm}\label{thm:spectral decomp of W}
We have the following:
\begin{equation*}
\mc{R}_{\bm{\sigma},{\bm\e}}(z)=
\begin{cases}
\sum_{k\in \Z_+} \prod_{j=1}^{k}\frac{1-q^{4j-2}z}{z-q^{4j-2}}\mc{P}^{\bm{\sigma}}_{(2k,0),{\bm\e}} & \text{if $\bm{\sigma}=(+,+)$}, 

\\[4mm] \mc{P}^{\bm{\sigma}}_{(1,1),\boldsymbol{\epsilon}}+
 \sum_{k\ge 1} \prod_{j=1}^{k}\frac{1-q^{4j-2}z}{z-q^{4j-2}}\mc{P}^{\bm{\sigma}}_{(2k,0),{\bm\e}} & \text{if $\bm{\sigma}=(-,-)$},
 
\\[4mm] \sum_{k\in \Z_+} \prod_{j=1}^{k}\frac{1-q^{4j}z}{z-q^{4j}}\mc{P}^{\bm{\sigma}}_{(2k+1,0),{\bm\e}} & \text{if $\bm{\sigma}=(+,-),(-,+)$},
\end{cases}
\end{equation*}
where the coefficient of $\mc{P}^{{\bm\sigma}}_{\la,\bm{\e}}$ is assumed to be 1 for $k=0$.
\end{thm}
\pf 
By \eqref{eq:decomp of two tensor of W}, $\mc{R}_{{\bm\sigma},{\bm\e}}(z)$ is a (not necessarily finite) linear combination of $\mc P^{{\bm\sigma}}_{\la,{\bm\e}}$ for $\la\in S^{\bm\sigma}_{\bm\e}$. In case of $\un{\bm\e}$, it is known that $\mc{R}_{\bm{\sigma},{\un{\bm\e}}}(z)$ is given by the formula above \cite[Proposition 11]{KO}, where $S^{\bm\sigma}_{\bm\e}=S^{\bm\sigma}_{\un{\bm\e}}$. Note that the scalar $\frac{z{\bf v}^{-1}}{z-1}$ (${\bf v}={\bf i}q^{\hf}$) for $\bm{\sigma}=(\pm,\mp)$ in \cite{KO} is removed here since 
$\mc{R}_{\bm{\sigma},{\bm\e}}(z)$ is normalized as in \eqref{eq:normalization}.

 By Theorem \ref{thm:trunc sends simple to simple} and Lemma \ref{lem:tr of norm R}, the coefficient of $\mc P^{{\bm\sigma}}_{\la,{\bm\e}}$ and $\mc P^{{\bm\sigma}}_{\la,\un{\bm\e}}$ in $\mc{R}_{{\bm\sigma},{\bm\e}}(z)$ and $\mc{R}_{{\bm\sigma},\un{\bm\e}}(z)$ respectively, are equal. This completes the proof.
\qed

\begin{rem}\label{rem:normalized R for finite dim}
{\rm
We may recover the spectral decomposition on $\W_{\ov{\bm\e}}^{\sigma_1}(z_1) \, {\ot}\,\W_{\ov{\bm\e}}^{\sigma_2}(z_2)$, a tensor product of spin representations of $U_{\td{q}}(D_{m+1}^{(1)})$ in \cite{O90}.
For this, we need another normalization. Let
\begin{equation}\label{eq:normalized R matrix'}
 \mc{R}^\diamond_{{\bm\sigma},{\bm\e}}(z)=D_{{\bm\sigma}}\mc{R}_{{\bm\sigma},{\bm\e}}(z),
\end{equation}
where
\begin{equation*}
 D_{{\bm\sigma}}=
\begin{cases}
 \prod_{i=1}^d\frac{z-q^{4i-2}}{1-q^{4i-2}} & \text{if $\bm{\sigma}=(+,+), (-,-)$},\\
 \prod_{i=1}^{d'}\frac{z-q^{4i}}{1-q^{4i}} & \text{if $\bm{\sigma}=(+,-), (-,+)$},
\end{cases}
\end{equation*}
with $d=\lfloor \frac{m+1}{2} \rfloor$ and $d'=\lfloor \frac{m}{2} \rfloor$. 
Then by applying $\mf{tr}^{\bm\e}_{\ov{\bm\e}}$ to $\mc{R}^\diamond_{{\bm\sigma},{\bm\e}}(z)$ and Lemma \ref{lem:tr of norm R}, we recover the formula in \cite[Proposition 5.4]{O90}, which is equal to the renormalized $R$ matrix \cite{KKKO15} (up to scalar multiplication). Note that the normalization \eqref{eq:normalization} is compatible with the spectral decomposition on $\W_{\un{\bm\e}}^{\sigma_1}(z_1) \, {\ot}\,\W_{\un{\bm\e}}^{\sigma_2}(z_2)$ in \cite{KO} under $\mf{tr}^{\bm\e}_{\un{\bm\e}}$.
} 
\end{rem}

\subsection{Irreducible representations in $\widehat{\mc O}^{\,\mf c}_{\rm osc,\e}$}\label{subsec:fusion type c}
Let
\begin{equation}\label{eq:poles}
 \mc{Z}^{\bm\sigma}=
\begin{cases}
 \{\,q^{4a+2}\,|\,a\in\Z_+\,\} &  \text{for $\bm{\sigma}=(+,+), (-,-)$},\\
 \{\,q^{4a+4}\,|\,a\in\Z_+\,\} &  \text{for $\bm{\sigma}=(+,-), (-,+)$}\\
\end{cases}
\end{equation}
be the set of poles in the coefficient of $\mc P^{{\bm\sigma}}_{\la,{\bm\e}}$ in $\mc{R}_{{\bm\sigma},{\bm\e}}(z)$ for $\la\in S^{\bm\sigma}_{\bm\e}$.

Let $\mc{P}^+$ be the set of $({\mb \sigma},{\mb c})$ such that
\begin{itemize}
\item[(1)]  ${\mb \sigma}=(\sigma_{1},\dots,\sigma_{\ell})\in\{+,-\}^{\ell}$ and $\mb{c}=(c_1,\dots,c_\ell)\in (\Bbbk^\times)^\ell$ for some $\ell\geq 1$,

\item[(2)] $c_{i}/c_{j}\not\in \mc{Z}^{(\sigma_i,\sigma_j)}$ for all $i<j$ when $\ell\ge 2$.
\end{itemize}

Let $({\mb \sigma},{\mb c})\in \mc{P}^+$ be given with $\ell\ge 2$.
Let $z_{1},\,\dots,\,z_{\ell}$ be indeterminates.
Let $\W_{\e}^{\sigma_{i}}(z_{i})$
denote the affinization $(\W_{\e}^{\sigma_{i}})_{{\rm aff}}=\Bbbk[z_{i}^{\pm1}]\ot \W_{\e}^{\sigma_{i}}$ for an indeterminate $z_i$ ($1\le i\le\ell$).

Since the normalized $R$ matrix \eqref{eq:norm R} on $\W^{\sigma_i}_\e(z_i)\ot \W^{\sigma_j}_\e(z_j)$ for $i<j$ satisfies the Yang-Baxter equation, we may define a $\Bbbk[z_{1}^{\pm1},\dots,z_{\ell}^{\pm1}]\otimes\U_\e$-linear
map 
\begin{equation*}\label{eq:norm R_l} 
\xymatrixcolsep{2.5pc}\xymatrixrowsep{3pc}\xymatrix{ 
\mathcal{W}^{\sigma_1}_{\epsilon}(z_1)\otimes\cdots\otimes\mathcal{W}^{\sigma_\ell}_{\epsilon}(z_\ell) \ \ar@{->}^{\hskip-2cm\mathcal{R}_{{\boldsymbol \sigma},\e}}[r] & \ \Bbbk(z_1,\dots,z_\ell)\otimes_{\Bbbk[z^{\pm 1}_1,\dots,z^{\pm 1}_\ell]} \mathcal{W}^{\sigma_\ell}_{\epsilon}(z_\ell)\otimes\cdots\otimes \mathcal{W}^{\sigma_1}_{\epsilon}(z_1), 
}
\end{equation*}
by taking the composition of the normalized $R$ matrices with respect to the longest element of the symmetric group $\mf{S}_\ell$.
By Theorem \ref{thm:spectral decomp of W}, the specialization of $\mc{R}_{\bm{\sigma},\e}$ at $(z_{1},\dots,z_{\ell})=(c_{1},\dots,c_{\ell})$
gives a well-defined $\U_\e$-linear map
\begin{equation}\label{eq:R matrix for fusion}  
\xymatrixcolsep{2.5pc}\xymatrixrowsep{3pc}\xymatrix{ 
\mc{R}_{\bm{\sigma},\e}(\bm{c}): \mathcal{W}^{\sigma_1}_{\epsilon}(c_1)\otimes\cdots\otimes\mathcal{W}^{\sigma_\ell}_{\epsilon}(c_\ell) \ \ar@{->}[r] & \ \mathcal{W}^{\sigma_\ell}_{\epsilon}(c_\ell)\otimes\cdots\otimes \mathcal{W}^{\sigma_1}_{\epsilon}(c_1)}. 
\end{equation} 
We denote by $\W_{\e}({\mb \sigma},{\mb c})$ the image of $\mc{R}_{\bm{\sigma},\e}(\bm{c})$.

\begin{prop}\label{prop:truncation and simples}
For $({\mb \sigma},{\mb c})\in \mc{P}^+$, we have 
\begin{equation*}
 \mf{tr}^{\bm\e}_{\e}\left(\W_{\bm\e}({\mb \sigma},{\mb c})\right)
 =\W_{\e}({\mb \sigma},{\mb c})\quad (\e=\un{\bm\e}, \ov{\bm\e}).
\end{equation*}
\end{prop}
\pf Since $\mc{R}_{\bm{\sigma},{\bm\e}}$ is a composition of $\mc{R}_{(\sigma_i,\sigma_j),{\bm\e}}(z_i/z_j)$ for $i<j$, it follows from Lemma \ref{lem:tr of norm R} and \eqref{eq:W for e} that $\mf{tr}^{\bm\e}_{\e}\left(\W_{\bm\e}({\mb \sigma},{\mb c})\right)$ coincides with $\W_{\e}({\mb \sigma},{\mb c})$.
\qed\smallskip

\begin{thm}\label{thm:fusion-spin}
For $({\mb \sigma},{\mb c})\in \mc{P}^+$, $\W_{\e}({\mb \sigma},{\mb c})$ is irreducible if it is not zero.
\end{thm}
\pf We note from Theorem \ref{thm:spectral decomp of W} that
$\mc{R}_{(\pm,\pm),\e}(z_{1}/z_{2})= \rm{id}_{\W_\e^{\pm}(c)\ot \W_\e^{\pm}(c)}$
at $z_{1}=z_{2}=c$ for any $c\in\Bbbk^{\times}$.
Then the irreducibility can be proved in the same way as in \cite[Theorems~4.3, 4.4]{KL20} with the normalized $R$ matrices. In our case, $\mc{R}_{\bm{\sigma},\e}(z)$ is an infinite sum of $\mc P^{{\bm\sigma}}_{\la,{\bm\e}}$ unless $\e=\ov{\bm\e}, \ov{\bm\e}'$. Hence it is not rationally renormalizable in the sense of \cite{KKKO15}. But we may still apply the arguments there using the normalized $R$ matrix as in \cite[Lemma~2.6]{KK} (cf.~\cite[Appendix D]{KwO21}).
\qed\smallskip

The category $\widehat{\mc O}^{\,\mf c}_{\rm osc,\un{\bm\e}}$ consists of infinite-dimensional representations of type $C_m^{(1)}$, which is an affine analogue of $q$-oscillator representation of type $C_m$. On the other hand, $\widehat{\mc O}^{\,\mf c}_{\rm osc,\ov{\bm\e}}$ is a category of finite-dimensional representation of type $D_{m+1}^{(1)}$ (see Section \ref{subsec:fundamental rep of type D}).

We remark that a special class of irreducible representations in $\widehat{\mc O}^{\,\mf c}_{\rm osc,\un{\bm \e}}$ is constructed in \cite{KO}, which corresponds to Kirillov-Reshetikhin modules of type $D_{m+1}^{(1)}$ associated to spin representations under truncations (cf.~Proposition \ref{prop:truncation on affine osc}). Other than this, the nonzero irreducible representations $\W_{\e}({\mb \sigma},{\mb c})$ for $\e={\bm\e}, \un{\bm\e}$ are new.

\begin{prop}\label{prop:fusion-c-nonzero}
Let $(\bm{\sigma},\bm{c})\in \mc{P}^+$ be given such that $\bm{\sigma}^+ \in \mathcal{P}(O_{\ell})_{\e}$, where $\bm{\sigma}^+ = (1^s)$ is the partition with $s$ the number of $-$'s in $\bm\sigma$.
Then $\W_{\e}({\mb \sigma},{\mb c})$ is nonzero. In particular, $\W_{\e}({\mb \sigma},{\mb c})$ is nonzero if the length of $\epsilon$ is sufficiently large.
\end{prop}
\pf Let $(\bm{\sigma},\bm{c})\in \mc{P}^+$ be given where ${\mb \sigma}=(\sigma_{1},\dots,\sigma_{\ell})$ and $\mb{c}=(c_1,\dots,c_\ell)\in (\Bbbk^\times)^\ell$.

First consider the case when $\e=\ov{\bm\e}$. Then we have ${\bm\sigma}^+\in \mc{P}(O_{\ell})_{\overline{\boldsymbol{\epsilon}}}$ for any $m$.
Note that $\W^+_{\ov{\bm\e}}(x) \cong V(\upvarpi_{m+1})_x$ and $\W^-_{\ov{\bm\e}}(x)\cong V(\upvarpi_{m})_x$ for $x\in \Bbbk^\times$. Let $w_{i}$ be the $\mathring{\U}_{\ov{\bm\e}}$-highest weight vector of $\W^{\sigma_i}_{\ov{\bm\e}}(c_i)$ for $1\le i\le \ell$.
By Theorem \ref{thm:spectral decomp of W}, $\mc{R}_{\bm{\sigma},\overline{\boldsymbol{\e}}}(\bm{c})$ maps $w_{1}\ot \cdots\ot w_{\ell}$ to $w_{\ell}\ot \cdots\ot w_{1}$, which implies in particular $\mc{V}^{\bm{\sigma}^+}_{\ov{\bm\e}}\subset \W_{\ov{\bm\e}}({\mb \sigma},{\mb c})$ and hence $\W_{\ov{\bm\e}}({\mb \sigma},{\mb c})$ is nonzero.
 
If $\e={\bm\e}=\e^{(m)}$, then we also have ${\bm\sigma}^+\in \mc{P}(O_{\ell})_{\bm\e}$ for any $m$. Since $\mf{tr}^{\bm\e}_{\ov{\bm\e}}\left(\mc{V}^{\bm{\sigma}^+}_{\bm\e}\right)=\mc{V}^{\bm{\sigma}^+}_{\ov{{\bm\e}}}$, we conclude that $\W_{\bm\e}({\mb \sigma},{\mb c})$ is nonzero. Finally, suppose that $\e=\un{\bm\e}$ and ${\bm\sigma}^+\in \mc{P}(O_\ell)_{\un{\bm\e}}$.
By applying $\mf{tr}^{\bm\e}_{\un{\bm\e}}$ to $\mc{V}^{\bm{\sigma}^+}_{{\bm\e}}\subset \W_{{\bm\e}}({\mb \sigma},{\mb c})$, we have  
$0\neq \mc{V}^{\bm{\sigma}^+}_{\un{\bm\e}}\subset \W_{\un{\bm\e}}({\mb \sigma},{\mb c})$, which implies that $\W_{\un{\bm\e}}({\mb \sigma},{\mb c})$ is nonzero.
\qed

\subsection{Fundamental representations}\label{subsec:fundamental rep of type D}
Let  
\begin{equation}\label{eq:normalized R matrix' any e}
 \mc{R}^\diamond_{{\bm\sigma},{\e}}(z) :=
 \mf{tr}^{\bm\e}_{\e}
 \left( \mc{R}^\diamond_{{\bm\sigma},{\bm\e}}(z)\right),
\end{equation}
where $\mc{R}^\diamond_{{\bm\sigma},{\bm\e}}(z)$ is given in \eqref{eq:normalized R matrix'}.

Let $\mc{Z}^{\bm\sigma}$ and $\mc{P}^+$ be defined in the same way as in Section \ref{subsec:fusion type c} with respect to the set of poles appearing in $\mc{R}^\diamond_{{\bm\sigma},{\bm\e}}(z)$, and define a map $\mc{R}^\diamond_{\bm{\sigma},\e}(\bm{c})$ as in \eqref{eq:R matrix for fusion} for $\e={\bm\e}, \un{\bm\e}, \ov{\bm\e}$.
Then $\W^\diamond_{\e}({\mb \sigma},{\mb c})$, the image of $\mc{R}^\diamond_{\bm{\sigma},\e}(\bm{c})$, is irreducible if it is not zero, by the same argument as in the proof of Theorem \ref{thm:fusion-spin}.

When $\e=\ov{\bm\e}$, it is more natural to find a connection with finite-dimensional irreducible representations of ${\U}_{\ov{\bm\e}}=U_{\td{q}}(D_{m+1}^{(1)})$ in $\widehat{\mc O}^{\,\mf c}_{\rm osc,\ov{\bm\e}}$ in terms of $\W^\diamond_{\ov{\bm\e}}({\mb \sigma},{\mb c})$.
Let $V(\upvarpi_k)$ denote the fundamental representation of $U_{\td{q}}(D_{m+1}^{(1)})$ corresponding to $\upvarpi_{k}$ for $1\le k\le m+1$. First, we see that
\begin{equation}\label{eq:spins for type D}
 V(\upvarpi_{m+1})\cong \W^+_{\ov{\bm\e}}(1),\quad 
 V(\upvarpi_{m}) \cong \W^-_{\ov{\bm\e}}(1).
\end{equation}

By Remark~\ref{rem:normalized R for finite dim} and \cite[Proposition 3.1]{Ko}, we have for $1\leq k\leq m-1$
\begin{equation}\label{eq:fund of type D}
  V(\upvarpi_k)\cong \W^\diamond_{\ov{\bm\e}}({\mb \sigma}_{k},{\mb c}_k),
\end{equation}
where
\begin{equation*}
  ({\mb \sigma}_k,{\mb c}_k)=
  \begin{cases}
   ((+,+),(q^{2m-2k},1)) & \text{if $m+1-k$ is even},\\
   ((+,-),(q^{2m-2k},1)) & \text{if $m+1-k$ is odd},\\
  \end{cases}
\end{equation*}
(cf. Remark~\ref{rem:q-osc finite type}). Recall that any finite-dimensional irreducible $U_{\widetilde q}(D_{m+1}^{(1)})$-module $V$ 
with weight spaces as in \eqref{eq:weight space} is isomorphic to the image of $V(\upvarpi_{a_1})_{x_1}\ot\cdots\ot V(\upvarpi_{a_\ell})_{x_\ell}$ under the composition of the normalized $R$ matrices on $V(\upvarpi_{a_i})_{x_i}\ot V(\upvarpi_{a_j})_{x_j}$ ($i<j$) for some $a_1,\dots,a_\ell$ and  $x_1,\dots,x_\ell$ \cite{Kas02}. 
Therefore, it follows from the hexagon property of the universal $R$ matrix and \eqref{eq:spins for type D},\eqref{eq:fund of type D} that 
\begin{equation*}
V \cong \W^\diamond_{\ov{\bm\e}}({\mb \sigma}_V,{\mb c}_V)\cong \mf{tr}_{\ov{\bm\e}}\left(\W^\diamond_{{\bm\e}}({\mb \sigma}_V,{\mb c}_V)\right),
\end{equation*}
for some $({\mb \sigma}_V,{\mb c}_V)\in \mathcal{P}^+$ (defined with respect to the renormalization $\mathcal{R}^{\diamond}_{\boldsymbol{\sigma},\boldsymbol{\epsilon}}(z)$).

According to \eqref{eq:fund of type D}, one may regard $\W^\diamond_{{\bm\e}}\left({\mb \sigma}_k,{\mb c}_k\right)$ as a fundamental representation of $\U_{\bm\e}$. But it is only defined for $1\le k\le m+1$ and depends on $m$.
One may also consider another family of irreducible representation $\mathcal{W}_{l,\epsilon}=\mathcal{W}_{\epsilon}(\boldsymbol{\sigma}^{\prime}_l,\boldsymbol{c}_{l}^{\prime})$ of $\U_{\e}$, where
\begin{equation}\label{eq:def of fundamental}
  ({\mb \sigma}^{\prime}_l,{\mb c}^{\prime}_l)=
\begin{cases}
 ((+,+),(q^{-2l-2},1)) & \text{if $l$ is even}, \\
 ((+,-),(q^{-2l-2},1)) & \text{if $l$ is odd},
\end{cases}\quad (l\in \Z_+)
\end{equation}
which is well-defined by Theorem \ref{thm:spectral decomp of W}.
As a $\ring{\U}_\e$-module, $\W_{l,\e}$ decomposes as follows: 
\begin{equation*}
 \W_{l,\e} = \bigoplus_{\substack{i\ge 0 \\ (l-2i,0)\in \mc{P}(O_2)_\e}} \V^{(l-2i,0)}_\e.
\end{equation*}
Since $(0,0),(1,0)\in\mathcal{P}(O_2)_\epsilon$, $\mathcal{W}_{l,\epsilon}$ is always nonzero 
and hence irreducible by Theorem~\ref{thm:fusion-spin}.

For $x\in \Bbbk^\times$, we let $\W_{l,\e}(x)=(\W_{l,\e})_x$.
As in the case of usual fundamental representations, one can check the irreducibility of a tensor product of two representations, with generic spectral parameters, among \eqref{eq:def of fundamental} and $\mathcal{W}_{\epsilon}^{\sigma}$.

\begin{prop}
For $l_1,\,l_2\in\mathbb{Z}_+$ and $\sigma\in \{+,-\}$,
$$\mathcal{W}_{l_1,\epsilon}\otimes\mathcal{W}_{l_2,\epsilon}(a),\quad
\mathcal{W}_{\epsilon}^{\sigma}\otimes\mathcal{W}_{l_1,\epsilon}(a),\quad
\mathcal{W}_{l_1,\epsilon}(a)\otimes\mathcal{W}_{\epsilon}^{\sigma}$$
are irreducible as $\U_{\e}$-modules for generic $a\in\Bbbk^\times$.
\end{prop}
\pf We remove $\e$ for simplicity since the argument does not depend on $\e$. We prove that $\W_{l_1}\ot \W_{l_2}(a)$ is irreducible for even $l_1$ and odd $l_2$, when $a\notin q^{2\mathbb{Z}}$. The proof for the other cases are similar.

By definition, $\W_{l_1}\ot \W_{l_2}(a)$ is the
image of $R_1\ot R_2$, where
\begin{align*}
&\xymatrixcolsep{3pc}\xymatrixrowsep{3pc}\xymatrix{
R_1=\mathcal{R}_{(+,+)}(q^{-2l_1-2}) : \mathcal{W}^{+}(q^{-2l_1-2})\otimes\mathcal{W}^{+} \ \ar@{->}[r] & \mathcal{W}^{+}\otimes\mathcal{W}^{+}(q^{-2l_1-2})},\\
&\xymatrixcolsep{3pc}\xymatrixrowsep{3pc}\xymatrix{
R_2=\mathcal{R}_{(+,-)}(q^{-2l_2-2}) : \mathcal{W}^{+}(aq^{-2l_2-2})\otimes\mathcal{W}^{-}(a) \ \ar@{->}[r] & \mathcal{W}^{-}(a)\otimes\mathcal{W}^{+}(aq^{-2l_2-2})}.
\end{align*}
Take a reduced expression $w_{0}=s_{2}s_{3}s_{1}s_{2}s_{3}s_{1}$
of the longest element $w_{0}$ of $\mathfrak{S}_{4}$, where $s_i$ is the transposition $(i,i+1)$,
and consider the corresponding composition of normalized $R$ matrices
{\small
\begin{equation*}
\xymatrixcolsep{1.5pc}\xymatrixrowsep{3pc}\xymatrix{
\mathcal{W}^{+}(q^{-2l_1-2})\otimes\mathcal{W}^{+}\otimes\mathcal{W}^{+}(aq^{-2l_2-2})\otimes\mathcal{W}^{-}(a)  \ 
\ar@{->}[r]^{R} & 
\mathcal{W}^{-}(a)\otimes\mathcal{W}^{+}(aq^{-2l_2-2})\otimes\mathcal{W}^{+}\otimes\mathcal{W}^{+}(q^{-2l_1-2})},
\end{equation*}}
which is well-defined by Theorem \ref{thm:fusion-spin}. 
By our choice of the reduced expression, we have $R=R^{\prime}\circ(R_{1}\otimes R_{2})$ where $R^{\prime}$
corresponds to the factor $s_{2}s_{3}s_{1}s_{2}$ in $w_0$ and so obtain the following commutative diagram:
\begin{equation*}
\xymatrixcolsep{4pc}\xymatrixrowsep{3pc}\xymatrix{
\W^+(q^{-2l_1-2})\ot \W^{+}\ot \W^{+}(aq^{-2l_2-2})\ot \W^{-}(a)  \  \ar@{->}[r]^{\quad\quad\quad R_1\ot R_2}  \ar@{->}[d]^{R} & \ {\rm Im}(R_1 \ot R_2) = \W_{l_1}\ot \W_{l_2}(a) \ar@{->}[d]^{R'}  \\
\W^{-}(a)\ot \W^{+}(aq^{-2l_2-2})\ot \W^{+}\ot \W^{+}(q^{-2l_1-2})  & \  R'(\W_{l_1}\ot \W_{l_2}(a))= {\rm Im} R  \ar@{_{(}->}[l] }.  
\end{equation*}

By Theorem \ref{thm:fusion-spin}, ${\rm Im} R$ is either zero
or irreducible. On the other hand, we see that $R'$ is a composition of the specializations of the normalized $R$ matrices, which are all isomorphisms when $a\notin q^{2\mathbb{Z}}$ by Theorem \ref{thm:spectral decomp of W}. Hence $R^{\prime}$ is also an isomorphism.
Since $\mathrm{Im}(R_{1}\otimes R_{2})$
is nonzero and $R^{\prime}$ is an isomorphism, $\mathrm{Im}R=R^{\prime}(\mathrm{Im}(R_{1}\otimes R_{2}))$
is nonzero, and hence irreducible. This implies the irreducibility of
$\mathrm{Im}(R_{1}\otimes R_{2})=\mathcal{W}_{l}\otimes\mathcal{W}_{m}(q)$,
again since $R^{\prime}$ is an isomorphism.
\qed

\section{Irreducible representations in $\widehat{\mc O}^{\,\mf d}_{\rm osc,\e}$}\label{sec:irr of type d}
In this section, we construct a family of irreducible representations in $\widehat{\mc O}^{\,\mf d}_{\rm osc,\e}$ and explain a connection between $q$-oscillator representations of type $D_{m+1}^{(1)}$ and finite-dimensional irreducible representations of type $C_{m}^{(1)}$. 

The idea of construction of irreducible representations is similar as in the case of $\widehat{\mc O}^{\,\mf c}_{\rm osc,\e}$. But an additional and crucial step necessary in this case is to construct a fundamental type representation in $\widehat{\mc O}^{\,\mf d}_{\rm osc,\mb{\e}}$ corresponding to a finite-dimensional fundamental representation of type $C_m^{(1)}$ under truncation and to find the spectral decomposition on their tensor products.
We assume the notations in Sections \ref{subsec:truncation d} and \ref{subsec:fin osc category type x}.

\subsection{Fundamental representation $\W_{l,{\e}}(x)$}\label{subsec:fundamental representation of osc D}

Suppose first that $\e=\mb{\e}', \un{\mb{\e}}'$. 
For $l\in \Z_+$ and $0\le k\le l$, 
let 
\begin{equation*}
\begin{split}
 v_{l,k}&= \ket{k\mathbf{e}_{n}}\otimes\ket{(l-k)\mathbf{e}_{n}}\in\W_\e^{\ot 2},\\\
 \W_{l,k,\e}&= \mathring{\U}_\e v_{l,k} \subset \W_\e^{\ot 2}.
\end{split}
\end{equation*}
Then $\W_{l,k,\e}$ is a $\mathring{\U}_\e$-submodule of $\W_\e^{\ot 2}$, and $\W_{l,k,\e}\cong \V^{(l)}_\e$ by Lemma \ref{lem:semisimplicity of osc} since $v_{l,k}$ is a $\mathring{\U}_\e$-highest weight vector of weight $\La_{(l),\e}$.

\begin{lem}\label{lem:classical decomposition}
As a $\mathring{\U}_\e$-module, we have
\begin{equation*}
 \W_\e^{\ot 2}=\bigoplus_{l\in \Z_+}\bigoplus_{0\le k \le l}\W_{l,k,\e}.
\end{equation*}
\end{lem}
\pf 
By Corollary \ref{cor:decomp of W for e}, we have as a $\mathring{\U}_\e$-module
\begin{equation*}
 \W^{\ot 2}_{\e} \cong \bigoplus_{(l)\in \mc{P}(Sp_2)_\e}\ \left({\mc{V}_{\e}^{(l)}}\right)^{\oplus (l+1)},
\end{equation*}
where $\mc{P}(Sp_2)_\e=\{\,(l)\,|\,l\in \Z_+\,\}$ if $\e=\mb{\e}',\underline{\mb{\e}}'$.
Hence the decomposition follows from the fact that $\W_{l,k,\e}\cong \V^{(l)}_\e$.
\qed

\begin{lem}\label{lem:decomp-osc-A}
As a $\dring{\U}_\e$-module, we have 
\begin{equation*}
\mathcal{W}_{l,k,\e}=
\begin{cases}
\bigoplus_{j\geq0}\dring{\mathcal{U}}_\e\left(f_{m+1}^{j}v_{l,k}\right) & \text{if }\epsilon=\un{\bm{\e}}',\\
\dring{\mathcal{U}}_\e v_{l,k}\oplus\left(\bigoplus_{j\geq 1} \dring{\mathcal{U}}_\e\left(f_{n}(f_{n-2}f_{n})^{j-1}v_{l,k}\right)\right) & \text{if }\epsilon=\bm{\e}'.
\end{cases}
\end{equation*}
\end{lem}
\pf Note that
\[
\W_\e=\bigoplus_{s\ge 0}V_{\e}((s)),
\]
as a $\dring{\mathcal{U}}_\e$-module, and $\W_\e^{\ot 2}$ decomposes into a direct sum of $V_\e(\la)$'s according to the usual Littlewood-Richardson rule (cf. Section \ref{subsec:notation for osc}).

On the other hand, consider the following vectors in $\W_{l,k,\e}$:
\begin{equation*}\label{eq:maximal of type A}
\begin{cases}
f_{m+1}^{j-1}v_{l,k} & \text{if }\e=\un{\bm{\e}}',\\
v_{l,k},\,f_{n}(f_{n-2}f_{n})^{j}v_{l,k} & \text{if }\e =\bm{\e}',
\end{cases}
\quad (j\geq 1).
\end{equation*}
We can check easily that they are $\dring{\mathcal{U}}_\e$-highest weight vectors, and conclude that they give all the $\dring{\U}_\e$-highest weight vectors in $\W_\e^{\ot 2}$ from the decomposition of $\W_\e^{\ot 2}$. Hence there is no other $\dring{\U}_\e$-highest weight vector in $\mathcal{W}_{l,k,\e}$, which implies the decomposition of $\W_{l,k,\e}$. 
\qed

\begin{prop} Let $l\in\mathbb{Z}_+$, $0\leq k\leq l$, and $x\in \Bbbk^\times$ be given. Suppose that $\W_{l,k,\e}\subset \W_\e^{\ot 2}(x)$.
\begin{enumerate}
\item $\W_{l,k,\e}$ is invariant under $e_0$ and $f_0$, and hence is a $\U_\e$-module which we denote by $\W_{l,k,\e}(x)$. In particular  $\W_{l,k,\e}(x)$ is irreducible.
\item As $\U_\e$-modules, $\W_{l,k,\e}(x)$ are isomorphic to each other for any $0\leq k\leq l$.
\end{enumerate}
\end{prop}
\pf
Let us prove that $\W_{l,k}:=\W_{l,k,\boldsymbol{\epsilon}'}$ is invariant under $e_0$ and $f_0$.

First we claim that $e_0\W_{l,k}\subset \W_{l,k}$.
Since $\W_{l,k}$ is a $\mathring{\U}_{\boldsymbol{\epsilon}'}$-highest weight module generated by $v_{l,k}$ and $e_{0}$ commutes with $f_i$ for $i\in I_{\boldsymbol{\epsilon}'}\setminus\{0\}$, it suffices to show that $e_{0}v_{l,k}\in\mathcal{W}_{l,k}$.
It follows from the identity
\begin{equation}\label{eq:e_0v_{l,k}}
\begin{split}
x^{-1}e_{0}v_{l,k} & =\ket{\mathbf{e}_{1}+k\mathbf{e}_{n}}\otimes\ket{\mathbf{e}_{2}+(l-k)\mathbf{e}_{n}}-q\ket{\mathbf{e}_{2}+k\mathbf{e}_{n}}\otimes\ket{\mathbf{e}_{1}+(l-k)\mathbf{e}_{n}}\\
 & =\frac{1}{[l+1]} \big\{ (f_{2}f_{3}\cdots f_{n-2})f_{n-1}(f_{1}f_{2}\cdots f_{n-2})f_{n}v_{l,k}\\ 
 &\quad\quad\quad\quad\quad\quad\quad\quad\quad\quad -(f_{2}f_{3}\cdots f_{n-2})f_{n}(f_{1}f_{2}\cdots f_{n-2})f_{n-1}v_{l,k})\big\},
\end{split} 
\end{equation}
which can be checked in a straightforward way.

Next, we claim that $f_0\W_{l,k}\subset \W_{l,k}$.
Recall the decomposition of $\mathcal{W}_{l,k}$ in Lemma \ref{lem:decomp-osc-A}.
Each $\dring{\U}_{\boldsymbol{\e}'}$-component is finite-dimensional and irreducible, whose lowest weight vector for the $j$-th component (we regard $\dring{\mathcal{U}}_{\boldsymbol{\epsilon}'} v_{l,k}$ as the $0$-th component) is given by
\[
\begin{cases}
\ket{k\mathbf{e}_1}\otimes\ket{(l-k)\mathbf{e}_1} & \text{if }j=0,\\
(e_{0}e_{2})^{j-1}e_{0} \ket{k\mathbf{e}_1}\otimes\ket{(l-k)\mathbf{e}_1}  & \text{if }j\geq 1,
\end{cases}
\]
which indeed belongs to $\mathcal{W}_{l,k}$ by the first claim.
Now the claim reduces to proving that $f_0 v'_{l,k,j}\subset \W_{l,k}$, which 
can be seen without difficulty using the induction on $j$, together with relations
\[
f_{0}\ket{k\mathbf{e}_{1}}\otimes\ket{(l-k)\mathbf{e}_{1}}=0,\quad e_{0}f_{0}-f_{0}e_{0}=\frac{k_{0}-k_{0}^{-1}}{q-q^{-1}},\quad f_{0}e_{1}=e_{1}f_{0}.
\]
This proves (1).

Finally, the above arguments are already enough to prove (2). Indeed, for distinct $k,\,k'$, there is a $\mathring{\mathcal{U}}_{\boldsymbol{\epsilon}'}$-module isomorphism between $\mathcal{W}_{l,k}$ and $\mathcal{W}_{l,k'}$ under which the $\mathring{\mathcal{U}}_{\boldsymbol{\epsilon}'}$-highest weight vectors (resp. $\dring{\mathcal{U}}_{\boldsymbol{\e}'}$-lowest weight vectors) correspond. Since $e_0$-action on $v_{l,k}$ (\ref{eq:e_0v_{l,k}}) is expressed in terms of $f_1,\,\dots,\,f_{n}$, the isomorphism commutes with $e_0$. The commutativity with $f_0$ can be seen by inspecting the action of $f_0$ on $\dring{\mathcal{U}}_{\boldsymbol{\e}'}$-lowest weight vectors and induction on $j$ as above.

The proof for $\epsilon=\underline{\boldsymbol{\epsilon}}'$ is identical, taking the $j$-th lowest weight vector $e^{j}_{0}\ket{k\mathbf{e}_1}\otimes\ket{(l-k)\mathbf{e}_1}$.
\qed
\smallskip

From now on, we put
\begin{equation}
 \W_{l,\e}(x) = \W_{l,l,\e}(x) \quad (l\in \Z_+, x\in \Bbbk^\times),
\end{equation}
and call it the \emph{$l$-th fundamental representation} of $\U_\e$. We put $ \W_{l,\e}= \W_{l,\e}(1)$.

\begin{prop}\label{prop:truncation of fundamental of osc d}
 For $l\in \Z_+$ and $x\in \Bbbk^\times$, we have $\mf{tr}_{\un{\bm\e}'}(\W_{l,{\bm\e}'}(x))=\W_{l,\un{\bm\e}'}(x)$.
\end{prop}
\pf We have $\mf{tr}_{\un{\bm\e}'}(\W_{l,{\bm\e}'}(x))\supset \W_{l,\un{\bm\e}'}(x)$ by definition of $\W_{l,\e}$. The equality holds by Lemma \ref{lem:classical decomposition}.
\qed\smallskip

For $l\in \Z_+$ and $x\in \Bbbk^\times$, we put
\begin{equation}
 \W_{l,\ov{\bm\e}'}(x)=\mf{tr}_{\ov{\bm\e}'}(\W_{l,{\bm\e}'}(x)).
\end{equation}

\begin{prop}\label{prop:truncation of fundamental of osc d-2}
We have $\W_{l,\ov{\bm\e}'}(x)\neq 0$ if and only if $l\le m$. 
For $0\le l<m$, we have
\begin{equation*}
 \W_{l,\ov{\bm\e}'}(x)\cong V(\upvarpi_{m-l})_x, 
\end{equation*}
as a $U_{\td{q}}(C_m^{(1)})$-module, where $V(\upvarpi_{m-l})_x$ is the fundamental representation of $U_{\td{q}}(C_m^{(1)})$ corresponding to $\upvarpi_{m-l}$, and $\W_{m,\ov{\bm\e}'}(x)$ is a trivial representation.
\end{prop}
\pf It follows directly from the fact that $\W_{l,{\bm\e}'}\cong \V^{(l)}_{{\bm\e}'}$ and Theorem \ref{thm:trunc sends simple to simple}.
\qed\smallskip

\subsection{Irreducibility of $\W_{l_1,{\e}}(x_1) \ot \W_{l_2,{\e}}(x_2)$}
For $(l)\in \mc{P}(Sp_2)_\e$, let $\W_{l,\e}$ be the fundamental representation of $\U_\e$ given in Section \ref{subsec:fundamental representation of osc D}.
Recall that 
\begin{equation*}
 \mc{P}(Sp_2)_\e=
\begin{cases}
 \{\,(l)\,|\,l\in\mathbb{Z}_+\} & \text{for $\e={\mb \e}', \un{\mb \e}'$},\\
 \{\,(l)\,|\,0\le l\le m\,\} & \text{for $\ov{\mb \e}'$}.
\end{cases}
\end{equation*}
In this subsection, we prove the following.

\begin{thm}\label{thm:irreducibility of tensor product of fundamentals of osc D}
For $(l_1),(l_2)\in \mc{P}(Sp_2)_\e$ and generic $x_1, x_2\in \Bbbk^\times$, $\W_{l_1,{\e}}(x_1) \ot \W_{l_2,{\e}}(x_2)$ is an irreducible $\U_\e$-module.
\end{thm}

We first consider the case when $\e=\un{\mb \e}'$.
Note that ${\U}_{\un{\mb \e}'}=U_q(D_{m+1}^{(1)})$ and $\mathring{\U}_{\un{\mb \e}'}=U_q(D_{m+1})$.
For simplicity, put $\W_{l}=\W_{l,\un{\mb \e}'}$ for $l\in \Z_+$ and $\V^\la=\V^\la_{\un{\mb \e}'}$ for $\la\in \mc{P}(Sp_{2\ell})_{\un{\mb \e}'}$.

Let $l_1,l_2\in\Z_+$ be given.
We will prove the irreducibility of $\W_{l_1}(x_1)\ot \W_{l_2}(x_2)$ by finding all $\mathring{\U}_{\un{\mb \e}'}$-highest weight vectors in $\W_{l_1}(x_1)\ot \W_{l_2}(x_2)$ and then a relation between them under the action of $\U_{\un{\mb \e}'}$. Below we always regard $\mathcal{W}_{l}(x)=\mathcal{W}_{l,l,\underline{\boldsymbol{\e}}'}(x)\subset \mathcal{W}_{\underline{\boldsymbol{\e}}'}^{\otimes 2}(x)$.

First, we see from \eqref{eq:decom of two tensor of fundamentals} that
\begin{equation}\label{eq:decom of two tensor of fundamentals'}
 \W_{l_1}(x_1)\ot \W_{l_2}(x_2)=\V^{(l_1)}\ot \V^{(l_2)}\cong
 \bigoplus_{r\ge 0}\bigoplus_{0\leq s\leq\min(l_1,l_2)}\V^{(l_1+l_2+r-s,r+s)},
\end{equation}
which is a multiplicity-free decomposition.
Set 
\begin{align*}
v_{l_1} & =\ket{l_1\mathbf{e}_{n}}\ot \ket {\bf 0} \in\W_{l_1}(x_1),\quad 
v_{l_2}   =\ket{l_2\mathbf{e}_{n}}\ot \ket {\bf 0} \in\W_{l_2}(x_2),\\
u_{r,s}^{i,j} & =
f_{m+1}^{(i)}f_{m}^{(j)}v_{l_1}\ot 
f_{m+1}^{(r-i)}f_{m}^{(s-j)}v_{l_2}
\in\mathcal{W}_{l_1}(x_1)\otimes\mathcal{W}_{l_2}(x_2),
\end{align*}
for $0\leq i\leq r$ and $0\leq j\leq s$. 
Note that $v_{l_1}\otimes v_{l_2}$
is a $\mathring{\U}_{\un{\mb \e}'}$-highest weight vector of the component $\V^{(l_1+l_2,0)}$ in \eqref{eq:decom of two tensor of fundamentals'}.
%


\begin{lem}\label{lem:maximal vectors of osc D}
For $r\ge 0$ and $0\le s\le \min(l_1,l_2)$, the vector
\begin{equation*}
u_{r,s}=
\sum_{\substack{0\leq i\leq r\\ 0\leq j\leq s}}
\left[(-1)^{i+j}\prod_{k=1}^{i}q^{2k-2-l_2-2r}\frac{[r+l_2-k+2]}{[l_1+k+1]}\prod_{k^{\prime}=1}^{j}q^{2k^{\prime}+l_2-2s}\frac{[l_2-s+k^{\prime}]}{[l_1-k^{\prime}+1]}\right]u_{r,s}^{i,j}\label{eq:hw-vect-osc}
\end{equation*}
is a $\mathring{\U}_{\un{\mb \e}'}$-highest weight vector of $\V^{(l_1+l_2+r-s,r+s)}$ in $\mathcal{W}_{l_1}(x_1)\otimes\mathcal{W}_{l_2}(x_2)$,
where the coefficient of $u_{r,s}^{0,0}$ is understood to be 1.
\end{lem}
\pf
It is enough to show that $e_{m}u_{r,s}=e_{m+1}u_{r,s}=0$, which is straightforward from the commutation relation (\ref{comm-rel-e,f}).
\qed\smallskip

We also need the following lemmas which give relations between $\mathring{\U}_{\un{\mb \e}'}$-highest weight vectors in Lemma \ref{lem:maximal vectors of osc D}. The proofs can be found in the Appendices~\ref{sec:app-pf-lem-EF} and \ref{sec:app-pf-E-hw}.

\begin{lem}\label{lem:EF action on u^ij}
For $0\leq i\leq r$ and $0\leq j\leq s$, 
we have the following identities:
\begin{align*}
\boldsymbol{F}_{m}u_{0,s}^{0,j} & =
x_{1}q^{l_2 -2s+2j}[j+1]u_{0,s+1}^{0,j+1}+x_{2}[s-j+1]u_{0,s+1}^{0,j},\\
\boldsymbol{E}_{m}u_{0,s}^{0.j} & =x_{1}^{-1}[l_1-j+1]u_{0,s-1}^{0,j-1}+x_{2}^{-1}q^{2j-l_1}[l_2-s+j+1]u_{0,s-1}^{0,j},\\
\boldsymbol{F}_{m+1}u_{r,s}^{i,j} & =
x_{1}q^{2i-2r-l_2-2}[i+1]u_{r+1,s}^{i+1,j}+x_{2}[r-i+1]u_{r+1,s}^{i,j},\\
\boldsymbol{E}_{m+1}u_{r,s}^{i,j} & =-x_{1}^{-1}[l_1+i+1]u_{r-1,s}^{i-1,j}-x_{2}^{-1}q^{l_1+2i+2}[l_2+r-i+1]u_{r-1,s}^{i,j},
\end{align*}
where
\begin{align*}
\boldsymbol{F}_{m} & =(e_{m-1}\cdots e_{2}e_{1})e_{m+1}(e_{m-1}\cdots e_{3}e_{2})e_{0},\\
\boldsymbol{F}_{m+1} & =(e_{m-1}\cdots e_{2}e_{1})e_{m}(e_{m-1}\cdots e_{3}e_{2})e_{0},\\
\boldsymbol{E}_{m} & =f_{0}(f_{2}f_{3}\cdots f_{m-1})f_{m+1}(f_{1}f_{2}\cdots f_{m-1}),\\
\boldsymbol{E}_{m+1} & =f_{0}(f_{2}f_{3}\cdots f_{m-1})f_{m}(f_{1}f_{2}\cdots f_{m-1}).
\end{align*}
Here, we assume $u_{a,b}^{c,d}=0$ unless $0\leq c$, $0\leq d \leq \min(l_1,l_2)$, $0\leq a\leq c$ and $0\leq b \leq d$.
\end{lem}

\begin{lem}\label{lem:F-action-on-hwvectors}
For $r \ge 0$, $0 \le s \le \min\{l_1,l_2\}$ and generic $x_1,x_2\in\Bbbk^\times$, we have
\begin{align*}
\boldsymbol{F}_{m} u_{0,s} & \in \Bbbk^\times u_{0,s+1}+\Bbbk f_m u_{0,s} + \Bbbk f^{(2)}_m u_{0,s-1}, \\
\boldsymbol{F}_{m+1} u_{r,s} & \in \Bbbk^\times u_{r+1,s} + \Bbbk f_{m+1} u_{r,s} + \Bbbk f^{(2)}_{m+1} u_{r-1,s}.
\end{align*}
\end{lem}

\begin{thm}\label{thm:irreducibility of tensor product of fundamentals of osc D-1}
For $l_1,l_2\in \Z_+$ and generic $x_1,x_2\in\Bbbk^\times$, $\W_{l_1}(x_1) \ot \W_{l_2}(x_2)$ is an irreducible $U_{q}(D_{m+1}^{(1)})$-module.
\end{thm}
\pf
Take a nonzero $U_{q}(D_{m+1}^{(1)})$-submodule $N$ of $\W_{l_1}(x_1)\ot \W_{l_2}(x_2)$.
We claim that
\begin{enumerate}
\item $N$ contains $u_{0,0}=v_{l_1}\otimes v_{l_2}$,
\item $u_{0,0}$ generates $\W_{l_1}(x_1)\ot \W_{l_2}(x_2)$.
\end{enumerate}

We first prove (1). 
Since $\W_{l_1}(x_1)\ot \W_{l_2}(x_2)$ is semisimple over $U_{q}(D_{m+1})$, 
$N$ contains at least one $U_q(D_{m+1})$-highest weight vector, say $u_{r,s}$.

Suppose $r>0$ and let us show that $\boldsymbol{E}_{m+1}u_{r,s}$
is a nonzero multiple of $u_{r-1,s}$.
First, using Lemma~\ref{lem:maximal vectors of osc D} and Lemma~\ref{lem:EF action on u^ij}, we can explicitly compute the coefficient of $u_{r-1,s}^{0,0}$ in the expansion of $\boldsymbol{E}_{m+1}u_{r,s}$. 
For generic $x$ and $y$, that coefficient is nonzero and hence so is $\boldsymbol{E}_{m+1}u_{r,s}$.

%
Next we show that $\boldsymbol{E}_{m+1}u_{r,s}$ is a $U_{q}(D_{m+1})$-highest weight vector.
Indeed, for $p\neq m,\,m+1$, we have
$e_{p}\boldsymbol{E}_{m+1}u_{r,s}=0$ by considering its weight. 
For $p=m$, we have
\begin{align*}
e_{m}\boldsymbol{E}_{m+1}u_{r,s} & =f_{0}(f_{2}\cdots f_{m-1})\frac{k_{m}-k_{m}^{-1}}{q-q^{-1}}(f_{1}\cdots f_{m-1})u_{r,s}\\
 & =[l_1+l_2-2s+1]f_{0}(f_{2}\cdots f_{m-1})(f_{1}\cdots f_{m-1})u_{r,s}=0.
\end{align*}
Here the last equality follows from $f_0(f_2 \cdots f_{m-1})(f_1 \cdots f_{m-1})u_{r,s}^{i,j}=0$, which can be checked directly by expanding $u^{i,j}_{r,s}$ in terms of standard basis (see \eqref{eqn:vanishing-appendix} in Appendix \ref{sec:app-pf-lem-EF}).
For $p=m+1$, it is clear that $e_{m+1}\boldsymbol{E}_{m+1}u_{r,s}=0$. 
To sum up, we have $\boldsymbol{E}_{m+1}u_{r,s}\in \Bbbk^{\times}u_{r-1,s}$
and hence $u_{r-1,s}\in N$. 

Inductively, we obtain $u_{0,s}\in N$. 
Similarly one can see that $\boldsymbol{E}_m u_{0,s}$ is a nonzero multiple of $u_{0,s-1}$, ending up with $u_{0,0}\in N$.

To see (2) hold, it suffices to show that $u_{r,s}\in N$ for any $r,s$ such that  $0\leq r,\,0\leq s\leq \min(l_1,l_2)$.
This is done by induction on $r+s$. Namely, under the induction hypothesis on $r+s=k$ and genericity of $x$ and $y$, 
we obtain $u_{0,k+1}\in N$ by the first assertion, and $u_{r',s'}\in N$ for $r'+s'=k+1,\,r'\neq0$ by the second one of Lemma~\ref{lem:F-action-on-hwvectors}.
The proof is completed.

%
\qed\smallskip

{\em Proof of Theorem \ref{thm:irreducibility of tensor product of fundamentals of osc D}}. 
In case of $\e=\ov{\bm\e}'$, $\W_{l,{\e}}(x)$ is the usual fundamental representation for $0\le l<m$ by Proposition \ref{prop:truncation of fundamental of osc d-2}, and hence $\W_{l_1,{\e}}(x_1) \ot \W_{l_2,{\e}}(x_2)$ is irreducible for $l_1,l_2<m$ by \cite{Kas02}.

It remains to prove the case when $\e=\bm{\e}'$.
By Proposition \ref{prop:truncation of fundamental of osc d} and Proposition \ref{prop:truncation} for $\mf{tr}_{\un{\bm \e}'}$, we have
\begin{equation*}
 \mf{tr}_{\un{\bm \e}'}(\W_{l_1,{\bm\e}'}(x_1) \ot \W_{l_2,{\bm\e}'}(x_2))
 =\W_{l_1,\un{\bm\e}'}(x_1) \ot \W_{l_2,\un{\bm\e}'}(x_2).
\end{equation*}
Since $\W_{l_1,{\bm\e}'}(x_1) \ot \W_{l_2,{\bm\e}'}(x_2)$ has the same decomposition as in \eqref{eq:decom of two tensor of fundamentals'} with respect to $\V_{\bm\e}^{\la}$ as a $\mathring{\U}_{{\bm\e}'}$-module, its irreducibility as a $\U_{{\bm\e}'}$-module follows from Theorem \ref{thm:irreducibility of tensor product of fundamentals of osc D-1} and Theorem \ref{thm:trunc sends simple to simple}.  
\qed

\subsection{Spectral decomposition}

Let $(l_1),(l_2)\in \mc{P}(Sp_2)_\e$ be given and put ${\bm l}=(l_1,l_2)$.
We have
\begin{equation}
 \W_{l_1,\e}\ot \W_{l_2,\e}=\V^{(l_1)}_\e\ot \V^{(l_2)}_\e\cong
 \bigoplus_{\la\in S_{{\bm l},\e}}\V^{\la}_\e,
\end{equation}
where $S_{{\bm l},\e}\subset \mc{P}(Sp_{4})_\e$ is given by  
\begin{equation*}
\begin{split}
 S_{{\bm l},{\bm \e}'}&=S_{{\bm l},\un{\bm \e}'}=\{\,(l_1+l_2+r-s,r+s)\,|\,r\ge 0, 0\leq s\leq\min(l_1,l_2)\,\},\\
 S_{{\bm l},\ov{\bm\e}'}&=\{\,\la=(\la_1,\la_2)\in S_{{\bm l},{\bm\e}'}\,|\,\la_1\le m\,\},
\end{split}
\end{equation*}
due to \eqref{eq:decom of two tensor of fundamentals} and Theorem \ref{thm:trunc sends simple to simple}.

As in Section \ref{subsec:norm R matrix for osc C}, we apply the universal $R$ matrix \eqref{eq:univ R matrix} to 
$\W_{l_1,\e}(z_1)\, \widehat{\ot}\,\W_{l_2,\e}(z_2)$ 
and then normalize it so that it acts as the identity map on $\V_\e^{\la^0}$, where $\la^0=(\max(l_1,l_2),\min(l_1,l_2))$. 
Then we get a $\Bbbk[z_1^{\pm 1},z_2^{\pm 1}]\ot \U_\e$-linear map
\begin{equation}\label{eq:R matrix on two tensor of fundamentals of osc D}
\xymatrixcolsep{2pc}\xymatrixrowsep{3pc}\xymatrix{
\mc{R}_{{\bm l},\e}(z) : 
\W_{l_1,\e}(z_1)\, {\ot}\, \W_{l_2,\e}(z_2)\ \ar@{->}[r] &  \
\Bbbk(z_1,z_2)\ot_{\Bbbk[z_1^{\pm 1},z_2^{\pm 1}]}\left(\W_{l_2,\e}(z_2)\, {\ot}\,\W_{l_1,\e}(z_1)\right)
}
\end{equation}
such that $\mc{R}_{{\bm l},\e}(z)\big\vert_{\V_\e^{\la^0}}= {\rm id}_{\V_\e^{\la^0}}$, which is uniquely given 
since the proof of Theorem \ref{thm:irreducibility of tensor product of fundamentals of osc D} also proves that $\Bbbk(z_1,z_2)\ot_{\Bbbk[z_1^{\pm 1},z_2^{\pm 1}]}\left(\W_{l_1,\e}(z_1)\, {\ot}\,\W_{l_2,\e}(z_2)\right)$ is irreducible over $\Bbbk(z_1,z_2)\ot \U_\e$. 
By the same arguments as in Lemma \ref{lem:tr of norm R}, we have
\begin{lem}\label{lem:tr of norm R'}
Under the above hypothesis, we have
\begin{equation*}
  \mf{tr}^{\bm\e'}_{\e}\left(\mc{R}_{{\bm l},\bm{\e}'}(z)\right)
 =\mc{R}_{{\bm l},\e}(z)\quad (\e=\un{\bm\e}', \ov{\bm\e}').
\end{equation*} 
\end{lem}\smallskip

For $\la\in S_{{\bm l},\e}$, we define a $\mathring{\U}_\e$-linear map
$\mc P^{{\bm l}}_{\la,\e} : 
\W_{l_1,\e}\ot \W_{l_2,\e} \longrightarrow \W_{l_2,\e}\ot \W_{l_1,\e}$ 
as follows:

\begin{itemize}
 \item[(1)] Let $v_\la$ and $v'_\la$ be the $\mathring{\U}_{\un{\bm\e}'}$-highest weight vectors of $\mc{V}^{\la}_{\un{\bm\e}'}$ in $\W_{l_1,\un{\bm\e}'}\ot \W_{l_2,\un{\bm\e}'}$ and $\W_{l_2,\un{\bm\e}'}\ot \W_{l_1,\un{\bm\e}'}$ given in Lemma \ref{lem:maximal vectors of osc D}, respectively for $\la\in S_{{\bm l},\un{\bm\e}'}$.
 
 \item[(2)] Let 
 $\mc P^{{\bm l}}_{\la,\un{\bm\e}'}$ be the map determined by 
$\mathcal{P}^{{\bm l}}_{\la,\un{\bm\e}'}(v_\mu)=\delta_{\lambda,\mu}v^{\prime}_\lambda$ for $\mu\in S_{{\bm l},\un{\bm\e}'}$.
 
 \item[(3)] Let $\mc P^{{\bm l}}_{\la,{\bm\e}'}$ and $\mc P^{{\bm l}}_{\la,\ov{\bm\e}'}$ be the unique maps satisfying 
\begin{equation*}
 \mf{tr}^{\bm{\e}'}_{\un{\bm\e}'}
 \left(\mc P^{{\bm l}}_{\la,{\bm\e}'}\right)
 =\mc P^{{\bm l}}_{\la,\un{\bm\e}'},\quad 
 \mf{tr}^{\bm{\e}'}_{\ov{\bm\e}'}\left(\mc P^{{\bm l}}_{\la,{\bm\e}'}\right)
 =\mc P^{{\bm l}}_{\la,\ov{\bm\e}'}.
\end{equation*}
\end{itemize}

\begin{thm}\label{thm:spectral decomp of W of osc D}
We have the following:
\[
\mathcal{R}_{{\bm l},\e}(z)=\sum_{\la\in S_{{\bm l},\e}}\rho_{\la}^{\bm l}(z)\mc{P}_{\la,{\e}}^{\bm l},
\]
where $\rho_{\la}^{\bm l}(z)$ with $\la=(\la_1,\la_2)=(l_1+l_2+r-s,r+s)\in  S_{{\bm l},\e}$ is given by
\[
\rho_{\la}^{\bm l}(z)=D_{\la}^{\bm l}\prod_{k=1}^{r}\frac{1-zq^{l_2+l_1+2k+2}}{z-q^{l_2+l_1+2k+2}}\prod_{k^{\prime}=1}^{\min(l_1,l_2)-s}\frac{1-zq^{\left|l_2-l_1\right|+2k^{\prime}}}{z-q^{\left|l_2-l_1\right|+2k^{\prime}}}
\]
for some $D_{\la}^{\bm l}\in\Bbbk^{\times}$. Here we assume that the product in $\rho_{\la}^{\bm l}(z)$ is 1 when $r=0$ or $s=\min\{l_1,l_2\}$.
\end{thm}
\pf We consider the case when $\e=\un{\bm\e}'$. Then the other cases follow from Lemma \ref{lem:tr of norm R'} and Theorem \ref{thm:trunc sends simple to simple}.

For $r\ge 0$ and $0\le s\le \min(l_1,l_2)$, let $u_{r,s}$ be the $\mathring{\U}_{\un{\bm\e}'}$-highest weight vector of $\W_{l_1}(z_1)\ot \W_{l_2}(z_2)$ given in Lemma \ref{lem:maximal vectors of osc D}. 
Then we have
\begin{equation*}\label{eq:lem for computation of spectral decomp}
\begin{split}
\boldsymbol{E}_{m}u_{0,s} & =[l_2-s+1](z_2^{-1}q^{-l_1}-z_1^{-1}q^{2+l_2-2s})u_{0,s-1},\\
\boldsymbol{E}_{m+1}u_{r,s} & =[l_2+r+1](-z_2^{-1}q^{l_1+2}+z_1^{-1}q^{-l_2-2r})u_{r-1,s}. 
\end{split}
\end{equation*}
Indeed, we have seen in the above proof that $\boldsymbol{E}_{m+1}u_{r,s}$ is a $U_{q}(D_{m+1})$-highest weight vector and hence is a scalar multiple of $u_{r-1,s}$. 
Then the second identity follows from Lemma \ref{lem:EF action on u^ij} where $z_1$ and $z_2$ appear under the action of $f_0$, since only $\boldsymbol{E}_{m+1}u_{r,s}^{0,0}$ and $\boldsymbol{E}_{m+1}u_{r,s}^{1,0}$ contribute to $u_{r-1,s}^{0,0}$ among $\boldsymbol{E}_{m+1}u_{r,s}^{i,j}$. The proof for the other one is similar.

For brevity we write $\mc{R}=\mc{R}_{{\bm l},{\e}}(z)$.
Let $\la=(\la_1,\la_2)=(l_1+l_2+r-s,r+s)\in  S_{{\bm l},\e}$ given and write $\rho_{r,s}(z)=\rho_{\la}^{\bm l}(z)$.

We compute the both sides of $\mc{R}(\boldsymbol{E}_{m+1}u_{r,s})=\boldsymbol{E}_{m+1}\mc{R}(u_{r,s})$:
\begin{align*}
\mc{R}(\boldsymbol{E}_{m+1}u_{r,s}) & =\mc{R}([l_2+r+1](-z_2^{-1}q^{l_1+2}+z_1^{-1}q^{-l_2-2r})u_{r-1,s})\\
 & =[l_2+r+1](-z_2^{-1}q^{l_1+2}+z_1^{-1}q^{-l_2-2r})\mc{R}(u_{r-1,s})\\
 & =[l_2+r+1](-z_2^{-1}q^{l_1+2}+z_1^{-1}q^{-l_2-2r})\rho_{r-1,s}(z)u_{r-1,s}^{\prime},
\end{align*}
while
\begin{align*}
\boldsymbol{E}_{m+1}\mc{R}(u_{r,s}) & =\rho_{r,s}(z)\boldsymbol{E}_{m+1}u_{r,s}^{\prime}\\
 & =\rho_{r,s}(z)[l_1+r+1](-z_1^{-1}q^{l_2+2}+z_2^{-1}q^{-l_1-2r})u_{r-1,s}^{\prime},
\end{align*}
and hence we obtain
\begin{align}\label{eq:recursion-sptrCoeff-r}
\rho_{r,s}(z) & =\frac{[l_2+r+1]}{[l_1+r+1]}\cdot\frac{-z_2^{-1}q^{l_1+2}+z_1^{-1}q^{-l_2-2r}}{-z_1^{-1}q^{l_2+2}+z_2^{-1}q^{-l_1-2r}}\rho_{r-1,s}(z) \nonumber\\
 & =q^{l_1-l_2}\frac{[l_2+r+1]}{[l_1+r+1]}\cdot\frac{1-zq^{l_1+l_2+2r+2}}{z-q^{l_1+l_2+2r+2}}\rho_{r-1,s}(z).
\end{align}
Similarly we get
\begin{equation}\label{eq:recursion-sptrCoeff-s}
\rho_{0,s}(z)=q^{l_2-l_1}\frac{[l_2-s+1]}{[l_1-s+1]}\cdot\frac{1-zq^{-l_1-l_2-2+2s}}{z-q^{-l_1-l_2-2+2s}}\rho_{0,s-1}(z).
\end{equation}
Now the desired formula follows from the normalization $\rho_{0,\min(l_1,l_2)}(z)=1$.
\qed

\begin{rem}{\rm 
The constants $D_{\la}^{\bm l}$ can be computed
explicitly using \eqref{eq:recursion-sptrCoeff-r} and \eqref{eq:recursion-sptrCoeff-s}.
We should remark that an explicit formula for the spectral decomposition of $\mathcal{R}_{{\bm l},\e}(z)$ even for $\e=\ov{\bm{\e}}'$ is new.
}
\end{rem}

\subsection{Irreducible representations in $\widehat{\mc O}^{\,\mf d}_{\rm osc,\e}$}\label{subsec:fusion type d}
For ${\bm l}=(l_1,l_2)\in \Z_+^2$, let
\begin{equation}\label{eq:poles}
 \mc{Z}^{\mb{l}}=
\{\,q^{|l_2-l_1|+2k}\,|\, k\ge 1\,\} \cup \{\,q^{l_2+l_1+2k+2}\,|\, k\ge 1\,\}
\end{equation}
be the set of poles in the coefficient of $\mc{P}_{\la,{\bm\e}}^{\bm l}$ in $\mathcal{R}_{{\bm l},{\bm\e}}(z)$ for $\la\in  S_{{\bm l},{\bm\e}}$.

Let $\mc{P}^+$ be the set of $({\mb l},{\mb c})$ such that
\begin{itemize}
\item[(1)]  ${\mb l}=(l_{1},\dots,l_{\ell})\in\Z_+^{\ell}$ and $\mb{c}=(c_1,\dots,c_\ell)\in (\Bbbk^\times)^\ell$ for some $\ell\geq 1$,

\item[(2)] $c_{i}/c_{j}\not\in \mc{Z}^{(l_i,l_j)}$ for all $i<j$ when $\ell\ge 2$.
\end{itemize}

Let $({\mb l},{\mb c})\in \mc{P}^+$ be given with $\ell\ge 2$.
Thanks to Theorem \ref{thm:spectral decomp of W}, we may define 
a well-defined $\U_\e$-linear map
\begin{equation}\label{eq:R matrix for fusion'}  
\xymatrixcolsep{2.5pc}\xymatrixrowsep{3pc}\xymatrix{ 
\mc{R}_{\bm{l},\e}(\bm{c}): \mathcal{W}_{l_1,\epsilon}(c_1)\otimes\cdots\otimes\mathcal{W}_{l_\ell,\epsilon}(c_\ell) \ \ar@{->}[r] & \ \mathcal{W}_{l_\ell,\epsilon}(c_\ell)\otimes\cdots\otimes \mathcal{W}_{l_1,\epsilon}(c_1)},
\end{equation}
as in \eqref{eq:R matrix for fusion}, whose image we denote by $\W_{\e}({\mb l},{\mb c})$.
Then we obtain the following as in Proposition \ref{prop:truncation and simples} and Theorem \ref{thm:fusion-spin}. 

\begin{prop}\label{lem:truncation and simples'}
For $({\mb l},{\mb c})\in \mc{P}^+$, we have 
\begin{equation*}
 \mf{tr}^{\bm\e'}_{\e}\left(\W_{\bm\e}({\mb l},{\mb c})\right)
 =\W_{\e}({\mb l},{\mb c})\quad (\e=\un{\bm\e}', \ov{\bm\e}').
\end{equation*}
\end{prop}
 
\begin{thm}\label{thm:fusion-fudamental}
For $({\mb l},{\mb c})\in \mc{P}^+$, $\W_{\e}({\mb l},{\mb c})$ is irreducible if it is not zero.
\end{thm}

The category $\widehat{\mc O}^{\,\mf d}_{\rm osc,\un{\bm\e}'}$ consists of infinite-dimensional representations of type $D_{m+1}^{(1)}$, while $\widehat{\mc O}^{\,\mf c}_{\rm osc,\ov{\bm\e}'}$ is a category of finite-dimensional representation of type $D_{m+1}^{(1)}$ and $\W_{\ov{\bm\e}'}({\mb l},{\mb c})$ gives all the finite-dimensional irreducible representations in $\widehat{\mc O}^{\,\mf c}_{\rm osc,\ov{\bm\e}'}$ by Proposition \ref{prop:truncation of fundamental of osc d-2}, Theorem \ref{thm:spectral decomp of W of osc D}, and \cite{Kas02}. We remark that the nonzero irreducible representations $\W_{\e}({\mb l},{\mb c})$ for $\e={\bm\e}', \un{\bm\e}'$ are new.
 
\begin{prop}\label{prop:fusion-d-nonzero}
Let $(\boldsymbol{l},\boldsymbol{c})\in \mc{P}^+$ be given such that $\boldsymbol{l}^+ \in \mathcal{P}(Sp_{2\ell})_{\epsilon}$, where $\boldsymbol{l}^+$ denotes the nonincreasing rearrangement of $\boldsymbol{l}$. 
Then $\W_{\e}({\mb l},{\mb c})$ is nonzero. In particular, $\W_{\e}({\mb l},{\mb c})$ is nonzero for all sufficiently large $m$.
\end{prop}
\pf  First suppose that $\epsilon=\overline{\boldsymbol{\epsilon}}'$. The assumption $\boldsymbol{l}^+ \in\mathcal{P}(Sp_{2\ell})_{\overline{\boldsymbol{\epsilon}}'}$ ensures that $l_i \leq m$ for all $i$, so that $\mathcal{W}_{l_i,\ov{\bm\e}}(c_i)\cong V(\upvarpi_{m-l_i})_{c_i}$ by Proposition~\ref{prop:truncation of fundamental of osc d-2}. If we let $w_i$ be the dominant extremal weight vector of $V(\upvarpi_{m-l_i})_{c_i}$, then $\mathcal{R}_{(l_{i}, l_{j}),\overline{\boldsymbol{\epsilon}}'}(c_{i}/c_{j})$ is normalized so that it sends $w_i \otimes w_j$ to $w_j \otimes w_i$ by Theorem \ref{thm:spectral decomp of W of osc D}. 
Hence $\mathcal{R}_{\boldsymbol{l},\overline{\boldsymbol{\epsilon}}'}(\boldsymbol{c})$ maps $w_1 \otimes \cdots \otimes w_\ell$ to $w_\ell \otimes \cdots \otimes w_1$, and $\W_{\bm{\e}'}({\mb l},{\mb c})$ is nonzero.

For the cases $\epsilon=\boldsymbol{\epsilon}',\,\underline{\boldsymbol{\epsilon}}'$, let us take $\widetilde{\epsilon}={\bm\e'}^{(\widetilde{m})}$ and $\overline{\widetilde{\epsilon}}=(1^{\widetilde{m}})$, for sufficiently large $\widetilde{m}$ so that $\boldsymbol{l}^{+}\in \mathcal{P}(Sp_{2\ell})_{\overline{\widetilde{\epsilon}}}$. By the above argument, $\W_{\ov{\widetilde{\epsilon}}}({\mb l},{\mb c})$ contains a $\mathring{\U}_{\ov{\widetilde{\epsilon}}}$-submodule generated by $w_\ell \otimes \cdots \otimes w_1$, namely $\mathcal{V}^{\boldsymbol{l}^{+}}_{\ov{\widetilde{\epsilon}}}$. 
This implies that 
$\mathcal{V}^{\boldsymbol{l}^{+}}_{\widetilde{\epsilon}}\subset 
\W_{\widetilde{\epsilon}}(\boldsymbol{l},\boldsymbol{c})$. Now applying $\mf{tr}^{\widetilde{\epsilon}}_{\e}$ for $\epsilon=\boldsymbol{\epsilon}',\,\underline{\boldsymbol{\epsilon}}'$, we conclude that $\W_{\e}(\boldsymbol{l},\boldsymbol{c})$ is nonzero.
\qed
 
\section{Remarks}\label{sec:remarks}
We give some remarks on generalizations of the results in this paper, and related problems.

\subsection{Duality}\label{subsec:duality}
Recall the following diagram:
\begin{equation}\label{eq:triangle of truncation osc x}
\xymatrixcolsep{3pc}\xymatrixrowsep{0.3pc}\xymatrix{
 & \widehat{\mc O}^{\,\mf x}_{\rm osc,{\e}}  \ar@{->}_{\mf{tr}^\e_{\un{\e}}}[dl]\ar@{->}^{\mf{tr}^\e_{\ov{\e}}}[dr] &  \\
 \widehat{\mc O}^{\,\mf x}_{\rm osc,\un{\e}} & &  \widehat{\mc O}^{\,\mf x}_{\rm osc,\ov{\e}}
}
\end{equation}
where $(\e,\un{\e},\ov{\e})=({\bm\e},\un{\bm\e},\ov{\bm\e})$ and $({\bm\e}',\un{\bm\e}',\ov{\bm\e}')$ for $\mf{x}=\mf{c}$ and $\mf{d}$, respectively (Proposition \ref{prop:truncation on affine osc}). 

Let ${\mb \la}=(\boldsymbol{l},\boldsymbol{c})\in \mc{P}^+$ be given, where $\mc{P}^+$ is defined in Sections \ref{subsec:fusion type c} and \ref{subsec:fusion type d}.
Then $\mf{tr}^{\e}_{\un\e}$ and $\mf{tr}^{\e}_{\ov\e}$ send the nonzero irreducible $\W_{\e}({\mb \la})$ to an irreducible representation or zero: 
\begin{equation*}\label{eq:triangle of truncation osc x with simples}
\xymatrixcolsep{3pc}\xymatrixrowsep{0.8pc}\xymatrix{
 & \W_{\e}({\mb \la})  \ar@{|->}_{\mf{tr}^\e_{\un{\e}}}[dl]\ar@{|->}^{\mf{tr}^\e_{\ov{\e}}}[dr] &  \\
 \W_{\un{\e}}({\mb\la}) & &  \W_{\ov{\e}}({\mb \la})
}
\end{equation*}
where  $\W_{\un{\e}}({\mb\la})$ and $\W_{\ov{\e}}({\mb \la})$ are interpolated by $\W_{\e}({\mb \la})$ when they are nonzero. 
Therefore, this provides an affine analogue of the correspondence between $q$-oscillator representations of $\mf{sp}_{2n}$ (resp. $\mf{so}_{2n}$) and finite-dimensional representations of $\mf{so}_{2n}$ (resp. $\mf{sp}_{2n}$) via Howe duality of the form $(\mf{g}, O_\ell)$ (resp. $(\mf{g}, Sp_{2\ell})$) in Sections \ref{subsec:Howe duality} and \ref{subsec:Howe duality-2}.

Given ${\bm\la}\in \mc{P}^+$, $\W_{\e^*}({\mb \la})$ ($\e^*=\e,\un{\e},\ov{\e}$) is nonzero for all sufficiently large $m$ by Propositions \ref{prop:fusion-c-nonzero} and \ref{prop:fusion-d-nonzero}, and the multiplicity of $\mc{V}_{\e^*}^{\la}$ in $\W_{\e^*}({\mb \la})$ for $\la\in \mc{P}(G_\ell)$ has a stable limit as $m\rightarrow \infty$ by Theorem \ref{thm:trunc sends simple to simple}(2). 
Hence we may consider a suitable limit of 
$\widehat{\mc O}^{\,\mf x}_{\rm osc,{\e^*}}$, where the correspondence between the limits of $\W_{\e^*}({\mb \la})$ becomes bijective in the corresponding diagram \eqref{eq:triangle of truncation osc x}.

The correspondence between oscillator representations of $\mf{sp}_{2n}$ (resp. $\mf{so}_{2n}$) and finite-dimensional representations of $\mf{so}_{2n}$ (resp. $\mf{sp}_{2n}$) via Howe duality can also be explained from the restriction of an equivalence called {\em super duality} \cite{CLW} between bigger categories including those representations. So it is natural to ask whether there is an equivalence between these limits of $\widehat{\mc O}^{\,\mf x}_{\rm osc,{\e^*}}$. This would imply that various important properties of $\widehat{\mc O}^{\,\mf x}_{\rm osc,{\ov\e}}$ also hold in $\widehat{\mc O}^{\,\mf x}_{\rm osc,{\e}}$ and $\widehat{\mc O}^{\,\mf x}_{\rm osc,{\un\e}}$ including the structure of the Grothendieck ring, the existence of $T$-systems and so on.

\subsection{Other cases than ${\bm \e}$ and ${\bm\e}'$}\label{subsec:general e} 
Let us give remarks on the generalization of \eqref{eq:triangle of truncation osc x} to the case of $\e=(\e_1,\dots,\e_n)$ other than ${\bm \e}$ and ${\bm\e}'$. We assume that $s=|\{\,i\,|\,\e_i=0\,\}|\ge 4$ and $t=|\{\,i\,|\,\e_i=1\,\}|\ge 4$.
Let $\un{\e}$ (resp. $\ov{\e}$) be the subsequence of $\e$ obtained by removing all $\e_i=1$ (resp. $\e_i=0$).

First, consider the case when $\e_1=\e_n$. 
We assume that $\mf{x}=\mf{c}$ for $\e_1=\e_n=1$ and $\mf{x}=\mf{d}$ for $\e_1=\e_n=0$. Following the arguments in Section \ref{subsec:fin osc category type x}, we can also define ${\mc O}^{\,\mf x}_{\rm osc,\e}$ and $\widehat{\mc O}^{\,\mf x}_{\rm osc,\e}$ using $\WW_\e(x)$ in \eqref{eq:WW(x)}, which is defined in Proposition \ref{prop:osc W for e} for $\mf{x}=\mf{c}$ and Proposition \ref{prop:osc W2 for e} for $\mf{x}=\mf{d}$.
The results in Sections \ref{subsec:fin osc category type x}, \ref{subsec:fusion type c} and \ref{subsec:fusion type d} also hold in the case of $\widehat{\mc O}^{\,\mf x}_{\rm osc,\e}$. Indeed, in order to prove them, we may use the following diagram:
\begin{equation*}\label{eq:diagram for general e}
\begin{tikzcd}
& \widehat{\mc O}^{\,\mf x}_{\rm osc,{{\tilde\e}}} \arrow[ld,"\mf{tr}^{{\tilde\e}}_{\e}"'] \arrow[rd, "\mf{tr}^{{\tilde\e}}_{{\tilde\e}^*}"] & \\
\widehat{\mc O}^{\,\mf x}_{\rm osc,\e}  \arrow[rd,"\mf{tr}^{\e}_{\e^*}"'] & & \widehat{\mc O}^{\,\mf x}_{\rm osc,{{\tilde\e}}^*} \arrow[ld,"\mf{tr}^{{{\tilde\e}}^*}_{\e^*}"]  \\
& \widehat{\mc O}^{\,\mf x}_{\rm osc,{\e}^*} &
\end{tikzcd}
\end{equation*}
where we take ${\tilde\e}={\bm\e}$ and ${\bm\e'}$ with sufficiently large $m$ for $\mf{x}=\mf{c}$ and $\mf{d}$, respectively, such that $\e$ is a subsequence of ${\tilde\e}$, and let $({\tilde\e}^*,{\e}^*)=(\un{{\tilde\e}},\un{\e})$ or $(\ov{{\tilde\e}},\ov{\e})$. Since $\e$ is not necessarily an alternating sequence, one may need to use the truncation of type $A$ in a general case of $\e$ (\cite[Theorem 4.3]{KY}) to define the above diagram.

Next consider the case when $\e_1\neq \e_n$. Let us assume that $\e_1=1$ and $\e_n=0$. Let $\WW_\e(x)=\W_\e^{\ot 2}(x)$, where we define the actions of $x_0$ ($x=e,f$) following Proposition \ref{prop:osc W for e} and the actions of $x_n$ ($x=e,f$) following Proposition \ref{prop:osc W2 for e}. Then we can define a semisimple tensor category ${\mc O}^{\,\mf d}_{\rm osc,{\e}}$ where an irreducible representation is parametrized by $\la\in \mc{P}(Sp_{2\ell})_\e$ $(\ell\ge 1)$, and then define $\widehat{\mc O}^{\,\mf d}_{\rm osc,{\e}}$ in a similar way.

We apply Propositions \ref{prop:reduction} and \ref{prop:reduction'} (and \cite[Theorem 4.3]{KY}) to the subsequences of $\e$ including its left and right (right and left) ends, respectively to have  homomorphisms of $\Bbbk$-algebras $\phi_{\ov{\e}} : \U_{\ov{\e}} \longrightarrow \U_{\e}$ and $\phi_{\un{\e}} : \U_{\un{\e}} \longrightarrow \U_{\e}$, where $\U_{\ov{\e}}$ is isomorphic to $U_{\td{q}}(A_{2t-1}^{(2)})$, and $\U_{\un{\e}}$ is isomorphic to $U_{q}(A_{2s-1}^{(2)})$ with the indices of simple roots arranged in a reverse way. Hence the diagram \eqref{subsec:duality} would interpolate the finite-dimensional representations and $q$-oscillator representations of type $A^{(2)}_{2s-1}$ and $A^{(2)}_{2t-1}$.

In order to have an analogue of the results in Section \ref{sec:irr of type d}
for $\widehat{\mc O}^{\,\mf d}_{\rm osc,{\e}}$, we first have to construct an irreducible representation in $\widehat{\mc O}^{\,\mf d}_{\rm osc,{\e}}$ which corresponds to a fundamental representation  $V(\upvarpi_k)_x$ of type $A_{2t-1}^{(2)}$ for $1\le k\le t$, and then compute the spectral decomposition on $V(\upvarpi_k)_{x}\ot V(\upvarpi_l)_y$, which is not known in general so far. The case when $\e_1=0$ and $\e_n=1$ is similar (see Section \ref{subsec:conjectures}).

\subsection{Quantum affine superalgebras of other types}\label{subsec:conjectures}
For a given $\e=(\e_1,\dots,\e_n)$, there exists a generalized quantum group of finite type $B$, $C$, $D$ \cite{Ma} (cf.~\cite{Ya94}). 
Using the defining relations of them, we can define the algebra $\U_{X\!Y}(\e)$ of  an affine type for $X,Y\in \{B,C,D\}$, where $X$ and $Y$ denote the types of the head and the tail in the associated Dynkin diagram, respectively.  
In particular, $\U_{D\!D}(\e)$ means $\U_D(\e)$ in Section \ref{subsec:def}.

We expect an analogue of the results in Sections \ref{sec:irr of type c} or \ref{sec:irr of type d} with respect to $\U_{X\!Y}(\e)$, which yields in particular various connections between finite-dimensional and $q$-oscillator representations of two quantum affine algebras. 
From the viewpoint of super duality \cite{CLW}, the correspondence between the Dynkin diagrams of finite type should be

{\small 
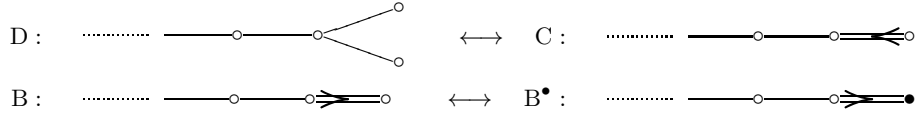
\begin{alignat*}{2}
{\rm D}:\quad 
&\vcenter{\xymatrix@R=1ex{
&&&& *{\circ}<3pt>\ar@{-}[dl]_<{} & \\
 \ar@{.}[r] 
&  \ar@{-}[r]   
&*{\circ}<3pt> \ar@{-}[r] 
&*{\circ}<3pt> 
&  \\
&&&&  *{\circ}<3pt>\ar@{-}[ul]_<{} & 
 }} & \!\!\!\!\!\!\!\!\!\! \longleftrightarrow \quad
{\rm C}: \quad 
&
\vcenter{\xymatrix@R=1ex{
 {} \ar@{.}[r] 
&  \ar@{-}[r]_<{} 
&*{\circ}<3pt> \ar@{-}[r]_<{} 
&*{\circ}<3pt> \ar@{=}[r] |-{\scalebox{2}{\object@{<}}}_<{}  
&*{\circ}<3pt>\ar@{}   
}}\\
{\rm B}:\quad
&
\vcenter{\xymatrix@R=1ex{
 {} \ar@{.}[r] 
&  \ar@{-}[r]_<{} 
&*{\circ}<3pt> \ar@{-}[r]_<{} 
&*{\circ}<3pt> \ar@{=}[r] |-{\scalebox{2}{\object@{>}}}_<{}  
&*{\circ}<3pt>\ar@{}   
}} 
& \!\!\!\!\!\!\!\!\!\! \longleftrightarrow \quad
{\rm B}^\bullet:\quad
&
\vcenter{\xymatrix@R=1ex{
 {} \ar@{.}[r] 
&  \ar@{-}[r]_<{} 
&*{\circ}<3pt> \ar@{-}[r]_<{} 
&*{\circ}<3pt> \ar@{=}[r] |-{\scalebox{2}{\object@{>}}}_<{}  
&*{\bullet}<3pt>\ar@{}   
}}
\end{alignat*}
}
where the Dynkin diagram of type $B$ corresponds to that of orthosymplectic Lie superalgebras without an odd isotropic root.

Therefore, the expected pairs of quantum affine (super)algebras without odd isotropic simple roots, which have a correspondence between ``finite-dimensional" and ``oscillator" representations or vice versa, are as follows:
\smallskip

{\small
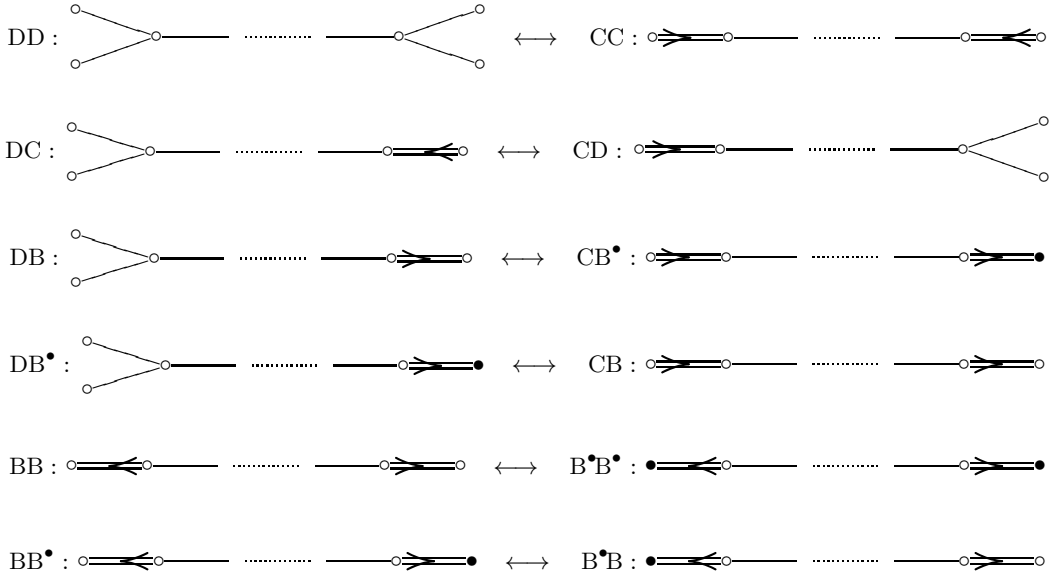
\begin{alignat*}{2}
{\rm DD}:\ 
&\vcenter{\xymatrix@R=1ex{
*{\circ}<3pt> \ar@{-}[dr]^<{}&&&&&*{\circ}<3pt> \ar@{-}[dl]^<{}\\
& *{\circ}<3pt> \ar@{-}[r]_<{} 
& {} \ar@{.}[r]&{} \ar@{-}[r]_>{} &
*{\circ}<3pt> & \\
*{\circ}<3pt> \ar@{-}[ur]_<{}&&&&&
*{\circ}<3pt> \ar@{-}[ul]_<{}}} & \quad \longleftrightarrow \quad 
{\rm CC}:\ 
&\vcenter{\xymatrix@R=1ex{
*{\circ}<3pt> \ar@{=}[r] |-{\scalebox{2}{\object@{>}}}_<{} 
&*{\circ}<3pt> \ar@{-}[r]_<{} 
& {} \ar@{.}[r]&{}  \ar@{-}[r]_>{} &
*{\circ}<3pt> \ar@{=}[r] |-{\scalebox{2}{\object@{<}}}
& *{\circ}<3pt>\ar@{}_<{}}}
\end{alignat*}
}\smallskip
{\small
\begin{alignat*}{2}
{\rm DC}:\ 
&\vcenter{\xymatrix@R=1ex{
*{\circ}<3pt> \ar@{-}[dr]^<{} \\
& *{\circ}<3pt> \ar@{-}[r]_<{} 
& {} \ar@{.}[r]&{}  \ar@{-}[r]_>{} &
*{\circ}<3pt> \ar@{=}[r] |-{\scalebox{2}{\object@{<}}}& *{\circ}<3pt>\ar@{}_<{} \\
*{\circ}<3pt> \ar@{-}[ur]_<{}}} & \quad \longleftrightarrow \quad
{\rm CD}:\ 
&\vcenter{\xymatrix@R=1ex{
 &&&&&*{\circ}<3pt> \ar@{-}[dl]^<{}\\
*{\circ}<3pt> \ar@{=}[r] |-{\scalebox{2}{\object@{>}}}_<{}& *{\circ}<3pt> \ar@{-}[r]_<{} 
& {} \ar@{.}[r]&{} \ar@{-}[r]_>{} &
*{\circ}<3pt> & \\
 &&&&&
*{\circ}<3pt> \ar@{-}[ul]_<{}}}%
\end{alignat*}
}
\smallskip
{\small
\begin{alignat*}{2}
{\rm DB}:\ 
&\vcenter{\xymatrix@R=1ex{
*{\circ}<3pt> \ar@{-}[dr]^<{} \\
& *{\circ}<3pt> \ar@{-}[r]_<{} 
& {} \ar@{.}[r]&{}  \ar@{-}[r]_>{} &
*{\circ}<3pt> \ar@{=}[r] |-{\scalebox{2}{\object@{>}}}& *{\circ}<3pt>\ar@{}_<{} \\
*{\circ}<3pt> \ar@{-}[ur]_<{}}} & \quad \longleftrightarrow \quad
{\rm CB}^\bullet:\ 
&\vcenter{\xymatrix@R=1ex{
*{\circ}<3pt> \ar@{=}[r] |-{\scalebox{2}{\object@{>}}}_<{} 
&*{\circ}<3pt> \ar@{-}[r]_<{} 
& {} \ar@{.}[r]&{}  \ar@{-}[r]_>{} &
*{\circ}<3pt> \ar@{=}[r] |-{\scalebox{2}{\object@{>}}}
& *{\bullet}<3pt>\ar@{}_<{}
}
}
\end{alignat*}
}
\smallskip
{\small
\begin{alignat*}{2}
{\rm DB}^\bullet:\ 
&\vcenter{\xymatrix@R=1ex{
*{\circ}<3pt> \ar@{-}[dr]^<{} \\
& *{\circ}<3pt> \ar@{-}[r]_<{} 
& {} \ar@{.}[r]&{}  \ar@{-}[r]_>{} &
*{\circ}<3pt> \ar@{=}[r] |-{\scalebox{2}{\object@{>}}}& *{\bullet}<3pt>\ar@{}_<{} \\
*{\circ}<3pt> \ar@{-}[ur]_<{}}} & \quad \longleftrightarrow \quad
{\rm CB}:\ 
&\vcenter{\xymatrix@R=1ex{
*{\circ}<3pt> \ar@{=}[r] |-{\scalebox{2}{\object@{>}}}_<{} 
&*{\circ}<3pt> \ar@{-}[r]_<{} 
& {} \ar@{.}[r]&{}  \ar@{-}[r]_>{} &
*{\circ}<3pt> \ar@{=}[r] |-{\scalebox{2}{\object@{>}}}
& *{\circ}<3pt>\ar@{}_<{}}}
\end{alignat*}
}
\smallskip
{\small
\begin{alignat*}{2}
{\rm BB}:\ 
&\vcenter{\xymatrix@R=1ex{
*{\circ}<3pt> \ar@{=}[r] |-{\scalebox{2}{\object@{<}}}_<{} 
&*{\circ}<3pt> \ar@{-}[r]_<{} 
& {} \ar@{.}[r]&{}  \ar@{-}[r]_>{} &
*{\circ}<3pt> \ar@{=}[r] |-{\scalebox{2}{\object@{>}}}
& *{\circ}<3pt>\ar@{}_<{}}} & \quad \longleftrightarrow \quad
{\rm B}^\bullet\! {\rm B}^\bullet:\ 
&\vcenter{\xymatrix@R=1ex{
*{\bullet}<3pt> \ar@{=}[r] |-{\scalebox{2}{\object@{<}}}_<{} 
&*{\circ}<3pt> \ar@{-}[r]_<{} 
& {} \ar@{.}[r]&{}  \ar@{-}[r]_>{} &
*{\circ}<3pt> \ar@{=}[r] |-{\scalebox{2}{\object@{>}}}
& *{\bullet}<3pt>\ar@{}_<{}}}
\end{alignat*}}
{\small
\begin{alignat*}{2}
{\rm BB}^\bullet:\ 
&\vcenter{\xymatrix@R=1ex{
*{\circ}<3pt> \ar@{=}[r] |-{\scalebox{2}{\object@{<}}}_<{} 
&*{\circ}<3pt> \ar@{-}[r]_<{} 
& {} \ar@{.}[r]&{}  \ar@{-}[r]_>{} &
*{\circ}<3pt> \ar@{=}[r] |-{\scalebox{2}{\object@{>}}}
& *{\bullet}<3pt>\ar@{}_<{}}} & \quad \longleftrightarrow \quad
{\rm B}^\bullet\!{\rm B}:\ 
&\vcenter{\xymatrix@R=1ex{
*{\bullet}<3pt> \ar@{=}[r] |-{\scalebox{2}{\object@{<}}}_<{} 
&*{\circ}<3pt> \ar@{-}[r]_<{} 
& {} \ar@{.}[r]&{}  \ar@{-}[r]_>{} &
*{\circ}<3pt> \ar@{=}[r] |-{\scalebox{2}{\object@{>}}}
& *{\circ}<3pt>\ar@{}_<{}}}\\
\end{alignat*}}
\noindent 
where the first case is studied in this paper (see \cite{Ya99} for the quantum affine superalgebras associated to the Dynkin diagrams above).

For example, we may obtain the same results in this paper using $\U_{CC}(\e)$: Let $\mathring{\U}_{C}(\e)$ and $\mathring{\U}_D(\e)$ be the generalized quantum groups of finite types $C$ and $D$ associated to $\e=(\e_1,\dots,\e_n)$. There is an isomorphism of $\Bbbk$-algebras from $\mathring{\U}_{C}(\e)$ to $\mathring{\U}_D(\e')$ where $\e'=(\dots,\e_n,\e_{n-1})$ \cite[Section 5.3]{Ma}. 
Using this isomorphism and the arguments in Section \ref{subsec:general e}, we obtain two families of irreducible $q$-oscillator representations for $\U_{CC}(\e)$ for $\e={\bm\e}, {\bm\e}'$, and the correspondence between types ${\rm DD}$ and ${\rm CC}$ in the above table.

But the other cases would be more challenging, since the spectral decomposition on the tensor product of any two finite-dimensional fundamental representations is not known, and the fundamental representations are not classically irreducible in general. We also remark that the denominator of a normalized $R$ matrix is enough to understand the tensor structure of the category of finite-dimensional representations, while the denominator is not well-defined in case of $q$-oscillator representations.

\appendix
\section{Proof of Proposition~\ref{prop:reduction}}\label{sec:app-pf-Serre-rel}
In this appendix, we verify the following identity from the proof of Proposition~\ref{prop:reduction}
\[
\widehat{e}^{2}_0 \widehat{e}_1 - (q^2+q^{-2})\widehat{e}_0 \widehat{e}_1 \widehat{e}_0 + \widehat{e}_1 \widehat{e}^{2}_0 =0,
\]
where $\widehat{e}_i$ ($i=0,1$) was given in \eqref{eq:hat{e} for c} as follows:
\[
\widehat{e}_0 = \frac{1}{[2]}(e_0 e_1 + q^2 e_1 e_0), \quad \widehat{e}_1 = e_2 e_3 +q e_3 e_2.
\]
Here we fix $\eta =1$, for the proof being the same for $\eta= -1$.
Recall that $e_j$ ($j=0,1,2,3$) is a generator of $\mathcal{U}_D (\boldsymbol{\epsilon})$ and so satisfies the following relations:
\begin{align}
e_{i}^{2} & =0,\quad e_{0}e_{3}=e_{3}e_{0},\quad e_{1}e_{3}=e_{3}e_{1},\nonumber \\
e_{0}e_{1}e_{2}-e_{1}e_{0}e_{2} & +[2](e_{1}e_{2}e_{0}-e_{0}e_{2}e_{1})+e_{2}e_{0}e_{1}-e_{2}e_{1}e_{0}=0,\label{eq:Serre_D}\\
e_{0}e_{2}e_{3}e_{2}-e_{3}e_{2}e_{0}e_{2} & -[2]e_{2}e_{3}e_{0}e_{2}-e_{2}e_{0}e_{2}e_{3}+e_{2}e_{3}e_{2}e_{0}=0,\nonumber \\
e_{1}e_{2}e_{3}e_{2}-e_{3}e_{2}e_{1}e_{2} & -[2]e_{2}e_{3}e_{1}e_{2}-e_{2}e_{1}e_{2}e_{3}+e_{2}e_{3}e_{2}e_{1}=0.\nonumber 
\end{align}
In particular, by the first relation, 
\[
\widehat{e}_{0}^{2}=\frac{1}{[2]^{2}}(e_{0}e_{1}+q^{2}e_{1}e_{0})(e_{0}e_{1}+q^{2}e_{1}e_{0})=\frac{1}{[2]^{2}}(e_{0}e_{1}e_{0}e_{1}+q^{4}e_{1}e_{0}e_{1}e_{0})
\]
and then the identity we are going to show can be written as
\begin{align}\label{eq:app-A-step-0}
 0= \,& (010123)+q^{4}(101023)+q(010132)+q^{5}(101032) \\
 & -(q^{2}+q^{-2})\begin{pmatrix}(012301)+q^{2}(102301)+q^{2}(012310)+q^{4}(102310) \\
+q(013201)+q^{3}(103201)+q^{3}(013210)+q^{5}(103210)
\end{pmatrix}\nonumber\\
 & +(230101)+q^{4}(231010)+q(320101)+q^{5}(321010),  \nonumber
\end{align}
where we denote by $(i_{1}i_{2}\dots i_{l}) = e_{i_{1}}e_{i_{2}}\cdots e_{i_{l}}$
for brevity. 

The proof consists of repeatedly replacing a monomial in $e_i$ whose coefficient is the lowest power of $q$ with the ones with higher $q$ powers, using the relation \eqref{eq:Serre_D}. Namely, in \eqref{eq:app-A-step-0} we begin from $-q^{-2}(012301)=-q^{-2}(012013)$, which is equal to
\begin{align*}
-q^{-2}(0\overline{120}13) =& -q^{-1}(0)\{(012)-(102)+q(120)-[2](021)+(201)-(210)\}(13)\\
 =& -q^{-1}(01\overline{021}3)-q^{-1}(\overline{021}013)+(012013)\\
 =& -(01)\left((012)-(102)+[2](120)-q(021)+(201)-(210)\right)(3)\\
 & \quad-\left((012)-(102)+[2](120)-q(021)+(201)-(210)\right)(013)+(012013)\\
=-(010123) & +q(010213)+(012103)-(012013)+(102013)+q(021013)-(201013).
\end{align*}
For reader's convenience, we indicate by an overline which factor of a monomial is substituted by \eqref{eq:Serre_D}. We have also used the identities $(00)=(11)=0$ and $(30)=(03),\,(31)=(13)$.
After substitution, the righthand side of \eqref{eq:app-A-step-0} becomes
\begin{align}\label{eq:app-A-step-1}
 & q^{4}(101023)+q(010132)+q^{5}(101032)\nonumber \\
 & -(q^{2}+1)(012301)-q^{4}(102301)-q^{4}(012310)-(q^{6}+q^{2})(102310)\nonumber \\
 & -(q^{3}+q^{-1})(013201)-(q^{5}+q)(103201)-(q^{5}+q)(013210)-(q^{7}+q^{3})(103210)\\
 & +q^{4}(231010)+q(320101)+q^{5}(321010)\nonumber \\
 & +q(010231)+q(023101).\nonumber 
\end{align}

Next, since
\begin{align*}
-q^{-1}(30\overline{120}1) =& -(30)\left((012)-(102)+q(120)-[2](021)+(201)-(210)\right)(1)\\
 =& -(301\overline{021})-(3\overline{021}01)+q(301201)\\
 =& -q(301)\left((012)-(102)+[2](120)-q(021)+(201)-(210)\right)\\
 & \quad-q(3)\left((012)-(102)+[2](120)-q(021)+(201)-(210)\right)(01)+q(301201)\\
=-q(301012) & +q^{2}(301021)-q(301201)+q(301210)+q(310201)+q^{2}(302101)-q(320101),
\end{align*}
\eqref{eq:app-A-step-1} is equal to
\begin{align}\label{eq:app-A-step-2}
 & q^{4}(101023)+q^{5}(101032)\nonumber \\
 & -(q^{2}+1)(012301)-q^{4}(102301)-q^{4}(012310)-(q^{6}+q^{2})(102310)\nonumber \\
 & -(q^{3}+q)(013201)-q^{5}(103201)-q^{5}(013210)-(q^{7}+q^{3})(103210)\\
 & +q^{4}(231010)+q^{5}(321010)\nonumber \\
 & +q(010231)+q(023101)+q^{2}(010321)+q^{2}(032101).\nonumber 
\end{align}
Here
\begin{align*}
-(q^{2}+1) & (012301)=-q[2](0\overline{120}13)\\
= & q(0)\left((012)-(102)-[2](021)+(201)-(210)\right)(13)\\
= & -q(010213)-q(021013)
\end{align*}
and similarly
\[
-(q^{3}+q)  (013201)=-q^{2}[2](30\overline{120}1) = -q^{2}(301021)-q^{2}(302101),
\]
so that \eqref{eq:app-A-step-2} reduces to
\begin{align}\label{eq:app-A-step-3}
 & q^{4}(101023)+q^{5}(101032)\nonumber \\
 & -q^{4}(102301)-q^{4}(012310)-(q^{6}+q^{2})(102310)\nonumber \\
 & -q^{5}(103201)-q^{5}(013210)-(q^{7}+q^{3})(103210)\\
 & +q^{4}(231010)+q^{5}(321010).\nonumber 
\end{align}
Repeating the above procedure, we conclude that \eqref{eq:app-A-step-3} is equal to zero.

\section{Proof of Lemma~\ref{lem:EF action on u^ij}}\label{sec:app-pf-lem-EF}
In this appendix, we prove Lemma~\ref{lem:EF action on u^ij}. 
We shall check only the last two identities, as the others are easier.
We will further assume that $x_1=x_2=1$, from which a general case follows by easy tracking of spectral parameters.

For brevity, we introduce the following notation:
\[
w^{a,b,c,d}_{x,y,z,w} \coloneqq \ket{a \mathbf{e}_1 + b \mathbf{e}_3 + c \mathbf{e}_{2m-1} + d \mathbf{e}_{2m+1}}\otimes \ket{x \mathbf{e}_1 + y \mathbf{e}_3 + z \mathbf{e}_{2m-1} + w \mathbf{e}_{2m+1}} \in \mathcal{W}_{\underline{\boldsymbol{\e}}'}^{\otimes 2}
\]
where we are identifying $\mathcal{W}_{\underline{\boldsymbol{\e}}'}^{\otimes 2}=\mathfrak{tr}^{\boldsymbol{\e}'}_{\underline{\boldsymbol{\e}}'}\mathcal{W}_{\boldsymbol{\e}'}^{\otimes 2}
\subset \mathcal{W}_{\boldsymbol{\e}'}^{\otimes 2}$.

\subsection{The formula for $\boldsymbol{F}_{m+1}u_{r,s}^{i,j}$}
We begin from
\[
f_{m}^{(j)}v_{l_1}=\begin{bmatrix}l_1\\
j
\end{bmatrix}w^{0,0,j,l_1 -j}_{0,0,0,0}.
\]
Since $e_0,\,e_1,\,\dots,e_m$ commute with $f_{m+1}$, we have
\begin{align*}
  e_{0}u_{r,s}^{i,j}= & \begin{bmatrix}l_{1}\\
  j
  \end{bmatrix}\left(f_{m+1}^{(i)}e_{0}w^{0,0,j,l_{1}-j}_{0,0,0,0}\right)\otimes(k_{0}^{-1}f_{m+1}^{(r-i)}f_{m}^{(s-j)}v_{l_{2}})\\
   & +\begin{bmatrix}l_{2}\\
  s-j
  \end{bmatrix}(f_{m+1}^{(i)}f_{m}^{(j)}v_{l_{1}})\otimes\left(f_{m+1}^{(r-i)}e_{0}w^{0,0,s-j,l_{2}-s+j}_{0,0,0,0}\right)\\
  = & \begin{bmatrix}l_{1}\\
  j
  \end{bmatrix}q^{-2}f_{m+1}^{(i)}\left(w^{1,0,j,l_{1}-j}_{0,1,0,0}-q^{-1}w^{0,1,j,l_{1}-j}_{1,0,0,0}\right)\otimes(f_{m+1}^{(r-i)}f_{m}^{(s-j)}v_{l_{2}})\\
   & +\begin{bmatrix}l_{2}\\
  s-j
  \end{bmatrix}(f_{m+1}^{(i)}f_{m}^{(j)}v_{l_{1}})\otimes f_{m+1}^{(r-i)}\left(w^{1,0,s-j,l_{2}-s+j}_{0,1,0,0}-q^{-1}w^{0,1,s-j,l_{2}-s+j}_{1,0,0,0}\right),
  \end{align*}
  and so
  \begin{align*}
  (e_{m-1} & \cdots e_{3}e_{2})e_{0}u_{r,s}^{i,j}\\
  = & \begin{bmatrix}l_{1}\\
  j
  \end{bmatrix}q^{-2}f_{m+1}^{(i)}e_{m-1}\cdots e_{2} \left( w^{1,0,j,l_{1}-j}_{0,1,0,0}-q^{-1}w^{0,1,j,l_{1}-j}_{1,0,0,0} \right) \otimes k_{m-1}^{-1}\cdots k_{2}^{-1}f_{m+1}^{(r-i)}f_{m}^{(s-j)}v_{l_{2}}\\
   & +\begin{bmatrix}l_{2}\\
  s-j
  \end{bmatrix}(f_{m+1}^{(i)}f_{m}^{(j)}v_{l_{1}})\otimes f_{m+1}^{(r-i)}e_{m-1}\cdots e_{2}
  \left( w^{1,0,s-j,l_{2}-s+j}_{0,1,0,0}-q^{-1} w^{0,1,s-j,l_{2}-s+j}_{1,0,0,0} \right) \\
  = & \begin{bmatrix}l_{1}\\
  j
  \end{bmatrix}q^{-2-s+j-r+i}f_{m+1}^{(i)}
  \left( w^{1,0,j,l_{1}-j}_{0,0,1,0}-q^{-1}w^{0,0,j+1,l_{1}-j}_{1,0,0,0} \right)
  \otimes f_{m+1}^{(r-i)}f_{m}^{(s-j)}v_{l_{2}}\\
   & +\begin{bmatrix}l_{2}\\
  s-j
  \end{bmatrix}(f_{m+1}^{(i)}f_{m}^{(j)}v_{l_{1}})\otimes f_{m+1}^{(r-i)} 
  \left( w^{1,0,s-j,l_{2}-s+j}_{0,0,1,0}-q^{-1}w^{0,0,s-j+1,l_{2}-s+j}_{1,0,0,0} \right).
  \end{align*}
Next, we apply $e_{m}$ to $(e_{m-1}\cdots e_3 e_2 )e_0 u_{r,s}^{i,j}$:
  \begin{align*}
  e_{m}( & e_{m-1}\cdots e_{3}e_{2})e_{0}u_{r,s}^{i,j}\\
  = & \begin{bmatrix}l_{1}\\
  j
  \end{bmatrix}q^{-2+s-j-r+i-l_{2}}f_{m+1}^{(i)}\begin{Bmatrix}q[j]w^{1,0,j-1,l_{1}-j+1}_{0,0,1,0}\\
  +w^{1,0,j,l_{1}-j}_{0,0,0,1}\\
  -q^{-1}[j+1]w^{0,0,j,l_{1}-j+1}_{1,0,0,0}
  \end{Bmatrix}\otimes f_{m+1}^{(r-i)}f_{m}^{(s-j)}v_{l_{2}}\\
   & +[l_{2}-s+j+1]\begin{bmatrix}l_{1}\\
  j
  \end{bmatrix}q^{-2-s+j-r+i}f_{m+1}^{(i)}\begin{Bmatrix}w^{1,0,j,l_{1}-j}_{0,0,1,0}\\
  -q^{-1}w^{0,0,j+1,l_{1}-j}_{1,0,0,0}
  \end{Bmatrix}\otimes f_{m+1}^{(r-i)}f_{m}^{(s-j-1)}v_{l_{2}}\\
   & +\begin{bmatrix}l_{2}\\
  s-j
  \end{bmatrix}[l_{1}-j+1]q^{2s-2j+1-2l_{2}}f_{m+1}^{(i)}f_{m}^{(j-1)}v_{l_{1}}\otimes f_{m+1}^{(r-i)}\begin{Bmatrix}w^{1,0,s-j,l_{2}-s+j}_{0,0,1,0}\\
  -q^{-1}w^{0,0,s-j+1,l_{2}-s+j}_{1,0,0,0}
  \end{Bmatrix}\\
   & +\begin{bmatrix}l_{2}\\
  s-j
  \end{bmatrix}f_{m+1}^{(i)}f_{m}^{(j)}v_{l_{1}}\otimes f_{m+1}^{(r-i)}\begin{Bmatrix}q[s-j]w^{1,0,s-j-1,l_{2}-s+j+1}_{0,0,1,0}\\
  +w^{1,0,s-j,l_{2}-s+j}_{0,0,0,1}\\
  -q^{-1}[s-j+1]w^{0,0,s-j,l_{2}-s+j+1}_{1,0,0,0}
  \end{Bmatrix},
  \end{align*}
  where we use the well-known formula
  \begin{equation}\label{comm-rel-e,f}
  e_{i}f_{i}^{(k)}v=[a-k+1]f_{i}^{(k-1)}v+f_{i}^{(k)}e_{i}v
  \end{equation}
  for any $v$ such that $k_{i}v=q^{a}v$. To complete the computation,
  note that the second summand in the right hand side above is annihilated by $e_{m-1}\cdots e_{2}e_{1}$:
  \[
  e_{m-1}\cdots  e_{1}\left(w^{1,0,j,l_{1}-j}_{0,0,1,0}-q^{-1}w^{0,0,j+1,l_{1}-j}_{1,0,0,0}\right)
   =q^{-1}w^{0,0,j+1,l_{1}-j}_{0,0,1,0}-q^{-1}w^{0,0,j+1,l_{1}-j}_{0,0,1,0}=0
  \]
  and so is the third one. Finally, we have
  \begin{align*}
  \boldsymbol{F}_{m+1} & u_{r,s}^{i,j}=(e_{m-1}\cdots  e_{2}e_{1})e_{m}(e_{m-1}\cdots e_{3}e_{2})e_{0}u_{r,s}^{i,j}\\
  = & \begin{bmatrix}l_{1}\\
  j
  \end{bmatrix}q^{-2-2r+2i-l_{2}}f_{m+1}^{(i)}\begin{Bmatrix}[j]w^{0,0,j,l_{1}-j+1}_{0,0,1,0}\\
  +w^{0,0,j+1,l_{1}-j}_{0,0,0,1}\\
  -q^{-1}[j+1]w^{0,0,j,l_{1}-j+1}_{0,0,1,0}
  \end{Bmatrix}\otimes f_{m+1}^{(r-i)}f_{m}^{(s-j)}v_{l_{2}}\\
   & +\begin{bmatrix}l_{2}\\
  s-j
  \end{bmatrix}f_{m+1}^{(i)}f_{m}^{(j)}v_{l_{1}}\otimes f_{m+1}^{(r-i)}\begin{Bmatrix}[s-j]w^{0,0,s-j,l_{2}-s+j+1}_{0,0,1,0}\\
  +w^{0,0,s-j+1,l_{2}-s+j}_{0,0,0,1}\\
  -q^{-1}[s-j+1]w^{0,0,s-j,l_{2}-s+j+1}_{0,0,1,0}
  \end{Bmatrix}\\
  = & \begin{bmatrix}l_{1}\\
  j
  \end{bmatrix}q^{-2-2r+2i-l_{2}}f_{m+1}^{(i)}
  \left( -q^{-j-1}w^{0,0,j,l_{1}-j+1}_{0,0,1,0}+w^{0,0,j+1,l_{1}-j}_{0,0,0,1} \right)
  \otimes f_{m+1}^{(r-i)}f_{m}^{(s-j)}v_{l_{2}} \\
   & +\begin{bmatrix}l_{2}\\
  s-j
  \end{bmatrix}f_{m+1}^{(i)}f_{m}^{(j)}v_{l_{1}}\otimes f_{m+1}^{(r-i)}
  \left( -q^{-s+j-1}w^{0,0,s-j,l_{2}-s+j+1}_{0,0,1,0}+w^{0,0,s-j+1,l_{2}-s+j}_{0,0,0,1} \right)\\
  = & q^{2i-2r-l_{2}-2}[i+1]u_{r+1,s}^{i+1,j}+[r-i+1]u_{r+1,s}^{i,j}
  \end{align*}
  as desired.

\subsection{The formula for $\boldsymbol{E}_{m+1}u_{r,s}^{i,j}$}
We first expand $f_{m+1}^{(i)}f_{m}^{(j)}v_{l_1}$. Recall that
\[
f_{m}^{(j)}v_{l_1}=\begin{bmatrix}l_1\\
j
\end{bmatrix}w^{0,0,j,l_1-j}_{0,0,0,0},
\]
and let us write
\[
f_{m+1}^{(i)}f_{m}^{(j)}v_{l_1}=\begin{bmatrix}l_1\\ j \end{bmatrix} f_{m+1}^{(i)}w^{0,0,j,l_1-j}_{0,0,0,0}=\frac{1}{[i]!}\begin{bmatrix}l_1\\j \end{bmatrix}\sum_{x=0}^{i}d_{i,x} w^{0,0,i-x+j,x+l_1 -j}_{0,0,x,i-x}
\]
for some $d_{i,x}\in\Bbbk$. Similarly we can write
\[
f_{m+1}^{(r-i)}f_{m}^{(s-j)}v_{l_2}=\frac{1}{[r-i]!}\begin{bmatrix}l_2\\s-j \end{bmatrix}\sum_{y=0}^{r-i}d'_{r-i,y}w^{0,0,r-i-y+s-j,y+l_2 -s+j}_{0,0,y,r-i-y}.
\]
Then we have
\begin{align*}
f_{1}\cdots &  f_{m-1}u_{r,s}^{i,j}\\
=&(f_{1}\cdots f_{m-1}f_{m+1}^{(i)}f_{m}^{(j)}v_{l_1})\otimes f_{m+1}^{(r-i)}f_{m}^{(s-j)}v_{l_2}\\
&+(k_{1}\cdots k_{m-1}f_{m+1}^{(i)}f_{m}^{(j)}v_{l_1})\otimes(f_{1}\cdots f_{m-1}f_{m+1}^{(r-i)}f_{m}^{(s-j)}v_{l_2})\\
 =& \frac{1}{[i]!}\begin{bmatrix}l_1\\
j
\end{bmatrix}\sum_{x=0}^{i}\begin{pmatrix}[i-x+j]d_{i,x}w^{1,0,i-x+j-1,x+l_1 -j}_{0,0,x,i-x}\\
+q^{i-x+j}[x]d_{i,x}w^{0,0,i-x+j,x+l_1-j}_{1,0,x-1,i-x}
\end{pmatrix}\otimes f_{m+1}^{(r-i)}f_{m}^{(s-j)}v_{l_2}\\
 & +\frac{q^{i+j}}{[r-i]!}\begin{bmatrix}l_2\\
s-j
\end{bmatrix}f_{m+1}^{(i)}f_{m}^{(j)}v_{l_1}\\
& \quad\quad\quad\quad\quad\quad\quad\otimes
\sum_{y=0}^{r-i}\begin{pmatrix}[r-i-y+s-j]d^{\prime}_{r-i,y}w^{1,0,r-i-y+s-j-1,y+l_2-s+j}_{0,0,y,r-i-y}\\
+q^{r-i-y+s-j}[y]d^{\prime}_{r-i,y}w^{0,0,r-i-y+s-j,y+l_2-s+j}_{1,0,y-1,r-i-y}
\end{pmatrix}.
\end{align*}
Then we apply $f_{m}$:
\begin{align*}
f_{m} & (f_{1}\cdots f_{m-1})u_{r,s}^{i,j}=S_1 + S_2 + S_3 + S_4,\\
S_1 & =\frac{1}{[i]!}\begin{bmatrix}l_1\\
j
\end{bmatrix}\sum_{x=0}^{i}\begin{Bmatrix}[i-x+j]d_{i,x}\begin{pmatrix}[x+l_1-j]w^{1,0,i-x+j,x+l_1 -j-1}_{0,0,x,i-x}\\
+q^{l_1 +2x-i-2j+1}[i-x]w^{1,0,i-x+j-1,x+l_1 -j}_{0,0,x+1,i-x-1}
\end{pmatrix}\\
+q^{i-x+j}[x]d_{i,x}\begin{pmatrix}[x+l_1 -j]w^{0,0,i-x+j+1,x+l_1 -j-1}_{1,0,x-1,i-x}\\
+q^{l_1 +2x-i-2j}[i-x]w^{0,0,i-x+j,x+l_1 -j}_{1,0,x,i-x-1}
\end{pmatrix}
\end{Bmatrix} \\
& \hspace{10cm}\otimes f_{m+1}^{(r-i)}f_{m}^{(s-j)}v_{l_2},\\
 S_2 & = \frac{1}{[i]!}\begin{bmatrix}l_1\\
j
\end{bmatrix}[s-j+1]q^{l_1 -2j+1}\sum_{x=0}^{i}\begin{pmatrix}[i-x+j]d_{i,x}w^{1,0,i-x+j-1,x+l_1 -j}_{0,0,x,i-x}\\
+q^{i-x+j}[x]d_{i,x}w^{0,0,i-x+j,x+l_1 -j}_{1,0,x-1,i-x}
\end{pmatrix} \\
&\hspace{10cm} \otimes f_{m+1}^{(r-i)}f_{m}^{(s-j+1)}v_{l_2},\\
 S_3 & = \frac{q^{i+j}}{[r-i]!}\begin{bmatrix}l_2\\
s-j
\end{bmatrix}[j+1]f_{m+1}^{(i)}f_{m}^{(j+1)}v_{l_1}\\
& \qquad\qquad\qquad\qquad\otimes \sum_{y=0}^{r-i}\begin{pmatrix}[r-i-y+s-j]d^{\prime}_{r-i,y}w^{1,0,r-i-y+s-j-1,y+l_2 -s+j}_{0,0,y,r-i-y}\\
+q^{r-i-y+s-j}[y]d^{\prime}_{r-i,y}w^{0,0,r-i-y+s-j,y+l_2 -s+j}_{1,0,y-1,r-i-y}
\end{pmatrix},\\
 S_4 & =  \frac{q^{l_1 +i-j}}{[r-i]!}\begin{bmatrix}l_2\\
s-j
\end{bmatrix}f_{m+1}^{(i)}f_{m}^{(j)}v_{l_1}\\
\otimes & \sum_{y=0}^{r-i}\begin{Bmatrix}[r-i-y+s-j]d^{\prime}_{r-i,y}\begin{pmatrix}[y+l_2 -s+j]w^{1,0,r-i-y+s-j,y+l_2 -s+j-1}_{0,0,y,r-i-y}\\
+q^{l_2 +2y-r+i-2s+2j+1}[r-i-y]w^{1,0,r-i-y+s-j-1,y+l_2 -s+j}_{0,0,y+1,r-i-y-1}
\end{pmatrix}\\
+q^{r-i-y+s-j}[y]d^{\prime}_{r-i,y}\begin{pmatrix}[y+l_2 -s+j]w^{0,0,r-i-y+s-j+1,y+l_2 -s+j-1}_{1,0,y-1,r-i-y}\\
+q^{l_2 +2y-r+i-2s+2j}[r-i-y]w^{0,0,r-i-y+s-j,y+l_2 -s+j}_{1,0,y,r-i-y-1}
\end{pmatrix}
\end{Bmatrix}.
\end{align*}
Next, consider $f_0 (f_{2}\cdots f_{m-1})f_{m}(f_{1}\cdots f_{m-1})u_{r,s}^{i,j}$, which we still divide into four summands as above.
Then we have $f_0 (f_2 \cdots f_{m-1})S_2 = f_0 (f_2 \cdots f_{m-1})S_3 = 0$.
For example,
\begin{align}\label{eqn:vanishing-appendix}
f_0 & (f_2  \cdots f_{m-1})S_2 \\
=& \left\{f_0 (f_2 \cdots f_{m-1})  \left([i-x+j]d_{i,x}w^{1,0,i-x+j-1,x+l_1 -j}_{0,0,x,i-x} 
+ q^{i-x+j}[x]d_{i,x}w^{0,0,i-x+j,x+l_1 -j}_{1,0,x-1, i-x} \right)\right\} \nonumber \\
& \hspace{10.5cm}\otimes f_{m+1}^{(r-i)}f_{m}^{(s-j+1)}v_{l_2} \nonumber\\
=& f_0 
\begin{Bmatrix} 
[i-x+j]d_{i,x}\left( [i-x+j-1]w^{1,1,i-x+j-2,x+l_1 -j}_{0,0,x,i-x} 
+ q^{i-x+j-1}[x]w^{1,0,i-x+j-1,x+l_1 -j}_{0,1,x-1,i-x} \right)
\\ 
+q^{i-x+j}[x]d_{i,x} \left( [i-x+j]w^{0,1,i-x+j-1,x+l_1 -j}_{1,0,x-1,i-x} 
+ q^{i-x+j}[x-1]w^{0,0,i-x+j,x+l_1 -j}_{1,1,x-2,i-x} \right)
\end{Bmatrix} \nonumber\\
& \hspace{10.5cm} \otimes f_{m+1}^{(r-i)}f_{m}^{(s-j+1)}v_{l_2} \nonumber\\
=& \left( d_{i,x}[i-x+j][x]q^{i-x+j-1}\cdot(-q)+d_{i,x}[x][i-x+j]q^{i-x+j} \right)  \nonumber\\
& \hspace{6cm} \cdot w^{0,0,i-x+j-1,x+l_1 -j}_{0,0,x-1,i-x}\otimes f_{m+1}^{(r-i)}f_{m}^{(s-j+1)}v_{l_2} \nonumber\\
= & 0. \nonumber
\end{align}
Let us look into the remaining two summands. For $f_0 (f_2 \cdots f_{m-1})S_1 $, we first compute the sum of the first and the
third rows:
\begin{align*}
f_{0}&(f_{2}\cdots f_{m-1})\frac{1}{[i]!}\begin{bmatrix}l_1\\
j
\end{bmatrix}\sum_{x=0}^{i}f_{0}\begin{Bmatrix}[i-x+j]d_{i,x}[x+l_1 -j]w^{1,0,i-x+j,x+l_1 -j-1}_{0,0,x,i-x}\\
+q^{i-x+j}[x]d_{i,x}[x+l_1 -j]w^{0,0,i-x+j+1,x+l_1 -j-1}_{1,0,x-1,i-x}
\end{Bmatrix} \\
& \hspace{11cm}\otimes f_{m+1}^{(r-i)}f_{m}^{(s-j)}v_{l_2}\\
= & \frac{1}{[i]!}\begin{bmatrix}l_1\\
j
\end{bmatrix}\sum_{x=0}^{i}d_{i,x}[x+l_1 -j]f_{0}\begin{Bmatrix}[i-x+j][x]q^{i-x+j}w^{1,0,i-x+j,x+l_1 -j-1}_{0,1,x-1,i-x}\\
+q^{i-x+j}[x][i-x+j+1]w^{0,1,i-x+j,x+l_1 -j-1}_{1,0,x-1,i-x}
\end{Bmatrix} \\
& \hspace{11cm} \otimes f_{m+1}^{(r-i)}f_{m}^{(s-j)}v_{l_2}\\
= & \frac{1}{[i]!}\begin{bmatrix}l_1\\
j
\end{bmatrix}\sum_{x=0}^{i}d_{i,x}[x+l_1 -j][x]q^{-i-x+j}(-[i-x+j]q+[i-x+j+1]) \\
&\hspace{7cm} \cdot w^{0,0,i-x+j,x+l_1 -j-1}_{0,0,x-1,i-x}\otimes f_{m+1}^{(r-i)}f_{m}^{(s-j)}v_{l_2}\\
= & \frac{1}{[i]!}\begin{bmatrix}l_1\\
j
\end{bmatrix}\sum_{x=0}^{i}d_{i,x}[x+l_1 -j][x]w^{0,0,i-x+j,x+l_1 -j-1}_{0,0,x-1,i-x}\otimes f_{m+1}^{(r-i)}f_{m}^{(s-j)}v_{l_2}.
\end{align*}
The sum of the second and the fourth rows is
\begin{align*}
f_{0}&(f_{2} \cdots f_{m-1})\frac{1}{[i]!}\begin{bmatrix}l_1\\
j
\end{bmatrix}\sum_{x=0}^{i}f_{0}\begin{Bmatrix}[i-x+j]d_{i,x}q^{l_1 +2x-i-2j+1}[i-x]w^{1,0,i-x+j-1,x+l_1 -j}_{0,0,x+1,i-x-1}\\
+q^{i-x+j}[x]d_{i,x}q^{l_1 +2x-i-2j}[i-x]w^{0,0,i-x+j,x+l_1 -j}_{1,0,x,i-x-1}
\end{Bmatrix} \\
&\hspace{11cm} \otimes f_{m+1}^{(r-i)}f_{m}^{(s-j)}v_{l_2}\\
= & \frac{1}{[i]!}\begin{bmatrix}l_1\\
j
\end{bmatrix}\sum_{x=0}^{i}d_{i,x}[i-x]q^{l_1 +2x-i-2j}f_{0}\begin{Bmatrix}[i-x+j]q^{i-x+j}[x+1]w^{1,0,i-x+j-1,x+l_1 -j}_{0,1,x,i-x-1}\\
+q^{i-x+j}[x][i-x+j]\\w^{0,1,i-x+j-1,x+l_1 -j}_{1,0,x,i-x-1}
\end{Bmatrix} \\
&\hspace{11cm} \otimes f_{m+1}^{(r-i)}f_{m}^{(s-j)}v_{l_2}\\
= & \frac{1}{[i]!}\begin{bmatrix}l_1\\
j
\end{bmatrix}\sum_{x=0}^{i}d_{i,x}[i-x][i-x+j]q^{l_1 +x-j}(-q[x+1]+[x])\\
&\hspace{6cm} \cdot w^{0,0,i-x+j-1,x+l_1 -j}_{0,0,x,i-x-1} \otimes f_{m+1}^{(r-i)}f_{m}^{(s-j)}v_{l_2}\\
= & -\frac{1}{[i]!}\begin{bmatrix}l_1\\
j
\end{bmatrix}\sum_{x=0}^{i}d_{i,x}[i-x][i-x+j]q^{l_1 +2x-j+1} w^{0,0,i-x+j-1,x+l_1 -j}_{0,0,x,i-x-1}\otimes f_{m+1}^{(r-i)}f_{m}^{(s-j)}v_{l_2}.
\end{align*}
To sum up, 
\begin{align*}
f_{0}(f_{2}\cdots & f_{m-1})S_1 \\
=&\frac{1}{[i]!}\begin{bmatrix}l_1\\
j
\end{bmatrix}\sum_{x=0}^{i} \left(d_{i,x}[x+l_1 -j][x]-d_{i,x-1}[i-x+1][i-x+j+1]q^{l_1 +2x-j-1}\right)\\
 & \hspace{6cm} \cdot w^{0,0,i-x+j,x+l_1 -j-1}_{0,0,x-1,i-x}\otimes f_{m+1}^{(r-i)}f_{m}^{(s-j)}v_{l_2}.
\end{align*}
On the other hand, one can easily check that this right hand side coincides with
\begin{align*}
(e_{m+1}f_{m+1}^{(i)}f_{m}^{(j)} & v_{l_1})\otimes f_{m+1}^{(r-i)}f_{m}^{(s-j)}v_{l_2}\\
&=\left( e_{m+1}\frac{1}{[i]!}\begin{bmatrix}l_1\\
j
\end{bmatrix}\sum_{x=0}^{i}d_{i,x}w^{0,0,i-x+j,x+l_1 -j}_{0,0,x,i-x}\right)
\otimes f^{(r-i)}_{m+1}f^{(s-j)}_{m}v_{l_2}.
\end{align*}
Hence we arrive at
\begin{align*}
f_{0}(f_{2}\cdots f_{m-1})S_1 & = (e_{m+1}f_{m+1}^{(i)}f_{m}^{(j)}v_{l_1})\otimes f_{m+1}^{(r-i)}f_{m}^{(s-j)}v_{l_2}\\
 & =-[l_1 +i+1]f_{m+1}^{(i-1)}f_{m}^{(j)}v_{l_1}\otimes f_{m+1}^{(r-i)}f_{m}^{(s-j)}v_{l_2}=-[l_1 +i+1]u_{r-1,s}^{i-1,j}.
\end{align*}
By a similar computation, we get
\[
f_{0}(f_{2}\cdots f_{m-1})S_4 =-q^{l_1 +2i+2}[l_2 +r-i+1]u_{r-1,s}^{i,j}
\]
and complete the verification.

\section{Proof of Lemma~\ref{lem:F-action-on-hwvectors}}\label{sec:app-pf-E-hw}
We set
\[
A^{i}_{r}=(-1)^i \prod_{k=1}^{i} q^{2k-2-l_2 -2r}\frac{[r+l_2 -k+2]}{[l_1+k+1]},\quad
B^{j}_{s}=(-q)^j \prod_{k'=1}^{j} q^{2k'+l_2 -2s}\frac{[l_2-s+k']}{[l_1-k'+1]},
\]
so that $u_{r,s}=\sum_{i,j} A^{i}_{r} B^{j}_{s} u_{r,s}^{i,j}$ (see Lemma~\ref{lem:maximal vectors of osc D}).

Fix $r\geq 0$ and $0\leq s \leq \min(l_1,l_2)$. We shall prove only that
\[
\boldsymbol{F}_{m+1}u_{r,s} \in \Bbbk^\times u_{r+1,s}+\Bbbk f_{m+1}u_{r,s}+\Bbbk f^{(2)}_{m+1}u_{r-1,s},
\]
and then the other assertion of the Lemma~\ref{lem:F-action-on-hwvectors} is easier.

From the classical decomposition of $\mathcal{W}_{l_1}(x_1)\otimes\mathcal{W}_{l_2}(x_2)$ \eqref{eq:decom of two tensor of fundamentals'} and weight consideration, we have
\[
\boldsymbol{F}_{m+1}u_{r,s} = \sum_{a,b\geq 0} C_{a,b}f^{(a)}_{m+1}f^{(b)}_{m} u_{r+1-a,s-b}
\]
for some $C_{a,b}\in\Bbbk$. We first show that $C_{a,b}=0$ whenever $b>0$, by verifying $e_m \boldsymbol{F}_{m+1} u_{r,s}=0$. The latter one is straightforward: as
\[
e_m u_{r,s}^{i,j} = q^{-l_2 +2s-2j}[l_1 -j+1]u_{r,s-1}^{i,j-1}+[l_2 -s+j+1]u_{r,s-1}^{i,j},
\]
we have
\begin{align*}
e_{m}&\boldsymbol{F}_{m+1}u_{r,s}  =e_{m}\sum_{i,j}A^{i}_{r}B^{j}_{s}\boldsymbol{F}_{m+1}u_{r,s}^{i,j}\\
=& \sum A^{i}_{r}B^{j}_{s} \left(x_1 q^{2i-2r-l_2 -2}[i+1]e_{m}u_{r+1,s}^{i+1,j} +x_2 [r-i+1]e_{m}u_{r+1,s}^{i,j}\right)\\
=&\sum A^{i}_{r}B^{j}_{s} \begin{Bmatrix}x_1 q^{-l_2 -2r+2i-2}[i+1]\left(q^{-l_2 +2s-2j}[l_1 -j+1] u_{r+1,s-1}^{i+1,j-1} +[l_2 -s+j+1]u_{r+1,s-1}^{i+1,j} \right)\\
 +x_2[r-i+1]\left(q^{-l_2 +2s-2j}[l_1 -j+1]u_{r+1,s-1}^{i,j-1} +[l_2 -s+j+1]u_{r+1,s-1}^{i,j}\right)
\end{Bmatrix}\\
=&\sum\begin{Bmatrix}
x_1 A_{r}^{i-1}B_{s}^{j+1} q^{-2l_2 -2r+2i+2s-2j-6}[i][l_1 -j]+x_1 A_{r}^{i-1}B_{s}^{j} q^{-l_2 -2r+2i-4}[i][l_2 -s+j+1]\\
+x_2 A_{r}^{i}B_{s}^{j+1} q^{-l_2 +2s-2j-2}[r-i+1][l_1 -j]+x_2 A_{r}^{i}B_{s}^{j} [r-i+1][l_2 -s+j+1]
\end{Bmatrix}\\
&\hspace{12cm}\cdot u_{r+1,s-1}^{i,j}.
\end{align*}
Then each coefficient of $u_{r+1,s-1}^{i,j}$ can be easily seen to vanish, using
\[
\frac{A_{r}^{i-1}}{A_{r}^{i}}=-q^{l_2 +2r-2i+2}\frac{[l_1+i+1]}{[r+l_2 -i+2]},\quad
\frac{B_{s}^{j+1}}{B_{s}^{j}}=-q^{l_2 -2s+2j+2}\frac{[l_2 -s+j+1]}{[l_1-j]}.
\]

Next, we claim $e^{2}_{m+1}\boldsymbol{F}_{m+1}u_{r,s}\in \Bbbk^{\times}u_{r-1,s}$, which implies $C_{a,0}=0$ whenever $a>2$.
Indeed, from
\begin{align*}
e_{m+1}^{2} &u_{r,s}^{i,j}=-q^{l_2 +2r-2i+2}[l_1 +i+1] e_{m+1} u^{i-1,j}_{r-1,s} - [l_2 +r-i+1] e_{m+1} u^{i,j}_{r-1,s}\\
&= q^{2l_2 +4r-4i+4}[l_1+i+1][l_1+i]u_{r-2,s}^{i-2,j} +q^{l_2 +2r-2i+2}[l_1+i+1][l_2 +r-i+1]u_{r-2,s}^{i-1,j}\\
&\quad\quad +q^{l_2 +2r-2i}[l_2 +r-i+1][l_1 +i+1]u_{r-2,s}^{i-1,j}+[l_2 +r-i+1][l_2+r-i]u_{r-2,s}^{i,j}\\
&=q^{2l_2 +4r-4i+4}[l_1+i+1][l_1+i]u_{r-2,s}^{i-2,j} + q^{l_2 +2r-2i+1}[2][l_1+i+1][l_2 +r-i+1]u_{r-2,s}^{i-1,j} \\
&\quad\quad +[l_2 +r-i+1][l_2+r-i]u_{r-2,s}^{i,j},
\end{align*}
we have
\begin{align*}
&e^{2}_{m+1}\boldsymbol{F}_{m+1}u_{r,s}  \\
=&e^{2}_{m+1}\sum_{i,j}A^{i}_{r}B^{j}_{s}\boldsymbol{F}_{m+1}u_{r,s}^{i,j}\\
=& \sum A^{i}_{r}B^{j}_{s} \left(x_1 q^{2i-2r-l_2 -2}[i+1]e^{2}_{m+1}u_{r+1,s}^{i+1,j} +x_2 [r-i+1]e^{2}_{m+1}u_{r+1,s}^{i,j}\right)\\
=& \sum A^{i}_{r}B^{j}_{s} \begin{Bmatrix}x_1 q^{-l_2 -2r+2i-2}[i+1]\begin{pmatrix}
q^{2l_2 +4r-4i+4}[l_1 +i+2][l_1 +i+1]u_{r-1,s}^{i-1,j}\\
+q^{l_2 +2r-2i+1}[2][l_1 +i+2][l_2 +r-i+1]u_{r-1,s}^{i,j}\\
+[l_2 +r-i+1][l_2 +r-i]u_{r-1,s}^{i+1,j}
\end{pmatrix}\\
+x_2 [r-i+1]\begin{pmatrix}
q^{2l_2 +4r-4i+8}[l_1 +i+1][l_1 +i]u_{r-1,s}^{i-2,j}\\
+q^{l_2 +2r-2i+3}[2][l_1 +i+1][l_2 +r-i+2]u_{r-1,s}^{i-1,j}\\
+[l_2 +r-i+2][l_2 +r-i+1]u_{r-1,s}^{i,j}
\end{pmatrix}\end{Bmatrix}\\
=&\sum B_{s}^{j}\begin{Bmatrix}
x_2 A_{r}^{i+2}q^{2l_2 +4r-4i}[r-i-1][l_2 +i+3][l_2 +i+2]\\
+A_{r}^{i+1}q^{l_2 +2r-2i}[l_1 +i+2]\left(x_1 [i+2][l_1 +i+3]+x_2 q[2][r-i][l_2 +r-i+1]\right)\\
+A_{r}^{i}[l_2 +r-i+1]\left(x_1 q^{-1}[2][i+1][l_1 +i+2] + x_2 [r-i+1][l_2 +r-i+2]\right)\\
+x_1 A_{r}^{i-1}q^{-l_2 -2r+2i-4}[i][l_2 +r-i+2][l_2 +r-i+1]
\end{Bmatrix}u_{r-1,s}^{i,j}\\
=&\sum B_{s}^{j}A_{r}^{i}[l_2 +r-i+1]\begin{Bmatrix}
x_2 q^{2}[r-i-1][l_2 +r-i]\\
-x_1 [i+2][l_1 +i+3] -x_2 q[2][r-i][l_2 +r-i+1]\\
+x_1 q^{-1}[2][i+1][l_1 +i+2] +x_2 [r-i+1][l_2 +r-i+2]\\
-x_1 q^{-2}[i][l_1 +i+1]
\end{Bmatrix}u_{r-1,s}^{i,j}\\
=& \sum B_{s}^{j}A_{r}^{i}[l_2 +r-i+1][2]\left(x_2 q^{-l_2 -2r+2i}-x_1 q^{l_1 +2i+2} \right)u_{r-1,s}^{i,j}\\
=& \sum B_{s}^{j}A_{r-1}^{i}[l_2 +r+1][2]\left(x_2 q^{-l_2 -2r} -x_1 q^{l_1 +2}\right)u_{r-1,s}^{i,j}\\
=&[l_2 +r+1][2]\left(x_2 q^{-l_2 -2r} -x_1 q^{l_1 +2}\right)u_{r-1,s}.
\end{align*}
In the sixth equality, we used the following identity
\[
q^{2}[a][b]-q[2][a+1][b+1]+[a+2][b+2]=q^{-a-b-1}[2]
\]
for $a,\,b\in\mathbb{Z}$. 

Finally, it remains to see $C_{0,0}\neq 0$. Note that since
\begin{align*}
e^{2}_{m+1}\boldsymbol{F}_{m+1}u_{r,s}&=e^{2}_{m+1}\left(C_{0,0}u_{r+1,s}+C_{1,0}f_{m+1}u_{r,s}+C_{2,0}f^{(2)}_{m+1}u_{r-1,s}\right) \\
&= C_{2,0}[l_1 +l_2 +2r][l_1 +l_2 +2r+1]u_{r-1,s},
\end{align*}
the above computation also determines
\[
C_{2,0}=\frac{[l_2 +r+1][2]}{[l_1 +l_2 +2r][l_1 +l_2 +2r+1]}\left(x_2 q^{-l_2 -2r} -x_1 q^{l_1 +2}\right).
\]
Comparing the coefficients of $u_{r,s}^{0,0}$ of both sides of the following identities
\begin{align*}
-C_{1,0} & [l_1 +l_2 +2r+2]u_{r,s}-C_{2,0}[l_1 +l_2 +2r+1]f_{m+1}u_{r-1,s}\\
&=e_{m+1}\left(C_{0,0}u_{r+1,s}+C_{1,0}f_{m+1}u_{r,s}+C_{2,0}f^{(2)}_{m+1}u_{r-1,s}\right) \\
&=e_{m+1}\boldsymbol{F}_{m+1}u_{r,s}\\
&= \sum A^{i}_{r}B^{j}_{s} \left(q^{2i-2r-l_2 -2}[i+1]e_{m+1}u_{r+1,s}^{i+1,j} +[r-i+1]e_{m+1}u_{r+1,s}^{i,j}\right),
\end{align*}
we obtain
\begin{align*}
  C_{1,0}&=\frac{1}{[l_{1}+l_{2}+2r+2]} ([l_{1}+2]+[r+1][l_{2}+r+2]-q^{2}[r][l_{2}+r+1] \\
  &\qquad\qquad\qquad\qquad\qquad\qquad\qquad\qquad -\frac{(x_2 q^{-l_{1}-l_{2}-2r-2}-x_1)[2][l_{2}+r+1][r]}{[l_{1}+l_{2}+2r]}).
\end{align*}
Similarly, $C_{0,0}$ can be computed from 
\[
\boldsymbol{F}_{m+1}u_{r,s}=C_{0,0}u_{r+1,s}+C_{1,0}f_{m+1}u_{r,s}+C_{2,0}f^{(2)}_{m+1}u_{r-1,s}
\]
in which the coefficients of $u_{r+1,s}^{0,0}$ are
\[
  [r+1]=C_{0,0}+q^{-l_{1}-2}[r+1]C_{1,0}+q^{-2l_{1}-4}\frac{[r][r+1]}{[2]}C_{2,0}.
\]
Substituting $C_{1,0}$ and $C_{2,0}$ above in this identity, one can easily check that $C_{0,0}$ is nonzero for generic $x_1$ and $x_2$.

{\small

\end{document}